

\nopagenumbers
\parindent= 15pt

\hsize=14cm
\vsize=21cm
\hoffset=0.9cm
\voffset=1cm

\input amssym.def
\input amssym.tex

\font\ssrm=cmr6
\font\srm=cmr8

\font\csc=cmcsc10
\font\title=cmr12 at 14pt

\font\teneusm=eusm10    
\font\seveneusm=eusm7  
\font\fiveeusm=eusm5    
\newfam\eusmfam
\def\eusm{\fam\eusmfam\teneusm}
\textfont\eusmfam=\teneusm 
\scriptfont\eusmfam=\seveneusm
\scriptscriptfont\eusmfam=\fiveeusm

\def\Re{{\rm Re}\,}
\def\Im{{\rm Im}\,}
\def\sgn{{\rm sgn}\,}
\def\txt#1{{\textstyle{#1}}}
\def\scr#1{{\scriptstyle{#1}}}
\def\c#1{{\cal #1}}
\def\f#1{{\goth{#1}}}
\def\r#1{{\rm #1}}
\def\e#1{{\eusm #1}}
\def\varGamma{{\mit\Gamma}}
\def\B#1{{\Bbb #1}}
\def\b#1{{\bf #1}}
\def\sgn{{\rm sgn}}

\def\rightheadline{\hfil{\srm Elements of automorphic 
representations}\hfil\tenrm\folio}
\def\leftheadline{\tenrm\folio\hfil{\srm Y. Motohashi}\hfil}
\def\emptyheadline{}
\headline{\ifnum\pageno=1 \emptyheadline\else
\ifodd\pageno \rightheadline \else \leftheadline\fi\fi}

\def\firstpage{\hss{\vbox to 1.5cm
{\vfil\hbox{\rm\folio}}}\hss}
\def\emptyfootline{\hfil}
\footline{\ifnum\pageno=1\firstpage\else
\emptyfootline\fi}

\centerline{\title Elements of Automorphic Representations}
\vskip 0.5cm
\centerline{\csc By Yoichi Motohashi}
\vskip 0.7cm
This is an attempt at 
a practical and essentially self-contained 
theory of automorphic representations 
in the framework
$$
\hbox{
$L^2(\varGamma\backslash\r{G})$ with
$\r{G}=\r{PSL}(2,\B{R})$ and 
$\varGamma=\r{PSL}(2,\B{Z})$.}
$$
The restriction of the underlying
discrete subgroup $\varGamma$ is imposed
solely for the sake of simplicity; our argument should
extend to fairly general arithmetic subgroups of $\r{G}$ 
without substantial alteration.
Our motivation lies in the observation that
applying the spectral theory 
of cusp forms on $\r{G}$ to problems in
analytic number theory we need a
lowbrow but in fact highly 
informative approach, 
as then it is of paramount importance
to be able to work in regions of 
absolute and uniform convergence.
We need to have explicit descriptions, including
the convergence issues, of integral 
transforms arising in representation theory; 
and in order to acquire various
means indispensable in dealing with relevant involved
technicalities, it is best to
trace  from scratch how those transforms come into play. 
Essential in applications are not only
the spectral structures but also the procedure to establish
them. Therefore we shall take an approach based on
explicit computation rather than the common dose of the theory
of commuting compact operators which we are, nevertheless,
aware is susceptible of generalisations to bigger Lie groups.
\par
With this, our
attention is specifically directed to an accessible explication of
the Kirillov map and the Bessel functions of representation.
Both concepts have played 
fundamental r\^oles in recent applications of
representation theory of Lie groups
to problems in analytic number theory
such as the spectral theory of sums of Kloosterman sums and
the theory of mean values of the Riemann zeta 
and a variety of $L$-functions. It is worth stressing that
the massive cancelations 
which take place in these subjects have
been detected only with analysis of
integral transforms indicated above; indeed such a detection
appears to be beyond the reach of the sole use of 
algebra--operator theoretic means.
\par
We shall naturally use terms from the
theory of Lie groups and algebras; but there is
no need to know the entire theory in order to
understand this article. For those who 
share interest and taste with us, the first three introductory chapters of
Vilenkin--Klimyk [36] are recommended to have
an overview of the theory, although it is of no
absolute necessity either.
What we are developing below may in fact be regarded as the
material that should be basic
in order to enter into the theory of Lie groups and
their representations in much the same sense as 
elementary number theory is meant for analytic number theory.
However, to restrict the text within a reasonable size,
we assume a familiarity with
the spectral theory of real analytic cusp forms
on the upper half-plane, i.e., the weight zero situation, a fairly
elementary account of which is
developed in [22, Chapter 1]. 
Salient points of
the situation with even integral weights
are given in later sections of the present article. 
Also, as an appendix,
Selberg's trace formula and zeta-function 
associated with the group $\varGamma$
are treated, although the subject is
not directly related to our principal aim; a novelty is mainly in a 
relatively fast approach to the quintessence of Selberg's theory.
Notations and conventions are introduced 
along with necessity, and will continue to 
be effective thereafter. References are limited to those
immediately related to our purpose. 
\smallskip
\noindent
{\csc Acknowledgements.} This article is partly based
on the author's talk delivered at the centennial celebration of the birth
of his late mentor Tur\'an P\'al 
(August 21--25, 2011, Budapest);  
he is indebted to the organisers of the conference
for the support and the hospitality. The author wishes to record that
in course of learning the elements of the
theory of representations of Lie groups
during the past years, especially while preparing
the present work, he benefited greatly from numerous
correspondences with Roelof W. Bruggeman of Utrecht. 
The author profited as well from expert comments 
by Nigel Watt of Edinburgh on preliminary versions.
\medskip
\noindent 
{\bf 1. Hyperbolic plane.} The group $\r{SL}(2,\B{R})$ acts on 
$z\in\B{C}\cup\{\infty\}$: 
$$
\r{h}(z)={az+b\over cz+d}\,,\quad
\left(\matrix{a&b\cr c& d}\right)\in
\r{SL}(2,\B{R}).\eqno(1.1)
$$
We have, with $\jmath(\r{h},z)=cz+d$,
$$
\displaystyle{{d\over dz}\r{h}(z)
={1\over\jmath(\r{h},z)^2},\quad
\Im\r{h}(z)={y\over|\jmath(\r{h}, z)|^2},}\atop
\displaystyle{\r{h}(z_1)-\r{h}(z_2)
={z_1-z_2\over\jmath(\r{h},z_1)
\jmath(\r{h},z_2)},\quad \jmath(\r{h}_1\r{h}_2, z)
=\jmath(\r{h}_1,\r{h}_2(z))\jmath(\r{h}_2,z).}\eqno(1.2)
$$
In particular, $\r{SL}(2,\B{R})$ is
a group of orientation-preserving 
conformal transformations of the upper
half plane 
$$
\B{H}^2=\big\{z=x+iy: x\in\B{R}, y>0\big\}.\eqno(1.3)
$$ 
This is, however, the same as dealing with the group
$\r{G}=\r{SL}(2,\B{R})/\{\pm1\}$, and we shall
use the notation
$$
\r{h}=\left[\matrix{a&b\cr c& d}\right]\in\r{G}\eqno(1.4)
$$
in place of $(1.1)$; see the notes below.
On $\B{H}^2$ we have the hyperbolic metric
$$
{|dz|\over y}={\big((dx)^2+(dy)^2\big)^{1/2}\over y}
={\big((du)^2+(dv)^2\big)^{1/2}\over v}
={|dw|\over v},\eqno(1.5)
$$
with $\r{h}(z)=w,\, z=x+iy,\, w=u+iv$; the invariance 
follows from the first line of $(1.2)$. 
This induces the hyperbolic measure
$$
d\mu(z)={dxdy\over y^2}={dudv\over v^2}=d\mu(w).
\eqno(1.6)
$$
Also induced are the hyperbolic outer-normal differential
and the Laplace--Beltrami operator: We have,
for any smooth function $f$ on $\B{H}^2$, 
$$
\eqalignno{
y\left({dy\over |dz|}{\partial\over\partial x}
-{dx\over |dz|}{\partial\over\partial y}\right)f(w)&=
v\left({dv\over |dw|}{\partial\over\partial u}
-{du\over |dw|}{\partial\over\partial v}\right)f(w),&(1.7)\cr
y^2\left(\left({\partial\over\partial x}\right)^2
+\left({\partial\over\partial y}\right)^2\right)f(w)&=
v^2\left(\left({\partial\over\partial u}\right)^2
+\left({\partial\over\partial v}\right)^2\right)f(w),&(1.8)
}
$$
where the differentiation on the left sides are performed on the
function $f(w)=f(\r{h}(z))$ of the variable $z$.
The invariance $(1.6)$--$(1.8)$ is a consequence of 
the Cauchy--Riemann equation for the function $\r{h}(z)$. They are
applied together with Green's formula in an 
obvious hyperbolic disguise: We have,
for any smooth functions $f,\,g$ on $\B{H}^2$,
$$
\eqalignno{
\int_D y^2\left(\left({\partial\over\partial x}\right)^2
+\left({\partial\over\partial y}\right)^2\right)& f\cdot g\,d\mu(z)
+\int_D y^2\left({\partial f\over\partial x}{\partial g\over\partial x}
+{\partial f\over\partial y}
{\partial g\over\partial y}\right)d\mu(z),\cr
=\int_{\partial D}&y\left({dy\over |dz|}{\partial\over\partial x}
-{dx\over |dz|}{\partial\over\partial y}\right)f\cdot g\,
{|dz|\over y},&(1.9)
}
$$
where $D$ is a domain in $\B{H}^2$ enclosed by a 
piece-wise smooth boundary curve $\partial D$ which is
positively oriented.
\par
Further, we have the fact that the action 
of the subgroup $\varGamma$ is discrete
and induces the tessellation
$$
\B{H}^2=\bigcup_{\gamma\in\varGamma}\gamma\e{F},
\quad\hbox{disjoint save for boundaries 
$\gamma\partial\e{F}$,}
\eqno(1.10)
$$
with
$$
\e{F}=\left\{z\in\B{H}^2: |\Re z|\le{1\over2},\,|z|\ge1
\right\}.\eqno(1.11)
$$
A function $f$ is said to be $\varGamma$-automorphic
if $f(\gamma(z))=f(z)$, $\forall(z,\gamma)\in
\B{H}^2\times\varGamma$, that is, $f$ is a function on
the Riemann surface $\varGamma\backslash\B{H}^2$.
We then introduce the Hilbert space
$$
L^2(\varGamma\backslash\B{H}^2)
=\left\{{\hbox{$\varGamma$-automorphic 
and }\atop\hbox{square integrable
over $\e{F}$ against $d\mu$}}\right\},\eqno(1.12)
$$
\par
The relations $(1.5)$--$(1.9)$ are basic
implements in developing the harmonic analysis
on $\varGamma\backslash\B{H}^2$, i.e.,
a spectral resolution of the Laplace--Beltrami operator 
which is symmetric on $L^2(\varGamma
\backslash\B{H}^2)$ as the invariance
structure of $(1.9)$ implies.
We shall extend them
to the Lie group $\r{G}$ that spreads
over $\B{H}^2$: Anticipating concepts
to be introduced in due course, we assert
that $(1.5)$ corresponds to 
point-pair invariants, $(1.6)$ to a Haar measure,
$(1.7)$ to the Maass operators, 
$(1.8)$ to the Casimir operator; and $(1.9)$ is
an instance of the effect of the decomposition of the 
Casimir operator in terms of the Maass operators,
although the details of these correspondences
are given only in Section 32, as they are irrelevant to our
main task to develop elements of automorphic
representations. The invariance
of these notions on $\r{G}$
are in the core of the differentiable structure of
$\r{G}$, and plays a central r\^ole in 
the development of the harmonic analysis
on $\varGamma\backslash\r{G}$, i.e., a 
spectral resolution of the Casimir operator that 
is symmetric on the Hilbert space 
$L^2(\varGamma\backslash\r{G})$ spreading over
$L^2(\varGamma\backslash\B{H}^2)$. 
\medskip
\noindent
{\csc Notes:} The notation $(1.4)$ causes a minor confusion:
The function $\jmath(\r{h},z)$ can equal 
either $cz+d$ or $-(cz+d)$.
Nevertheless, within the present article this
ambiguity should not cause any trouble, 
since the relevant instances are all consequences
of the basic relations $(3.3)$--$(3.4)$ below. 
See the notes to the third section.
\medskip
\noindent
{\bf 2. Hyperbolic distance.} The
distance $d(z_1,z_2)$ between $z_1$ and $z_2$ on $\B{H}^2$
is defined to be the minimum of lengths
of curves connecting the
points, measured against the hyperbolic line element $|dz|/y$.
It holds that $d\big(\r{h}(z_1),\r{h}(z_2)\big)
=d(z_1,z_2)$ for any $\r{h}\in\r{G}$, as is implied by 
the invariance $(1.5)$. We have
$$
\displaystyle{d(z_1,z_2)=2\,\r{arcsinh}\sqrt{\varrho(z_1,z_2)},
\quad \varrho(z_1,z_2)={|z_1-z_2|^2\over4(\Im z_1)
(\Im z_2)},}\atop
\displaystyle{\varrho\big(\r{h}(z_1),\r{h}(z_2)\big)
=\varrho(z_1,z_2),\quad\forall\r{h}\in\r{G},}\eqno(2.1)
$$
in which the second line follows from the second and the third
identities in $(1.2)$.
In fact, mapping $i$ to $z_2=x_2+iy_2$, we have
$$
d(z_1,z_2)=d\big(\r{h}^{-1}_0(z_1),i\big)
=d\big((z_1-x_2)/y_2,i\big),\quad \r{h}_0=
\left[\matrix{\sqrt{y_2}& x_2/\sqrt{y_2}
\cr&1/\sqrt{y_2}}\right]\in\r{G}.\eqno(2.2)
$$
Then, via
$$
z\mapsto w=\r{c}(z), 
\quad\r{c}=\left(\matrix{1&-i\cr1&\hfill i}\right),\eqno(2.3)
$$
we map $\B{H}^2$ onto the unit disk $|w|<1$, and apply
a rotation $w\mapsto \exp(i\tau)w$ so that the point
$(z_1-x_2)/y_2$ is mapped to 
$|z_1-z_2|/|z_1-\overline{z}_2|$. 
With the inverse of $(2.3)$ we return to $\B{H}^2$, finding
$$
d(z_1,z_2)=d\big(r(z_1,z_2)i, i\big)=\log(r(z_1,z_2)),
\eqno(2.4)
$$
where
$$
r(z_1,z_2)={|z_1-z_2|+|z_1-\overline{z}_2|\over
|z_1-\overline{z}_2|-|z_1-z_2|}
=\big(\varrho(z_1,z_2)^{1/2}
+(\varrho(z_1,z_2)+1)^{1/2}\big)^2.\eqno(2.5)
$$
We get $(2.1)$.
\medskip
\noindent
{\csc Notes:} Elements of the hyperbolic geometry is
 given in, e.g.,  Maass [20]. An interesting historical account
can be found in Penrose [28, Section 2.6].
\medskip
\noindent
{\bf 3. Iwasawa decomposition.} Ascension to $\r{G}$ starts.
Classifying 
results of the
action of $\r{G}$ on the point $i$, we obtain a 
co-ordinate system on $\r{G}$:
$$
\r{G}=\r{NAK}\ni\r{g}
=\r{n}[x]\r{a}[y]\r{k}[\theta],\eqno(3.1)
$$
with
$$
\displaystyle{\r{N}=\left\{\r{n}[x]=\left[\matrix{1&x\cr&1}\right]:
\,x\in{\Bbb R}\right\},\quad
\r{A}=\left\{\r{a}[y]
=\left[\matrix{\sqrt{y}&\cr&1/\sqrt{y}}\right]:
\,y>0\right\},}\atop\displaystyle{
\r{K}=\left\{\r{k}[\theta]=
\left[\matrix{\phantom{-}\cos\theta&
\sin\theta\cr-\sin\theta&\cos\theta}\right]:\,\theta\in{\Bbb
R}/\pi{\Bbb Z}\right\},\atop\displaystyle{
\r{n}[x]\r{a}[y]\r{k}[\theta]=\left[\matrix{\sqrt{y}\cos\theta-
x\sin\theta/\sqrt{y}&\sqrt{y}\sin\theta+x\cos\theta/\sqrt{y}\cr
-\sin\theta/\sqrt{y}&\cos\theta/\sqrt{y}}\right].}
}\eqno(3.2)
$$
The relations
$$
\r{g}=\r{n}[x]\r{a[y]}\r{k}[\theta]=
\left[\matrix{a&b\cr c&d}\right],\atop\displaystyle
{x={ac+bd\over c^2+d^2},\quad
 y=\big(c^2+d^2\big)^{-1},\quad
\exp(2i\theta)={\jmath(\r{g},-i)
\over \jmath(\r{g},i)},}\eqno(3.3)
$$
imply that $(3.1)$ with $(3.2)$ is indeed 
a co-ordinate system on $\r{G}$. 
We have, for $\r{h}\in\r{G}$,
$$
{\r{h}\cdot\r{n}[x]\r{a}[y]\r{k}[\theta]=\r{n}[x_1]\r{a}[y_1]
\r{k}[\theta_1],}\atop\displaystyle{\r{h}(x+iy)=x_1+iy_1,
\quad \exp(2i\theta_1)=\exp(2i\theta)
{\jmath(\r{h},x-iy)\over \jmath(\r{h},x+iy)},}\eqno(3.4)
$$
in which the second identity is the result of applying the first to
the point $i$ and the third follows from the last identity in $(1.2)$
and that in $(3.3)$. Thus the action of 
$\r{G}$ on $\B{H}^2$ is equivalent to that on the set 
of all right $\r{K}$-cosets of $\r{G}$; 
that is,
$$
\B{H}^2\cong \r{G}/\r{K},\eqno(3.5)
$$
with respect to the action of $\r{G}$. We note that 
under our formulation $\r{G}$ acts
on $\r{G}/\r{K}$ through the left multiplication or translation 
$$
l_\r{h}:\r{g}\mapsto\r{h}\r{g}.\eqno(3.6)
$$
\medskip
\noindent
{\csc Notes:} Hereafter it is understood that the variables 
$(x,y,\theta)$ are the co-ordinates $(3.1)$ unless otherwise stated.
Alternatively, one may ascend from
$\B{H}^2$ to $\r{SL}(2,\B{R})$. It does not entail the ambiguity
mentioned in the notes to the first section. On the other hand,
we would encounter side effects as well, specifically
in dealing with assertions involving the notion of
weights to be introduced in Section 12. Our choice of $\r{G}$ means
that we restrict our discussion to even functions 
on $\r{SL}(2,\B{R})$, i.e., sums of
functions of even integral weights on $\r{G}$; 
this suffices for basic applications
of automorphic representations to problems in analytic
number theory.
\medskip
\noindent
{\bf 4. Cartan decomposition.} This notion is 
closely related to point-pair invariants on $\r{G}$,
that is, any function $ F(\r{g}_1,\r{g}_2)$ on
$\r{G}\times\r{G}$ such that $F(\r{h}\r{g}_1,\r{h}\r{g}_2)
=F(\r{g}_1,\r{g}_2)$, $\forall\r{h}\in\r{G}$, which is 
in fact a function of $\r{g}_1^{-1}\r{g}_2$. We shall see 
that such a function is a natural extension of the hyperbolic distance.
Thus, along with $(3.1)$, we have
the decomposition
$$
\hbox{$\r{G}=\r{KAK}\ni\r{g}
=\r{k}[\tau_1]\r{a}[u]\r{k}[\tau_2]$,
uniquely provided $u>1$. }\eqno(4.1)
$$
We note that
$$
\r{g}^{-1}=\r{k}[-\tau_2]\r{w}\r{a}
[u]\r{w}^{-1}\r{k}[-\tau_1],
\quad \hbox{$\r{w}=\r{k}\big[\txt{1\over2}\pi\big]$
: the Weyl element of $\r{G}$.}\eqno(4.2)
$$
To show $(4.1)$, we consider $\r{g}\cdot\r{g}^t$ which 
is symmetric and positive. We have
$(\r{g}\cdot\r{g}^t)^{-1/2}\r{g}\in\r{K}$, and
get the decomposition. 
As to the uniqueness, we assume that
$\r{k}[\tau_1]\r{a}[u_1]\r{k}[\tau_2]=\r{a}[u_2]$.
Applying it to the point $i$, we have
$$
{(1-u_1^2)\sin(\tau_1)
\cos(\tau_1)+u_1i\over u_1^2\sin^2(\tau_1)
+\cos^2(\tau_1)}=u_2i,\eqno(4.3)
$$
which implies $\sin(2\tau_1)=0$ as we assume $u_1>1$.
If $\tau_1\equiv0\bmod\pi\B{Z}$, then $u_1=u_2$, and
$\tau_2\equiv0\bmod\pi\B{Z}$. On the other hand, if
$\tau_1\equiv{1\over2}\pi\bmod\pi\B{Z}$, then $u_1^{-1}
=u_2$, which is excluded.
\par
For $\r{g}=\r{n}[x]\r{a}[y]\r{k}[\theta]$ 
as in $(4.1)$, we have $u=\exp(d(z,i))$, $z=x+iy$.
In fact, we have $d(z,i)=d(\r{g}(i),i)=
d(\r{k}[\tau_1]\r{a}[u]\r{k}[\tau_2](i),i)
=d(\r{a}[u](i),i)=\log u$ as we may assume that
$u\ge1$.
Also, considering $\r{c}(z)=\r{c}\r{g}\r{c}^{-1}(0)$,
we get 
$$
{z-i\over z+i}=e^{2i\tau_1}{|z-i|\over|z+i|}.\eqno(4.4)
$$
More generally, on noting $(1.2)$, $(3.3)$ and
$d(\r{g}_1^{-1}\r{g}_2(i),i)
=d(\r{g}_1(i),\r{g}_2(i))$, we derive from
$(4.4)$ that 
$$
{\hbox{$\r{g}_1^{-1}\r{g}_2=\r{k}
[\eta_1]\r{a}[v]\r{k}[\eta_2]$
with $v=\exp(d(z_1,z_2))$ and}}\atop
\displaystyle{{z_2-z_1\over z_2-\overline{z}_1}
=e^{2i(\theta_1+\eta_1)}
{|z_2-z_1|\over |z_2-\overline{z}_1|},\quad
{z_1-z_2\over z_1-\overline{z}_2}
=-e^{2i(\theta_2-\eta_2)}
{|z_1-z_2|\over |z_1-\overline{z}_2|},}\eqno(4.5)
$$
where $\r{g}_j=\r{n}[x_j]\r{a}[y_j]\r{k}[\theta_j]$,
$z_j=x_j+iy_j$. To get the second identity, we write
$z_2=\r{g}_1(\r{g}_1^{-1}\r{g_2}(i))$, $z_1=\r{g}_1(i)$,
and apply the fourth expression of $(1.2)$, the
last of $(3.3)$ as well as $(4.4)$.
The third identity in $(4.5)$ follows from the
second,  since we have $(4.2)$ and
$\r{g}_2^{-1}\r{g}_1=(\r{g}_1^{-1}\r{g}_2)^{-1}$.
\medskip
\noindent
{\csc Notes:} The mode of decomposition $(4.1)$ 
is the same as introducing a polar co-ordinate system on
$\r{G}$; see Bruggeman [4, Section 2.2.6], for instance.
This is included here because 
its natural extension plays a salient r\^ole
in our discussion of automorphic representations
of $\r{PSL}(2,\B{C})$ which is under preparation. See also
Section 33.
\medskip
\noindent
{\bf 5. Invariant measure on $\r{G}$.} Skipping the discussion
of the notion of Haar measures in general, 
we put, in an a priori manner,
$$
\hbox{$\displaystyle{d\r{g}={dxdyd\theta\over\pi y^2}},\,$
with Lebesgue measures $\, dx, dy, d\theta$.}\eqno(5.1)
$$
The group $\r{G}$ is unimodular in the sense that
it admits a left and right invariant Haar measure; that is, we have
$$
d\r{g}=d\r{hg},\quad d\r{g}=d\r{gh}, 
\quad\forall\r{h}\in\r{G}.\eqno(5.2)
$$
The invariance of $d\r{g}$ against the left translation 
is a consequence of $(1.6)$ and $(3.4)$.
On the other hand, if we put
$\r{k}[\theta]\r{h}=\r{n}[\xi(\theta)]\r{a}[u(\theta)]
\r{k}[\vartheta(\theta)]$, then $\r{gh}=
\r{n}[x_1]\r{a}[y_1]
\r{k}[\theta_1]$,  with $x_1=x+\xi(\theta)y$, 
$y_1=u(\theta)y$, $\theta_1=\vartheta(\theta)$.
The Jacobian of the right translation 
$$
r_\r{h}:\r{g}\mapsto\r{gh}\eqno(5.3)
$$
equals $u(\theta)\vartheta'(\theta)$. Applying the
last two identities in $(3.3)$ to $\r{k}[\theta]\r{h}$, 
$\r{h}=\r{n}[\alpha]\r{a}[\beta]\r{k}[\tau]$, we have
$$
u(\theta)={\beta\over
(\cos\theta-\alpha\sin\theta)^2+(\beta\sin\theta)^2},\quad 
e^{2i\vartheta(\theta)}=e^{2i\tau}
{\cos\theta-\alpha\sin\theta+i\beta\sin\theta\over
\cos\theta-\alpha\sin\theta-i\beta\sin\theta},\eqno(5.4)
$$
and find that $\vartheta'(\theta)=u(\theta)$. Thus the
Jacobian is in fact equal to $u^2(\theta)=(y_1/y)^2$, 
which proves the second identity in $(5.2)$.
\medskip
\noindent
{\bf 6. Invariant differential operators.} We are 
to introduce a differentiable structure on $\r{G}$. To this end,
we observe that because of the relations $(3.4)$--$(3.5)$
the harmonic analysis on $\r{G}$ should be an extension of
that on $\B{H}^2$.  As suggested in the first section,
the latter is, by general potential theory, to be
based on the fact $(1.8)$ that the Laplace--Beltrami operator
commutes with the action of $\r{G}$ over $\B{H}^2$; and
this action corresponds to the left translation applied
on $\r{G}/\r{K}$ as already remarked after $(3.5)$. 
Namely, the differentiable structure 
under question should be invariant against the left translation;
in other words it should be defined by an
implement that is independent of any of 
$l_\r{h}$, $\r{h}\in\r{G}$. 
At the same time, as is well indicated by the definition
of the differentiation on the
additive Lie group $\B{R}$, the device ought to be
constructed via infinitesimal action of $\r{G}$ on itself.
In order to fulfil these prerequisites we are naturally
led to exploiting right translations $r_\r{h}$.
\par
With this, we introduce
$$
\b{X}_1=
\left(\matrix{&1\cr\phantom{-1}&}\right),\quad
\b{X}_2=\left(\matrix{1&\cr&-1}\right),\quad\b{X}_3
=\left(\matrix{&1\cr-1&}\right),\eqno(6.1)
$$
and observe that
$$
{\r{N}=\big\{\exp(t\b{X}_1): t\in\B{R}\big\},\quad
\r{A}=\big\{\exp(t\b{X}_2): t\in\B{R}\big\},}
\atop{\r{K}=\big\{\exp(t\b{X}_3):
t\in\B{R}/\pi\B{Z}\big\}}.\eqno(6.2)
$$
The matrices $\exp(t\b{X}_j)$ can obviously 
be identified as the corresponding elements $(1.4)$ of $\r{G}$.
In view of $(3.1)$ these three one-parameter subgroups 
or rather curves on $\r{G}$ give rise to 
the Iwasawa co-ordinate system. Hence
we introduce the right Lie differentials
$$
\b{x}_jf(\r{g})=\Big[{d\over dt}\Big]_{t=0}
f\big(\r{g}\cdot\exp(t\b{X}_j)\big),
\quad f\in C^\infty(\r{G}),\eqno(6.3)
$$
or the right differentiation 
at $\r{g}\in\r{G}$ in the direction $\b{X}_j$; here
$C^\infty(\r{G})$ is the set of all functions which are
differentiable infinitely many times with respect to 
$(x,y,\theta)$.
As $\b{x}_j$ is defined in terms of the right
translation $r_{\exp(t\b{X}_j)}$,
we have obviously
$$
l_\r{h}\b{x}_j=\b{x}_j l_\r{h},\quad 
\forall\r{h}\in\r{G}.\eqno(6.4)
$$
Decomposing $\r{g}\cdot\exp(t\b{X}_j)$ by $(3.3)$,
one may compute $\b{x}_j$ 
in terms of $(x,y,\theta)$: 
$$
\eqalign{
{\bf x}_1&=y\cos(2\theta)\,{\partial\over\partial x}
+y\sin(2\theta)\,{\partial\over\partial y}
+\sin^2\!\theta\,{\partial\over\partial\theta},\cr
{\bf x}_2&=-2y\sin(2\theta)\,{\partial\over\partial x}
+2y\cos(2\theta)\,{\partial\over\partial y}
+\sin(2\theta)\,{\partial\over\partial\theta},\quad
{\bf x}_3={\partial\over\partial\theta}.
}\eqno(6.5)
$$
As an orientation, we indicate how to compute
$\b{x}_1$; the operator $\b{x}_2$ is treated
similarly, and $\b{x}_3$ does not need any explanation:
Essential is to have the first order approximation
$\r{n}[x+\lambda_1 t]\r{a}[y+\lambda_2 t]
\r{k}[\theta+\lambda_3t]$ for $\r{n}[x]\r{a}[y]
\r{k}[\theta]\cdot\r{n}[t]$, that is, 
the part containing $t^2$ or higher powers can be ignored.
By the last expression in $(3.3)$ or rather by the
second formula in $(5.4)$
for $\alpha=t$, $\beta=1$, $\tau=0$, we get readily
$\lambda_3=\sin^2\!\theta$ which gives the third term
of $\b{x}_1$. Similarly, by the first formula in $(5.4)$ we
get $\lambda_2=y\sin(2\theta)$, which corresponds to the
second term of $\b{x}_1$. To compute $\lambda_1$, we
put $\r{g}=\left[{a\atop c}{b\atop d}\right]$ and
thus $\r{g}\cdot\r{n}[t]=\left[{a\atop c}{b+at\atop d+ct}\right]$.
By the second identity in $(3.3)$ we get $\lambda_1$ in
terms $a,b,c, d$ which is in turn expressed 
in terms of $x,y,\theta$ by the last expression in $(3.2)$. 
This ends the computation of $\b{x}_1$. 
\par
The set $\big\{\b{x}_1,\b{x}_2,\b{x}_3\big\}$
generates, over $\B{R}$, 
the Lie algebra $\f{g}$ of $\r{G}$ under the
operation $[\b{x}_i,\b{x}_j]=\b{x}_i\cdot\b{x}_j
-\b{x}_j\cdot\b{x}_i$. The Jacobi identity holds
obviously, and we have
$$
[{\bf x}_1,{\bf x}_2]=-2{\bf x}_1,\quad 
[{\bf x}_1,{\bf x}_3]=-{\bf x}_2,
\quad[{\bf x}_2,{\bf x}_3]=4{\bf x}_1-2{\bf x}_3,
\eqno(6.6)
$$
as is implied by $(6.5)$; in fact it suffices to
compute only the coefficients of the three first order differentials.
Further, in terms of ordinary
operator addition and multiplication
the same set generates the universal
enveloping algebra $\e{U}$ consisting of
left invariant differential operators on $\r{G}$.
It should be noted that the basic field of
$\e{U}$ is naturally $\B{C}$. Namely,
the Lie algebra $\f{g}$, which is originally defined
ove $\B{R}$, is complexified. This remark will become
relevant to the definition of the Maass operators given
in the next section. Thus, hereafter it is always understood that
$\f{g}$ is a Lie algebra over $\B{C}$ generated by
the operators $\big\{\b{x}_1,\b{x}_2,\b{x}_3\big\}$.
\par
Here is a trivial remark: 
Let $\b{x}\in\f{g}$, and let
$\b{x}f_k$ be continuous for all $k$. Then
we have
$$
\b{x}\sum_k f_k(\r{g})=\sum_k\b{x}f_k(\r{g}),
\eqno(6.7)
$$
provided, for instance, both sums converge uniformly.
This can be confirmed on noting that $\b{x}$ is in fact
a differential with respect to a single real variable.
The exchange may hold with a given $\b{u}\in\e{U}$ in place of
$\b{x}\in\f{g}$ as well, if the relevant 
chain of applications of $(6.7)$
can be performed. This procedure will be used
without mentioning details.
\medskip
\noindent
{\csc Notes:}  The most essential in the discussion of
Lie groups is the concept of one parameter subgroups,
a fundamental discovery made by Sophus Lie in 1888.
See Hawkins [10, Section 3.2]; this monograph 
gives an account of an early history 
of the theory of Lie groups. Our definition $(6.3)$ 
is to be regarded as an adaptation of
the general notion of Lie differentials to linear
Lie groups which are composed of matrices. Since
$\det\exp(A)=\exp(\hbox{trace of $A$})$ 
for any square matrix $A$, the set
$\{\exp(t\b{X}): t\in\B{R}\}$ is a curve on $\r{G}$ if and only if
the trace of a $2\times2$ matrix $\b{X}$ equals zero. The system
$(6.1)$ is a basis of the vector space spanned by those $\b{X}$
over $\B{R}$. The correspondence $\b{x}\leftrightarrow\b{X}$
via the same definition as $(6.3)$ is linear, as is easily seen; and
the assertion $(6.6)$ means that this is in fact a Lie algebra
isomorphism between $\f{g}$ and the matrix Lie algebra
generated by $(6.1)$, which is a fact that extends 
to general Lie groups. Needless to say, this correspondence viewed
via $(6.3)$ is applicable only prior to the 
complexification mentioned above.
\medskip
\noindent
{\bf 7. Maass operators.} Following Maass [20], we introduce
$$
\eqalign{
{\bf e}^+&=2i{\bf x}_1+{\bf x}_2-i{\bf x}_3
=e^{2i\theta}\left(2iy{\partial\over\partial x}
+2y{\partial\over\partial y}
-i{\partial\over\partial \theta}\right),\cr
{\bf e}^-&=-2i{\bf x}_1
+{\bf x}_2+i{\bf x}_3=e^{-2i\theta}
\left(-2iy{\partial\over\partial x}
+2y{\partial\over\partial y}
+i{\partial\over\partial \theta}\right),\quad
{\bf w}={\bf x}_3={\partial\over\partial\theta}.
}\eqno(7.1)
$$
These three operators 
will play an important r\^ole in our later discussion.
Since we have
$$
[{\bf w},{\bf e}^+]=2i{\bf e}^+,\quad 
[{\bf w},{\bf e}^-]=-2i{\bf e}^-,
\quad [{\bf e}^+,{\bf e}^-]=-4i{\bf w},\eqno(7.2)
$$
as $(6.6)$ implies,
Maass operators
generate $\f{g}$ as well as $\e{U}$. The operator
$\bf w$ should not be confused with the 
Weyl element $\r{w}$. 
\medskip
\noindent
{\csc Notes:} The Maass operators $\b{e}^\pm$ are extensions of
the hyperbolic outer normal differential. See the notes to the
next section as well as Section 32. It should be stressed again that
the definition $(7.1)$ contains
in fact the complexification of the original Lie algebra $\f{g}$ 
over $\B{R}$.
\medskip
\noindent
{\bf 8. Casimir operator.} Next, we fix
the centre of $\e{U}$. To this end, we introduce the Killing 
form on $\f{g}\times\f{g}$:
$$
\r{Tr}\big((\r{ad}\,{\bf x})\!\cdot\!(\r{ad}\,{\bf y})\big),
\eqno(8.1)
$$
with $(\r{ad}\,\bf{x})(\b{a})=[\b{x},\b{a}]$. 
Computing the
coefficient matrix $(k_{ij})$ of $(8.1)$ via $(6.6)$,
we see that the form is non-degenerate; that is, $\r{G}$ is
semi-simple. More explicitly,
$$
\big(k_{ij}\big)=
\left(\matrix{\hfill0&0&-4\cr\hfill0&8&\hfill0\cr
-4\hfil&0&-8}\right),\quad
(k_{ij})^{-1}=\left(\matrix{\hfill{1\over2}&\hfill0&-{1\over4}\cr
\hfill0&{1\over8}&\hfill0\cr
-{1\over4}\hfil&0&\hfill0}\right).\eqno(8.2)
$$
We write the second matrix as
$(k^{ij})$. Then $c\sum_{i,j}k^{ij}\b{x}_i\b{x}_j$, 
with any $c\in\B{C}$, is an element 
in the centre of $\e{U}$. In this way we are led to 
the Casimir element of $\e{U}$:
$$
\eqalign{
\Omega&=-{\bf x}_1^2-{1\over4}{\bf x}_2^2
+{1\over2}{\bf x}_1{\bf x}_3
+{1\over2}{\bf x}_3{\bf x}_1\cr
&=-{1\over4}{\bf e}^+{\bf e}^-
+{1\over4}{\bf w}^2-{1\over2}i{\bf w}\cr
&=-y^2\left(\Big({\partial\over\partial x}\Big)^2
+\Big({\partial\over\partial y}\Big)^2\right)
+y{\partial^2\over\partial x\partial\theta},
}\eqno(8.3)
$$
with the choice $c=-2$. It holds that
$$
\b{u}\cdot\Omega=\Omega\cdot\b{u},\quad\forall
\b{u}\in\e{U}.\eqno(8.4)
$$
In fact, one may confirm by
using $(6.6)$ that
$[\Omega,\b{x}_j]=0$ holds for $j=1,2,3$; see the notes below.
The expression in the middle of $(8.3)$, which follows 
from the definition $(7.1)$, 
will play a central r\^ole in what follows. 
The third line is an easy consequence of the second.
\par
The Casimir operator is not only left invariant 
but also right invariant:
$$
l_\r{h}\Omega=\Omega l_\r{h},\quad
r_\r{h}\Omega=\Omega r_\r{h},
\quad \forall\r{h}\in\r{G}.\eqno(8.5)
$$ 
The first identity is trivial. As to the second, we note that
the definition $(6.3)$ implies that it holds, for any
$f\in C^\infty(\r{G})$, that
$$
\eqalignno{
\b{x}_jr_\r{h}f(\r{g})&=\Big[{d\over dt}\Big]_{t=0}
f\big(\r{g}\cdot\exp(t\,\b{X}_j)\r{h}\big)\cr
&=\Big[{d\over dt}\Big]_{t=0}
f\big(\r{gh}\cdot\exp(t\,\r{h}^{-1}\b{X}_j\r{h})\big)
=r_\r{h}\b{x}_j^\r{h}f(\r{g}),&(8.6)
}
$$
say. Thus, $\Omega r_\r{h}=r_\r{h}\Omega^\r{h}$ with
$\Omega^\r{h}$ being the result of 
replacing $\b{x}_j$ by $\b{x}_j^\r{h}$
in the first line of $(8.3)$. 
On the other hand, $\Omega^\r{h}=
\Omega$, since $\Omega$ does not depend 
on any base change of $\f{g}$ as is confirmed in what
follows; note that $\{\b{x}^\r{h}_1,
\b{x}^\r{h}_2,\b{x}^\r{h}_3\}$ is a basis, since
$ \r{h}^{-1}\b{X}_j\r{h}$, $j=1,2,3$, are linearly independent
in the vector space of $2\times2$ matrices of zero trace: 
Let $K$ be the matrix of the Killing form 
with respect to the original base, and $B$
a base change matrix. Then $BKB^t$ and thus 
$-2(B^t)^{-1}K^{-1}B^{-1}$ correspond, respectively,
to the Killing form and the Casimir element
on the new base. This obviously yields the
assertion. 
\medskip
\noindent
{\csc Notes:} Maass
defined his operators in a more general
fashion than $(7.1)$; see [20, Chapter 4, (12)--(13)]. 
 Roelcke [31] discussed them very thoroughly. 
Maass' basic motivation seems
to have been in finding a decomposition of 
the Casimir operator, which itself is also generalised in [20], into
a product of two elements of $\f{g}$, 
if formulated within our context, 
so that integration by parts works effectively; 
see [20, Chapter 4, (14)] as well as Section 32 below.
The decomposition in the second line
of $(8.3)$, which is a discovery of Maass, seems, however, to be 
a fortuitous situation special to the 
group $\r{PSL}(2,\B{R})$. 
It may be worth remarking that
Maass [20] does not exploit the fact that $\r{G}$ is
a Lie group; thus his approach is different from
ours. The commutativity
of the Casimir operator with right translations 
depends solely on the unimodularity $(5.2)$ of $\r{G}$; 
hence the assertion $(8.5)$ extends to any
Lie group which has a Haar measure that is left and right invariant.
We add also that
the confirmation $(8.4)$ can of course be made by using the
definition of the matrices $(k_{ij}),\,(k^{ij})$, and this
argument readily extends to general situation.
\medskip
\noindent
{\bf 9. Hilbert space $L^2(\varGamma
\backslash\r{G})$.} This is
defined to be the set of all
functions $f$ or vectors on $\r{G}$ which are left
$\varGamma$-automorphic or simply 
automorphic, 
i.e., $l_\gamma f=f$, $\forall\gamma\in\varGamma$,
and square integrable against the measure $d\r{g}$ over a fundamental 
domain $D$ of $\varGamma$ on $\r{G}$:
$$
\int_D|f(\r{g})|^2d\r{g}<+\infty.\eqno(9.1)
$$
Here $D$ is a $d\r{g}$-measurable subset of $\r{G}$
such that
$$
\r{G}=\bigcup_{\gamma\in\varGamma}\gamma D,\quad
\int_{\gamma D\cap D}d\r{g}=0,\;\gamma\ne1.
\eqno(9.2)
$$
For instance, the set 
$$
[\e{F}]=\big\{\r{n}[x]\r{a}[y]\r{k}[\theta]:
x+iy\in\e{F},\, 0\le\theta\le\pi\big\},\eqno(9.3)
$$
with $\e{F}$ defined by $(1.11)$,
serves the purpose. Obviously we may replace 
$D$ in $(9.1)$ by any
fundamental domain. As a consequence, we have
$$
\int_D|f(\r{gh})|^2d\r{g}=\int_D|f(\r{g})|^2d\r{g},
\quad\forall\r{h}\in\r{G}.\eqno(9.4)
$$
In fact, the unimodularity of $\r{G}$ asserted at $(5.2)$
implies that the left side is the same as integrating
$|f(\r{g})|^2$ against $d\r{g}$
over $D\r{h}$ which is also
a fundamental domain. 
\par
The set $L^2(\varGamma\backslash\r{G})$ is 
a Hilbert space equipped with the inner-product
$$
\big\langle f_1, f_2\big\rangle
=\int_{\varGamma\backslash\r{G}}
f_1(\r{g})\overline{f_2(\r{g})}d\r{g},\eqno(9.5)
$$
where the integration range is the whole quotient space 
$\varGamma\backslash\r{G}$ and the measure is
induced via $(5.2)$; the value of $(9.5)$ is naturally the same as
the result of integrating over any fundamental
domain.
\medskip
\noindent
{\csc Notes:} The definition $(9.5)$ stems from
Petersson's fundamental discovery [29]; 
see the notes to Section 21. The space $L^2(\varGamma
\backslash\r{G})$ is, more precisely, defined to be the
set of classes of $\varGamma$-automorphic 
functions, under the convention that two  functions
$f_1,f_2$ are in the same class
if and only if $\Vert f_1-f_2\Vert=0$, where the norm is associated
with $(9.5)$. In what follows we shall mostly deal with
relations between functions continuous throughout $\r{G}$; 
otherwise we shall mention explicitly.
\medskip
\noindent
{\bf 10. Automorphic representation.} The identity $(9.4)$
means that right translations are all unitary 
maps of $L^2(\varGamma\backslash\r{G})$ onto
itself: We have, for any vector $f$,
$$
\Vert r_\r{h}f\Vert =\Vert f\Vert,
\quad\forall\r{h}\in\r{G},
\eqno(10.1)
$$
with the norm as above. The map
$$
r:\r{h}\mapsto r_\r{h},\eqno(10.2)
$$
which is a homomorphism of $\r{G}$ into the
unitary transformation group of 
$L^2(\varGamma\backslash\r{G})$, is termed the right regular 
$\varGamma$-automorphic representation or just an
automorphic representation
of $\r{G}$. 
\par
Any closed subspace $W$ of 
$L^2(\varGamma\backslash\r{G})$ which satisfies
$r_\r{h}W\subseteq W$ for all $\r{h}\in\r{G}$, is 
called an invariant subspace. The orthogonal complement
of $W$ in terms of the metric $(9.5)$
is also an invariant subspace. We shall use
a representation and an invariant subspace as interchangeable
notions, under an obvious abuse of terminology.
If $W$ does not contain any non-trivial invariant subspace,
then it is said to be an irreducible subspace or 
representation. A major task of ours is to
establish a complete decomposition of 
$L^2(\varGamma\backslash\r{G})$ into a direct sum of
irreducible subspaces in an explicit fashion. 
It is immediate to notice that such a
decomposition should be closely
related to the spectral decomposition
of $\Omega$, since any of its eigenspaces is invariant as
the second identity in
$(8.5)$ implies. However, it turns out that
eigenvectors of $\Omega$
do not span the full space. The complement, which is fairly 
large, is filled with the contribution of the continuous spectrum of
$\Omega$, whose precise description is a salient aspect
of the harmonic analysis on $\varGamma\backslash\r{G}$. 
\medskip
\noindent
{\csc Notes:} There exists vast literature on
unitary representations
of Lie groups. An introductory account is given 
in Vilenkin--Klimyk [36,
Chapter 2]. However, the above definition
suffices for our purpose. It  should be noted that
the general theory of
representations of Lie groups requires that
homomorphisms corresponding to $(10.2)$ be 
strongly continuous. With our situation, this is inherent
in the definition $(10.2)$ itself : For any given vector $f$ and 
for arbitrary $\varepsilon>0$,
there exists a neighbourhood $\r{Q}$ of 
the unit element of $\r{G}$ such that
$\Vert r_\r{h}f-f\Vert<\varepsilon$, $\forall\r{h}\in\r{Q}$.
Indicating the proof, we choose a $\varGamma$-automorphic
$f_0\in C^\infty(\r{G})$ which is compactly
supported if restricted to $[\e{F}]$ and such that 
$\Vert{f-f_0}\Vert<\varepsilon$. Obviously
we have also $\Vert{r_\r{h}f-r_\r{h}f_0}\Vert<\varepsilon$,
$\forall\r{h}\in\r{G}$. 
Then, we consider the
continuity of the map $f_0\mapsto r_\r{h}f_0$ in the ordinary sense,
i.e., with respect to the norm $\Vert\cdot\Vert_\infty$. 
\medskip
\noindent
{\bf 11. Symmetry of $\Omega$.} Naturally it should be
made precise in which domain we consider the action of
$\Omega$. To this end as well as for the sake of
convenience of our discussion, we introduce
the linear set
$$
B^\infty(\varGamma\backslash\r{G})
=\left\{f\in C^\infty(\varGamma
\backslash\r{G}):\, {\hbox{$\b{u} f$ decays rapidly}
\atop\hbox{for any
fixed $\b{u}\in\e{U}$}}\right\},\eqno(11.1)
$$
where $C^\infty(\varGamma\backslash\r{G})
=C^\infty(\r{G})\cap L^2(\varGamma\backslash\r{G})$,
and decaying rapidly means 
that $\b{u}f(\r{g})\ll y^{-M}$ for any fixed
$M>0$ as $y\to+\infty$. 
This set is dense in $L^2(\varGamma\backslash\r{G})$,
for it contains any $f_0$ employed in the last notes. 
\par
Since the unimodularity $(5.2)$ of $\r{G}$
and the definition $(6.3)$
imply 
$$
\eqalignno{
\langle{\b{x}_jf_1,f_2}\rangle&
=\Big[{d\over dt}\Big]_{t=0}\int_{\varGamma
\backslash\r{G}}f_1\big(\r{g}\exp(t{\bf X}_j)\big)\overline
{f_2(\r{g})}d\r{g}\cr
&=\Big[{d\over dt}\Big]_{t=0}\int_{\varGamma
\backslash\r{G}}f_1\big(\r{g}\big)\overline
{f_2(\r{g}\exp(-t{\bf X}_j))}d\r{g}
=\langle{f_1,-\b{x}_jf_2}\rangle,&(11.2)
}
$$
we have, for any $\b{u}\in\e{U}$
and $f_1,f_2\in B^\infty(\varGamma\backslash\r{G})$,
$$
\eqalign{
\big\langle{\b{u}f_1,f_2}\big
\rangle&\,=\,\big\langle{f_1,\b{u}^*f_2}
\big\rangle,\cr \b{u}=\sum c\,
\b{x}_{j_1}\b{x}_{j_2}
\cdots\b{x}_{j_k}&\mapsto
\b{u}^*=\sum (-1)^k\overline{c}\,\b{x}_{j_k}
\b{x}_{j_{k-1}}
\cdots\b{x}_{j_1},
}\eqno(11.3)
$$
with obvious abbreviations.
In particular, the Casimir operator is symmetric
over $B^\infty(\varGamma\backslash\r{G})$:
$$
\big\langle{\Omega f_1,f_2}\big\rangle
=\big\langle{f_1,\Omega f_2}\big\rangle.\eqno(11.4)
$$
\medskip
\noindent
{\csc Notes:} The symmetry of the Laplace--Beltrami operator is
proved with Green's formula $(1.9)$. As the Casimir operator is
an extension of the Laplace--Beltrami operator, there should be
a proof of $(11.4)$ via an extension of Green's formula.
For this point see Section 32.
\medskip
\noindent
{\bf 12. Weights.} To achieve a decomposition of
$L^2(\varGamma\backslash\r{G})$ into irreducible
subspaces or rather
a spectral resolution of $\Omega$, there are at least
two ways for us to take: One is to exploit the existence of the
Maass operators. The other is more generally applicable and 
to rework our elementary argument
of [22, Chapter 1]. The first is to be rendered hereafter, and
the second to be summarised in Sections 32--37.
\par
With this, we first introduce the notion of weights: If a
function $f$ on $\r{G}$ is such that
there exists an $\ell\in\B{Z}$ satisfying
$$
f(\r{g}\r{k}[\tau])=e^{2\ell i\tau}f(\r{g}),
\eqno(12.1)
$$
then $f$ is said to be of even integral weight $2\ell$. 
 Since
functions on $\r{G}$ are of period $\pi$ 
with respect to the Iwasawa co-ordinate
$\theta$, weights ought to be even
integers. We have the orthogonal decomposition
$$
L^2(\varGamma\backslash\r{G})
=\mathop{\oplus}_{\ell=-\infty}^\infty 
L^2_\ell(\varGamma\backslash\r{G}),\eqno(12.2)
$$
where the $\ell$-th summand is the set composed of 
all vectors satisfying $(12.1)$, and the right side
is in fact the closure of the sum. This is the same
as the Fourier expansion in $\theta$ of vectors.
Let $g$ be in the $\ell$-th
summand, and let $h(z)=g(\r{n}[x]\r{a}[y])$. We have
$h(\gamma(z))=g(\r{n}[x_1]\r{a}[y_1])$
with $\gamma\cdot\r{n}[x]\r{a}[y]=
\r{n}[x_1]\r{a}[y_1]\r{k}[\theta_1]$; thus
$h(\gamma(z))=g(\gamma\cdot\r{n}[x]\r{a}[y])
\exp(-2\ell i\theta_1)$.
 In view of the last identity in $(3.4)$, we have
$$
h(\gamma(z))=h(z)\bigg({\jmath(\gamma,z)\over
\jmath(\gamma,\bar{z})}\bigg)^\ell,
\quad\forall\gamma\in\varGamma.\eqno(12.3)
$$
That is, $h$ is a $\varGamma$-automorphic function 
on $\B{H}^2$ of weight $2\ell$. Hence
$$
L^2_\ell(\varGamma\backslash\r{G})
= L^2_\ell(\varGamma\backslash\B{H}^2)
\cdot\exp(2\ell i\theta),\eqno(12.4)
$$
with
$$
L^2_\ell(\varGamma\backslash\B{H}^2)
=\left\{{\hbox{$\varGamma$-automorphic of weight $2\ell$
and }\atop\hbox{square integrable
over $\e{F}$ against $d\mu$}}\right\},\eqno(12.5)
$$
which is obviously an extension of $(1.12)$.
Applied to $\Omega$, this separation of variables 
gives
$$
\Omega_\ell=-y^2\left(\Big({\partial\over\partial x}\Big)^2
+\Big({\partial\over\partial y}\Big)^2\right)
+2i\ell y{\partial\over\partial x}.\eqno(12.6)
$$
\medskip
\noindent
{\csc Notes:} The extension to $\Omega_\ell$ 
of the invariance assertions and Green's formula given
in the first section is developed in Section 32. 
It is possible to consider not only
even integral but also arbitrary complex weights, which
is due originally to Maass [20] and thoroughly 
investigated by Roelcke [31]. The view of such extensions
from the standpoint of the universal covering group of $\r{G}$, 
which is in fact a generalisation of the basic stance 
of the present article, is given in Bruggeman [4].
\medskip
\noindent
{\bf 13. Descending and ascending.} We employ 
the Maass operators in discussing 
relations among weight strata $(12.2)$ of 
$L^2(\varGamma\backslash\r{G})$; note that
we may restrict ourselves to the linear set 
$B^\infty(\varGamma\backslash\r{G})$.
The operator $\b{w}$ detects a stratum, and $\b{e}^\pm$ 
shift it: We pick up a vector $g$ in
$$
B^\infty_\ell(\varGamma\backslash\r{G})=
L^2_\ell(\varGamma\backslash\r{G})\cap
B^\infty(\varGamma\backslash\r{G}).\eqno(13.1)
$$ 
Then, we have $\b{w}g=2\ell ig$, and $(7.2)$ implies
$\big(\b{w}-2(\ell\pm1)i\big)\b{e}^\pm g=0$; 
namely, the weight of 
$\b{e}^\pm g$ is $2(\ell\pm1)$. When $\ell>0$, the
application of $(\b{e}^-)^{\ell}$ to 
$B^\infty_\ell(\varGamma\backslash\r{G})$ 
takes the set down to 
$L^2_0(\varGamma\backslash\r{G})
=L^2_0(\varGamma\backslash\B{H}^2)$, 
where we have a spectral
resolution of $\Omega_0$, 
a fairly elementary proof of which
is achieved in [22, Chapter 1]. 
Hence, the image
$(\b{e}^-)^{\ell}B^\infty_\ell(\varGamma\backslash\r{G})$ 
admits a spectral decomposition. To lift it up to the original
space, we apply $(\b{e}^+)^\ell$. In this way, one might think that
a spectral decomposition of $L^2_\ell(\varGamma
\backslash\r{G})$ could be achieved. 
However, as a matter of fact the
image $(\b{e}^+)^\ell (\b{e}^-)^\ell B^\infty_\ell
(\varGamma\backslash\r{G})$ 
does not span the space 
$L^2_\ell(\varGamma\backslash\r{G})$ 
in general. The discrepancy is closely related to the
notion of holomorphic
cusp forms on $\B{H}^2$, as we shall make 
precise in Sections 19--21. 
\par
It should be added that the case $\ell<0$ is analogous, 
because of the involution
$$
J:\,\r{n}[x]\r{a}[y]\r{k}[\theta]
\mapsto \r{n}[-x]\r{a}[y]\r{k}[-\theta],
\eqno(13.2)
$$
which is the same as the map
$\left[{a\atop c}{b\atop d}\right]\mapsto
\left[{\hfill a\atop -c}{-b\atop\hfill d}\right]$
in $\r{G}$.
\medskip
\noindent
{\csc Notes:} This mechanism among the weight strata is a
discovery of Maass [20, Chapter 4]; 
see the two formulas following (13) there, although
he formulated it in a quite generalised setting and without the notion of
Lie differentials. See also Roelcke [31, Teil I, \S3].
\medskip
\noindent
{\bf 14. Jacquet operator.} In order to make
the procedure in the last section explicit, we need to
render the spectral decomposition of $L^2_0(\varGamma
\backslash\r{G})$ in such a fashion that applications of
Maass operators may be 
performed in a conspicuous way. To this end, we
introduce the operator of Jacquet [13]: For a function $\phi$ on
$\r{G}$, we put
$$
\displaystyle{\e{A}^\delta\!\phi(\r{g})=\int_{-\infty}^\infty 
e(-\delta\xi)\phi(\r{w}\r{n}[\xi]\r{g})d\xi,}\atop{ 
e(\xi)=\exp(2\pi i\xi),
\quad \delta=\pm,\quad \r{w}=\r{k}
\big[\txt{1\over2}\pi\big].}\eqno(14.1)
$$
In what follows, we shall apply $\e{A}^\delta$
to those $\phi$ with which
ordinary convergence holds, and analytic continuation with respect
to parameters involved in $\phi$
will be taken into account, if needed.
Computation of the integral is carried out on
noting that by $(3.4)$ the map $\r{g}
\mapsto\r{w}\r{n}[\xi]\r{g}$ is the same as
$$
\eqalign{
&x\mapsto {-x-\xi\over (x+\xi)^2+y^2},\quad
y\mapsto {y\over (x+\xi)^2+y^2},\cr
&\hskip 1cm e^{2 i\theta}\mapsto e^{2 i\theta}\cdot
{x+\xi-iy\over x+\xi+iy}.
}\eqno(14.2)
$$
We have
$$
{r_\r{h}\e{A}^\delta=\e{A}^\delta r_\r{h},
\quad \forall\r{h}\in\r{G},}\atop
{\b{u}\e{A}^\delta=\e{A}^\delta\b{u},\quad
\forall\b{u}\in\e{U}}.\eqno(14.3)
$$
The first is obvious, and the second is a simple consequence
of the definition $(6.3)$, although we need to
have adequate smoothness of $\phi$ in $(14.1)$. The merit of
having the Jacquet operator will be felt in Section 18 and 
thereafter. 
\medskip
\noindent
{\csc Notes:} The definition $(14.1)$ stems
from Fourier expansions of Poincar\'e series 
on the big cell, i.e., the term $\r{NwNA}$, 
of the Bruhat decomposition:
$$
\r{G}=\r{NA}\sqcup\r{NwNA}.\eqno(14.4) 
$$
Thus the operator $\e{A}^\delta$ is a fairly natural device.
\medskip
\noindent
{\bf 15. Weight functions.}
With the specialisation 
$$
\phi_\ell(\r{g},\nu)=y^{\nu+1/2}\exp(2\ell i\theta),
\quad \ell\in\B{Z},
\eqno(15.1)
$$
which is termed a weight function, we have, by $(14.2)$,
$$
\e{A}^\delta\!\phi_\ell(\r{g},\nu)
=\exp(2\ell i\theta)e(\delta x)
y^{{1/2}-\nu}\int_{-\infty}^\infty{e(-y\xi)
\over(\xi^2+1)^{\nu+1/2}}
\left({\xi-i\over\xi+i}\right)^{\delta \ell}d\xi.
\eqno(15.2)
$$
This will play a fundamental r\^ole throughout the rest of
our discussion. Applying integration by
parts sufficiently many times, we see that
the integral converges for any $\nu$, and
$\e{A}^\delta\phi_\ell(\r{g},\nu)$ is in fact entire 
in $\nu$. More conspicuously, we have
$$
\eqalignno{
\e{A}^\delta\phi_\ell(\r{g},\nu)&=
{\pi^\nu\exp(2\ell i\theta)e(\delta x)
\over\Gamma\big(\nu+|\ell|+{1\over2}\big)}
\sum_{j=0}^{|\ell|}(-\pi)^j{2|\ell|\choose 2j}
\Gamma\big(|\ell|-j+\txt{1\over2}\big)\cr
\times\, y^{j+1/2}&\int_0^\infty u^{\nu-1}
\big(\sqrt{u}+\delta\sgn(\ell)/\sqrt{u}\big)^{2j}
\exp\left(-\pi y\big(u+1/u\big)\right)du,&(15.3)
}
$$
which is related to the familiar
Schl\"afli integral representation for
$K_\nu$, the $K$-Bessel function of order $\nu$. 
To show this, we assume $\Re\nu>0$, and express 
$\Gamma\big(\nu+|\ell|+{1\over2}\big)(\xi^2+1)^
{-\nu-|\ell|-1/2}$ in terms of the Euler integral 
over $u>0$ for the
Gamma-function. We insert it into $(15.2)$ and 
exchange the order of
integration, getting
$(\xi-\delta\sgn(\ell)i)^{2|\ell|}
\exp(-2\pi iy\xi-u\xi^2)$ as the new inner integrand. 
We shift the inner contour to
$\Im\xi=-\pi y/u$. 
Then $(15.3)$ follows after a 
rearrangement. The regularity assertion with respect
to $\nu$ is now immediate. 
Also, we note that a shift
of the contour in $(15.2)$ to $\Im\xi=-\infty$ gives,
for $l\in\B{N}$,
$$
\e{A}^\delta\phi_\ell\big(\r{g},
l-\txt{1\over2}\big)=\cases{\hfil 0& if 
$\delta\ell\le -l$,\cr
\displaystyle{(-1)^l{(2\pi)^{2l}\over\Gamma(2l)}
y^l \exp(2\ell i\theta)e(\delta x+iy)}& if $\delta\ell=l$.}
\eqno(15.4)
$$
\par
The formula $(15.3)$ implies that $\e{A}^\delta\phi_\ell
(\r{g},\nu)$ is of exponential decay in $y$, which is, however,
often inadequate and we need 
uniform bounds: For instance,
we have, for $\Re\nu>-{1\over2}$,
$$
{\e{A}^\delta\phi_\ell(\r{g},\nu)}
\ll (|\nu|+|\ell|+1)\cdot
\cases{\hfil y^{1/2-|\Re\nu|}(|\log y|+1)&if
$0<y<1$,\cr y^{-1/2-\Re\nu}
\exp\big(-y/(|\nu|+|\ell|+1)\big)
&if $1\le y$,}\eqno(15.5)
$$
where the implied constant depends on $\Re\nu$ only. 
In fact, for $y\ge1$, it suffices to apply the integration by parts 
to $(15.2)$ and shift the contour to 
$\Im\xi=-(|\nu|+|\ell|+1)^{-1}$.
For $0<y<1$, we use the trivial identity
$$
\eqalignno{
\exp(-2\ell i\theta)&e(-\delta x)
\e{A}^\delta\phi_\ell(\r{g},\nu)=
\e{A}^\delta\phi_0(\r{a}[y],\nu)\cr
&+y^{{1/2}-\nu}\int_{-\infty}^\infty{e(-y\xi)
\over(\xi^2+1)^{\nu+1/2}}
\left(\Big({\xi-i\over\xi+i}\Big)^{\delta \ell}
-1\right)d\xi.&(15.6)
}
$$
The second term on the right side is obviously 
$\ll y^{1/2-\Re\nu}
(|\ell|+1)$. The first term is dealt with
$(15.3)$, for $\ell=0$, $\r{g}=\r{a}[y]$. Assuming
$|\nu|$ is large, we turn the
contour through the angle 
${1\over2}(\pi-{1/|\nu|})\sgn(\Im\nu)$ around
the origin. The rest of the argument may be skipped.
\medskip
\noindent
{\bf 16.  Whittaker functions.}
In literature, the transform $\e{A}^\delta\phi_\ell(\r{g},\nu)$
is expressed in terms of 
Whittaker functions $W_{\mu,\nu}$:
$$
\e{A}^\delta\phi_\ell(\r{g},\nu)
=(-1)^\ell\pi^{\nu+1/2}\exp(2\ell i\theta)e(\delta x)
{W_{\delta \ell,\nu}(4\pi y)\over\Gamma
\big(\nu+\delta \ell+{1\over2}\big)}.\eqno(16.1)
$$
We shall prove this for a practical reason as well as
in respect for tradition. 
We begin with the definition
$$
W_{\mu,\nu}(y)=e^{\mu\pi i}
{\Gamma\big(\mu+\nu+{1\over2}\big)
\over \pi(y/4)^{\nu-1/2}}\int_{-\infty}^\infty
{e^{-iy\xi/2}\over (\xi^2+1)^{\nu+{1/2}}}
{(\xi-i)^\mu\over(\xi+i)^\mu} d\xi,
\eqno(16.2)
$$
where $y>0$, $\mu\in\B{C}$, and
$\arg(\xi\pm i)$ varies from $\pm\pi$ to $0$, respectively, 
along the contour. It will be justified in due
course that $(16.2)$ serves the purpose of 
defining Whittaker functions. Thus
we transform $(16.2)$
into an expression of the Mellin--Barnes type:
$$
W_{\mu,\nu}(y)
={y^\mu e^{-y/2}\over2\pi i}
\int_{-i\infty}^{i\infty}
{\Gamma\big(s-\mu+\nu+{1\over2}\big)
\Gamma\big(s-\mu-\nu+{1\over2}\big)
\over\Gamma\big({1\over2}-\mu+\nu\big)
\Gamma\big({1\over2}-\mu-\nu\big)}
\Gamma(-s)y^{-s}ds,\eqno(16.3)
$$
where the path separates the poles of
$\Gamma\big(s-\mu+\nu+{1\over2}\big)
\Gamma\big(s-\mu-\nu+{1\over2}\big)$ and
$\Gamma(-s)$ to the left and the right, respectively;
and it is assumed that parameters are such that the
path can be drawn. To show this,
we consider the integral
$$
\int_{Y_\varepsilon}
{e^{-iy\xi/2}d\xi\over (\xi+i)^{\mu+\nu+{1/2}}
(\xi-i)^{-\mu+\nu+1/2}},\eqno(16.4)
$$
where the
path $Y_\varepsilon$ with a small $\varepsilon>0$
starts at $-1/\varepsilon$, proceeds to $1/\varepsilon$ on
the real axis, then to $-i/\varepsilon$ along the 
circle $|\xi|=1/\varepsilon$ and goes up to
$-(1+\varepsilon)i$ on the imaginary axis;
it encircles $-i$ counter-clockwise, goes 
down to $-i/\varepsilon$ on the imaginary axis and
returns to $-1/\varepsilon$ along the circle
$|\xi|=1/\varepsilon$.  Under the temporary
condition $0<\Re\nu<-\Re\mu+{1\over2}$, the
contribution of the circular part of the path vanishes as
$\varepsilon\to +0$, and one may see that
the integral in $(16.2)$ equals
$$
2^{1-2\nu}e^{-i\pi\mu}e^{-y/2}
\sin\left(\pi(\mu+\nu+\txt{1\over2})\right)
\int_0^\infty{e^{-y\rho}d\rho\over \rho^{\mu+\nu+{1/2}}
(1+\rho)^{-\mu+\nu+1/2}}. \eqno(16.5)
$$
We apply the Mellin inversion to the numerator of
the integrand and exchange the order of integrals.
The new inner integral is Euler's Beta function, and 
after a rearrangement together with
analytic continuation we reach $(16.3)$. 
Shifting the contour in $(16.3)$
sufficiently far to the right we see immediately that
$$
W_{\mu,\nu}(y)=\big(1+o(1)\big)y^\mu e^{-y/2},
\quad y\to +\infty,\eqno(16.6)
$$
uniformly for any bounded $\mu,\nu$. Also it holds,
for $y>0$, that
$$
\eqalign
{&\eqalign{
&\big[\omega^++\mu\big]W_{\mu,\nu}(y)
=-W_{\mu+1,\nu}(y),\cr
&\big[\omega^--\mu\big]W_{\mu,\nu}(y)
=\big((\mu-\txt{1\over2})^2-\nu^2\big)
W_{\mu-1,\nu}(y),}
\quad \omega^\pm=y(d/dy)\mp{1\over2}y,\cr
&\bigg[\!-\Big({d\over dy}\Big)^2+{1\over4}
-{\mu\over y}+\left(\nu^2-\txt{1\over 4}\right)
{1\over y^2}\bigg]W_{\mu,\nu}(y)=0.}\eqno(16.7)
$$
The two recurrence equations can be confirmed by
applying the respective differentials
to $(16.3)$, which may be
performed freely as the integral converges
rapidly because of Stirling's formula
for the $\Gamma$-function. 
The third line in $(16.7)$ is a consequence of the
first two, and this confluent hypergeometric differential 
equation is customarily 
attributed to Whittaker. In view of $(16.6)$ we find that
$(16.2)$ can indeed be employed to define
Whittaker functions; see [38, Chapter XVI].
\par
Here is a beautiful integral formula:
For $\alpha,\beta\in{\Bbb C}$ and $|\Re\nu|<{1\over2}$,
$$
\eqalignno{
&\int_0^\infty W_{\alpha,\nu}(y) 
W_{\beta,\nu}(y){dy\over y}=
{\pi\over(\alpha-\beta)\sin(2\pi\nu)}\cr
&\times\Bigg[{1\over\Gamma({1\over2}-\alpha+\nu)
\Gamma({1\over2}-\beta-\nu)}
-{1\over\Gamma({1\over2}-\alpha-\nu)
\Gamma({1\over2}
-\beta+\nu)}\Bigg].\quad\qquad&(16.8)
}
$$ 
To show this, we note first that it holds
uniformly for $|\Re\nu|<{1\over2}$,
bounded $\mu$ and sufficiently small $y>0$ that
$$
\eqalign{
W_{\mu,\nu}(y)
&=\!\left({\Gamma(-2\nu)
\over\Gamma\big({1\over2}-\mu-\nu\big)}
y^{\nu+{1/2}}+{\Gamma(2\nu)\over\Gamma\big({1\over2}
-\mu+\nu\big)}y^{-\nu+{1/2}}\right)\big(1+O(y)\big),\cr
W'_{\mu,\nu}(y)
&=\!\left({\big(\nu+{1\over2}\big)
\Gamma(-2\nu)\over\Gamma\big({1\over2}-\mu-\nu\big)}
y^{\nu-{1/2}}-{\big(\nu-{1\over2}\big)\Gamma(2\nu)
\over\Gamma\big({1\over2}-\mu+\nu\big)}
y^{-\nu-{1/2}}\right)\big(1+O(y)\big).
}\eqno(16.9)
$$
In fact, we replace $s$ by $s+\mu-{1\over2}$ in $(16.3)$
and shift the path to $\Re(s)=-M+{1\over2}$
with a sufficiently large $M\in\B{N}$. An examination of residues
gives $(16.9)$. Then, with
$\omega_\nu=-(d/dy)^2+{1\over4}
+(\nu^2-{1\over 4})y^{-2}$ and by integration by parts, we
have, for $\alpha\ne\beta$,
$$
\eqalignno{ 
&(\alpha-\beta)\int_0^\infty
W_{\alpha,\nu}(y)W_{\beta,\nu}(y){dy\over y}\cr
&=\lim_{\varepsilon\to+0}
\int_\varepsilon^\infty \Big[\omega_\nu
W_{\alpha,\nu}(y)W_{\beta,\nu}(y)-W_{\alpha,\nu}(y)
\omega_\nu W_{\beta,\nu}(y)\Big]dy\cr
&=\lim_{\varepsilon\to+0}\Big[W_{\alpha,\nu}'(\varepsilon)
W_{\beta,\nu}(\varepsilon)-W_{\alpha,\nu}
(\varepsilon)W_{\beta,\nu}' (\varepsilon)\Big],&(16.10)
}
$$ 
since the first integral converges absolutely if
$|\Re\nu|<{1\over2}$, because of $(16.6)$ and 
the first line of $(16.9)$. The last limit 
can be computed by combining the two formulas of $(16.9)$.
\medskip 
\noindent
{\csc Notes:} The formula $(16.8)$ is tabulated
at [9, 7.611(3)]; however, the precious factor $\pi$ is 
missing there, which would cause a discrepancy concerning a
unitarity assertion given in Section 27.
Our proof is taken from [23, Part XII].  
The formula $(16.3)$ is given in 
Whittaker--Watson [38, p.\ 343] and the
recursive formulas $(16.7)$ in Vilenkin--Klimyk
[36, p.\ 218], but our treatments are different from theirs.
The assertion $(16.7)$ is related to the actions of the
Maass and the Casimir operators.
\medskip
\noindent
{\bf 17. Hilbert space $L_0^2(\varGamma\backslash
\r{G})$.} We read [22, Theorem 1.1] in the present
context, and have
$$
L_0^2(\varGamma\backslash\r{G})=\B{C}\cdot1
\oplus {}^0\!L_0^2(\varGamma\backslash\r{G})
\oplus {}^e\!L_0^2(\varGamma\backslash\r{G}).
\eqno(17.1)
$$
Here
$$
{}^0\!L_0^2(\varGamma\backslash\r{G})=
\mathop{\oplus}_V\B{C}\cdot
\lambda_V^{(0)},\quad 
\big\langle{\lambda_V^{(0)},\lambda_{V'}^{(0)}}
\big\rangle=\delta_{V,V'}\,,\eqno(17.2)
$$
with
$$
\displaystyle{\lambda_V^{(0)}(\r{g})=\sum_{n\ne0}
{\varrho_V(n)\over\sqrt{|n|}}\e{A}^{\sgn(n)}\phi_0
(\r{a}[|n|]\r{g},\nu_V),}\atop{ 
\Omega\lambda_V^{(0)}
=\left(\txt{1\over4}-\nu_V^2\right)\lambda_V^{(0)},
\quad \nu_V\in i\B{R},}
\eqno(17.3)
$$
where $V$'s are just  labels indexing 
the discrete set $\big\{{1\over4}-\nu_V^2\big\}$ 
of eigenvalues of $\Omega_0$ acting 
over $L^2_0(\varGamma\backslash\r{G})$;  
in Sections 23 and later
they will stand for a series of
invariant subspaces of $L^2(\varGamma\backslash\r{G})$.
Naturally, the right side of $(17.2)$ is to be understood as the
closure of the sum.
 By the identity $(15.3)$ with
$\ell=0$, we have in fact
$$
\lambda_V^{(0)}(\r{g})={2\pi^{\nu_V+1/2}\over
\Gamma\big(\nu_V+{1\over2}\big)}
\sum_{n\ne0}\varrho_V(n)y^{1/2}K_{\nu_V}(2\pi|n|y)
e(nx),\eqno(17.4)
$$
which converges absolutely.
Namely, we have written the Fourier
expansion [22, $(1.1.41)$] in such a way that
$$
\kappa_j=-i\nu_V,\quad \rho_j(n)
={2\pi^{\nu_V+1/2}\over\Gamma(\nu_V+{1\over2})}
\varrho_V(n).\eqno(17.5)
$$
As to ${}^e\!L_0^2(\varGamma\backslash\r{G})$, which is
the contribution of the continuous spectrum of $\Omega_0$,
we have
$$
{}^e\!L_0^2(\varGamma\backslash\r{G})
=\left\{\int_{(0)}h(\nu)E_0(\r{g},\nu)
d\nu:h\in L^2(i\B{R})\right\},\eqno(17.6)
$$
where the kernel 
$E_0$ is the Eisenstein series 
of weight $0$ to be defined in the next section and
$L^2(i\B{R})$  the $L^2$-space with respect to the
Lebesgue measure placed on the imaginary axis $(0)
=i\B{R}$; the integral converges in the mean
in the space $L_0^2(\varGamma\backslash\r{G})$. 
The decomposition $(17.1)$ with $(17.2)$ and $(17.6)$
is equivalent to the spectral expansion of any
$f\in L^2_0(\varGamma\backslash\r{G})$:
$$
f(\r{g})={3\over\pi}\big\langle f, 1\big\rangle+
\sum_V\big\langle f,\lambda_V^{(0)}\big\rangle
\lambda_V^{(0)}(\r{g})
+{1\over4\pi i}\int_{(0)}\e{E}_0(f,\nu)
E_0(\r{g},\nu)d\nu,\eqno(17.7)
$$
where
$$
\e{E}_0(f,\nu)=\int_{\varGamma\backslash\r{G}}f(\r{g})
\overline{E_0(\r{g},\nu)}d\r{g}
\quad\hbox{in $\; L^2(i\B{R})$}.\eqno(17.8)
$$
The identity $(17.7)$ holds, with the sum and the
integral  converging in the mean. Further,
$(17.8)$ is in fact the limit in the mean 
in the space $L^2(i\B{R})$ of the
integrals over $[\e{F}]_Y=[\e{F}]\cap\{y\le Y\}$,
with $Y$ tending to $+\infty$.
\par
We may assume that with $J$ defined by $(13.2)$
$$
J\lambda_V^{(0)}=\epsilon_V\lambda_V^{(0)},
\quad\epsilon_V=\pm1,\eqno(17.9)
$$
which is the same as [22, $(3.1.15)$] and equivalent
to $\varrho_V(-n)=\epsilon_V\varrho_V(n)$ for
any $n\in\B{N}$. We have the bounds
$$
\sum_{|\nu_V|\le K} 1\ll K^A,\eqno(17.10)
$$
with an absolute constant $A>0$,
and
$$
\sum_{|\nu_V|\le K}|\varrho_V(n)|^2
\ll K^2+n^{4/5},\eqno(17.11)
$$
where the implied constants are absolute. The former
is a simple consequence of [22, $(1.4.2)$]; 
or rather $(37.21)$ below gives $A\le8$. 
In Section 43 
we shall show an asymptotic formula, with $A=2$, via Selberg's
trace formula. Although weaker than
[22, (2.3.2)], the bound $(17.11)$ is adequate 
for our purpose; it depends
on a fairly elementary bound for Kloosterman sums instead
of Weil's. 
\par
\medskip
\noindent
{\csc Notes:} The concept of automorphic eigenfunctions
is due to Delsarte [7] and Maass [19]. 
The spectral decomposition
$(17.1)$ is an instance of Selberg's general 
statement [32], although the detailed proof
of the spectral resolution of the Casimir operator, in Maass'
extended context, is done by Roelcke in [31] with 
an essential appeal to functional analysis,
especially to the theory
of unbounded symmetric operators  
as well as to the theory of
elliptic differential operators. The proof of $(17.1)$ given in [22] is
elementary in the sense that it is virtually independent
of functional analysis, save
for a use of the classical Hilbert--Schmidt theory
on integral operators with bounded continuous symmetric
kernels. We should mention also that
there exists an elementary
proof of Weil's bound for Kloosterman sums; the theory
is initiated by Stepanov [34].
\medskip
\noindent
{\bf 18. Eisenstein series.} For an arbitrary even integral
weight $2\ell$, this is defined by
$$
E_\ell(\r{g},\nu)=\sum_{\gamma\in\varGamma_\infty
\backslash\varGamma}\phi_\ell(\gamma\r{g},\nu),
\quad\Re\nu>{1\over2},\eqno(18.1)
$$
with $\varGamma_\infty=\big\{\left[{1\atop}{n\atop
1}\right]:\, n\in\B{Z}\big\}$, which converges absolutely.
We have the Fourier expansion
$$
\eqalignno{
E_\ell(\r{g},\nu)&=\phi_\ell(\r{g},\nu)
+{(-1)^\ell\Gamma^2\big(\nu+{1\over2}\big)
\varphi_\varGamma(\nu)
\over\Gamma\big(\nu+|\ell|+{1\over2}\big)
\Gamma\big(\nu-|\ell|+{1\over2}\big)}
\phi_\ell(\r{g},-\nu)\cr
&+{1\over\zeta(2\nu+1)}
\sum_{n\ne0}{|n|^{-\nu}
\sigma_{2\nu}(|n|)\over\sqrt{|n|}}
\e{A}^{\sgn(n)}\phi_\ell(\r{a}[|n|]\r{g},\nu),&(18.2)
}
$$
with $\sigma_\alpha(n)=\sum_{d|n}d^\alpha$ and
$$
\varphi_\varGamma(\nu)=\sqrt{\pi}{\Gamma(\nu)\zeta(2\nu)
\over\Gamma(\nu+{1\over2})\zeta(2\nu+1)}
=\pi^{2\nu}{\Gamma\big(\txt{1\over2}-\nu\big)\zeta(1-2\nu)
\over\Gamma\big(\nu+\txt{1\over2}\big)
\zeta(2\nu+1)}.\eqno(18.3)
$$
The identity $(15.3)$ implies that $E_\ell(\r{g},\nu)$ is 
a meromorphic function of $\nu$ over $\B{C}$. 
The functional equation for the modified Eisenstein series
$$
{E_\ell^*(\r{g},\nu)=E^*_\ell(\r{g},-\nu),}\atop
{E_\ell^*(\r{g},\nu)=
\pi^{-\nu-1/2}\Gamma\big(\nu+|\ell|+\txt{1\over2}\big)
\zeta(2\nu+1)E_\ell(\r{g},\nu),}
\eqno(18.4)
$$
is a consequence of 
the symmetry $W_{\mu,\nu}=W_{\mu,-\nu}$,
which follows from $(15.3)$ and $(16.1)$. We have that
$E^*_\ell(\r{g},\nu)$, $\ell\ne0$, is entire in 
$\nu$, while $E^*_0(\r{g},\nu)$ is regular except for the
simple poles at $\nu=\pm{1\over2}$. 
\par
The proof of these facts 
on $E_\ell$ may be skipped as it is a standard application of
Poisson's sum formula via the double coset decomposition
$\varGamma_\infty\backslash\varGamma/
\varGamma_\infty$. It should, however, be stressed that
it suffices, as a matter of fact, to have the expansion
$(18.2)$ for $E_0$ only, since we have 
$$
({\bf e}^{\sgn(\ell)})^{|\ell|}\phi_0(\r{g},\nu)=2^{|\ell|}
{\Gamma\big({1\over2}+\nu+|\ell|\big)\over
\Gamma({1\over2}+\nu)}
\phi_\ell(\r{g},\nu),\eqno(18.5)
$$
that is, in view of $(14.3)$ we obtain $(18.2)$ for general
$\ell\ne0$ by an application of  
$(\bf{e}^{\sgn(\ell)})^{|\ell|}$ 
to the expansion for $E_0$. 
More precisely, we argue as follows: Under the convention
introduced after $(6.7)$, we first apply 
$(\bf{e}^{\sgn(\ell)})^{|\ell|}$ to the defining expression 
for $E_0$ term-wise, which is
legitimate, since the result converges
absolutely and uniformly when $\Re\nu>{1\over2}$ in view of
the left invariance of $\b{e}^\pm$ and $(18.5)$. We get
the defining expression for $E_\ell$, save for a
factor. On the other hand,
the expansion $(18.2)$ for $E_0$ admits the 
term-wise application of $(\bf{e}^{\sgn(\ell)})^{|\ell|}$, for
the result converges absolutely and uniformly
in view of $(14.3)$, $(15.5)$ and $(18.5)$. In this way
we may derive $(18.2)$ from the
expansion for $E_0$. This practical mechanism 
is a merit of having the Jacquet operator, and will be exploited
throughout the rest of our discussion.
\medskip
\noindent
{\bf 19. Discrepancy.} We resume the scheme
started in Section 13; thus we let $\ell$ be positive. 
We apply $(17.7)$ to 
$({\bf e}^{-})^\ell g(\r{g})$, and get, on noting $(11.3)$,
$$
({\bf e}^{-})^\ell g(\r{g})
=(-1)^\ell\sum_{V}\big\langle{g,
({\bf e}^{+})^\ell\lambda_V^{(0)}}\big\rangle
\lambda_V^{(0)}(\r{g})
+{(-1)^\ell\over4\pi i}\int_{(0)}
\eta_\ell(g,\nu)
E_0\left(\r{g},\nu\right)d\nu,\eqno(19.1)
$$
with
$$
\eta_\ell(g,\nu)=
\int_{\varGamma\backslash\r{G}}g(\r{g})
\overline{\big(({\bf e}^+)^\ell E_0\big)
\left(\r{g},\nu\right)}d\r{g}.\eqno(19.2)
$$
We have used that $({\bf e}^+)^\ell1=0$ and
that $\e{E}_0(({\bf e}^{-})^\ell g,\nu)=(-1)^\ell
\eta_\ell(g,\nu)$. The latter can be shown in 
the same way as $(11.3)$, since 
$g\in B^\infty_\ell(\varGamma\backslash\r{G})$. 
\par
Before applying $(\b{e}^+)^\ell$ to $(19.1)$,
we note that for any $\ell\ge0$
$$
({\bf e}^-)^\ell({\bf e}^+)^\ell\phi_k(\r{g},\nu)=
(-4)^\ell{\Gamma\big({1\over2}+\nu+k+\ell\big)
\Gamma\big({1\over2}-\nu+k+\ell\big)
\over\Gamma\big({1\over2}+\nu+k\big)
\Gamma\big({1\over2}-\nu+k\big)}
\phi_k(\r{g},\nu).\eqno(19.3)
$$
In fact, a combination of 
$(7.2)$, $(8.3)$ and $(8.4)$
gives
$$
\eqalignno{
&({\bf e}^-)^{\ell+1}({\bf e}^+)^{\ell+1}\phi_k(\r{g},\nu)\cr
=&\,({\bf e}^-)^\ell\big(-4\Omega
+{\bf w}^2+2i{\bf w}\big)
({\bf e}^+)^\ell\phi_k(\r{g},\nu)\cr
=&\,-4\big(\txt{1\over2}+\nu+k+\ell\big)
\big(\txt{1\over2}-\nu+k+\ell\big)
({\bf e}^-)^\ell({\bf e}^+)^\ell\phi_k(\r{g},\nu).&(19.4)
}
$$
We apply $(19.3)$, for $k=0$, to $(17.3)$ and $(18.2)$
specialised with $E_0$:
$$
\eqalign{
({\bf e}^-)^\ell({\bf e}^+)^\ell\lambda_V^{(0)}(\r{g})
&=(-4)^\ell{\big|\Gamma
\big({1\over2}+\nu_V+\ell\big)\big|^2
\over\big|\Gamma\big({1\over2}+\nu_V\big)\big|^2}
\lambda_V^{(0)}(\r{g}),\cr
({\bf e}^-)^\ell({\bf e}^+)^\ell E_0(\r{g},\nu)&=
(-4)^\ell{\big|\Gamma\big({1\over2}+\nu+\ell\big)\big|^2
\over\big|\Gamma\big({1\over2}+\nu\big)\big|^2}
E_0(\r{g},\nu),\quad \nu\in i\B{R},
}\eqno(19.5)
$$
in which the convention following $(6.7)$ is in effect.
\par
Then we put
$$
\eqalignno{
g(\r{g})&=\skew3\tilde g(\r{g})
+2^{-2\ell}\sum_V
{\big|\Gamma\big({1\over2}+\nu_V\big)\big|^2
\over\big|\Gamma\big({1\over2}+\nu_V+\ell\big)\big|^2}
\big\langle{g,({\bf e}^+)^\ell\lambda_V^{(0)}}
\big\rangle({\bf e}^+)^\ell\lambda_V^{(0)}(\r{g})\cr
&+{2^{-2\ell}\over4\pi i}\int_{(0)}
{\big|\Gamma\big({1\over2}+\nu\big)\big|^2
\over\big|\Gamma\big({1\over2}+\nu
+\ell\big)\big|^2} \eta_\ell(g,\nu)
({\bf e}^+)^\ell E_0(\r{g},\nu)d\nu.
&(19.6)
}
$$
We apply $(\b{e}^-)^\ell$ on both sides and use $(19.5)$.
Comparing the result with $(19.1)$, we find that
$$
({\bf e}^-)^\ell\skew3\tilde g=0.\eqno(19.7)
$$
Namely, $\skew3\tilde{g}$ stands for the discrepancy indicated 
in Section 13.
\medskip
\noindent
{\csc Notes:} This discrepancy or rather $(19.7)$ is treated
in Roelcke [31, Teil II] from the standpoint of analysing the
situation where $(18.5)$ vanishes, if formulated in our notation.
We have taken the above approach in order to further stress the
relevance to the structure of the weight strata.
\medskip
\noindent
{\bf 20. Convergence issue.} However, we have yet to
see whether the last operation on $(19.6)$ is legitimate or not. 
To this end, we shall first prove that the sum
on the right side of $(19.1)$ belongs to $B^\infty(
\varGamma\backslash\r{G})$: We observe that
for any $q\in\B{N}$
$$
\big|\big\langle{g,({\bf e}^+)^\ell\lambda_V^{(0)}}
\big\rangle\big|
={1\over({1\over4}-\nu_V^2)^q}
\big|\big\langle{\Omega^q({\bf e}^-)^\ell g,
\lambda_V^{(0)}}\big\rangle\big|
\le{\Vert{\Omega^q({\bf e}^-)^\ell g}\Vert
\over\big({1\over4}-\nu_V^2\big)^q}
\ll |\nu_V|^{-2q}.\eqno(20.1)
$$
Thus, we have, for any $k\in\B{N}$,
$$
\eqalignno{
&\sum_V\big|\big\langle{g,({\bf e}^+)^\ell
\lambda_V^{(0)}}\big\rangle 
(\b{e}^\pm)^k\lambda_V^{(0)}(\r{g})\big|
\ll\sum_{n\ne0}{1\over\sqrt{|n|}}\sum_K K^{k-2q}\cr
\times&\left(\sum_{K\le|\nu_V|\le2K}
|\varrho_V(n)|^2\right)^{1/2}
\left(\sum_{K\le|\nu_V|\le2K}
|\e{A}^{\sgn(n)}\phi_{\pm k}(\r{a}[|n|]\r{g},\nu_V)|^2
\right)^{1/2},\qquad&(20.2)
}
$$
where $K$ runs over dyadic numbers, and $(18.5)$, for $\ell=
\pm k$, is applied. The first sum over $\nu_V$ 
on the right side is estimated 
by $(17.11)$, and the second by $(15.5)$ while taking 
account of $(17.10)$. We then
divide the sum over $K$ at $K\approx(|n|y)^{1/2}$. 
This yields the desired assertion, since
$\{\b{e}^-,\b{e}^+,\b{w}\}$ generates $\e{U}$ and
the action of $\b{w}$ is immaterial. 
\par
We need further to
deal with the integrated term of $(19.1)$: 
We have a bound of $\eta_\ell(g,\nu)$
similar to $(20.1)$, since $\Omega^q(\b{e}^-)^\ell g$
is in $B^\infty(\varGamma\backslash\r{G})$ and
$E_0(\r{g},\nu)$, $\Re\nu=0$, is of polynomial
order in both $\r{g}$ and $\nu$ as is 
implied by $(18.2)$ together with
$(15.5)$ and a well-known lower bound 
for $\zeta(s)$ on $\Re s=1$. By the same token,
$(\b{e}^\pm)^k E_0(\r{g},\nu)$
is of polynomial order in $\nu$, $\Re\nu=0$, uniformly
for bounded $\r{g}$.  Hence it is immaterial
whether $(\b{e}^\pm)^k$ is applied inside
or outside the integral. The result
of the outside application is in $B^\infty(\varGamma
\backslash\r{G})$, as has been proved already.
\medskip
\noindent
{\bf 21. Holomorphic cusp forms.} We shall show that
$$
\skew3\tilde g=\sum_{l=0}^\ell ({\bf e}^+)^{\ell-l}\varphi_l,
\eqno(21.1)
$$
where
$$
B_l^\infty(\varGamma\backslash\r{G})\ni
\varphi_l(\r{g})=e^{2li\theta}y^l\sum_{n=1}^\infty
a(n)e(nz),\quad z=x+iy.\eqno(21.2)
$$
We employ the induction in terms of $\ell$, as 
the case $\ell=0$ is trivial. Thus, let
$p\in B^\infty_{\ell+1}(\varGamma\backslash\r{G})$ be
such that $({\bf e}^-)^{\ell+1}p=0$; this
smoothness of $p$ can be assumed because of 
the assertion of the previous section. By the
inductive hypothesis we have
${\bf e}^-p=\sum_{l=0}^\ell 
({\bf e}^+)^{\ell-l}\varphi_{1,l}$ with the specification
same as $(21.1)$. We apply $-{1\over4}{\bf e}^+$ 
to both sides, and get, on noting the second line of $(8.3)$,
$$
\eqalignno{
(\Omega+\ell(\ell+1))p
&=-{1\over4}\sum_{l=0}^\ell ({\bf e}^+)^{\ell+1-l}
\varphi_{1,l}\cr
&=-{1\over4}\sum_{l=0}^\ell ({\bf e}^+)^{\ell+1-l}
{(\Omega+\ell(\ell+1))\over
\ell(\ell+1)-l(l-1)}\varphi_{1,l}\,,&(21.3)
}
$$
which means that 
$$
p=q-{1\over4}\sum_{l=0}^\ell
{({\bf e}^+)^{\ell+1-l}\varphi_{1,l}\over
\ell(\ell+1)-l(l-1)},
\quad (\Omega+\ell(\ell+1))q=0.\eqno(21.4)
$$
We note that $q\in B_{\ell +1}^{\infty}
(\varGamma\backslash\r{G})$ by the construction.
Again by the second line of $(8.3)$, we have
$0=\big\langle{(\Omega+\ell(\ell+1))q,\,q}
\big\rangle={1\over4}\Vert{\bf e}^-q\Vert^2$;
namely ${\bf e}^-q=0$. Then, in the expansion
$$
q(\r{g})=e^{2(\ell+1)i\theta}
\sum_{n=-\infty}^\infty c(n)k_n(y)e(nx),\eqno(21.5)
$$
we have $\big(d/dy+2\pi n-(\ell+1)/y\big)k_n=0$, and thus
$k_n(y)=b(n)y^{\ell+1}\exp(-2\pi ny)$.
The rapid decay of $q$ implies
$b(n)=0$, $n\le0$, which ends the proof of $(21.1)$.
\par
Now, let $\varpi(z)=(e^{2il\theta} y^l)^{-1}\varphi_l(\r{g})$.
Then by $(12.3)$, for $\ell=l$, we see that $\varpi$ is regular 
throughout $\B{H}^2$ and satisfies
$$
\varpi(i\infty)=0,\quad \varpi(\gamma(z))
(\jmath(\gamma,z))^{-2l}
=\varpi(z),\quad\gamma\in\varGamma.\eqno(21.6)
$$
That is, $\varpi$ is a holomorphic cusp form of
weight $2l$ with respect to $\varGamma$.
Let $\e{C}_\varGamma(l)$ be the set of all such functions.
Then $\e{C}_\varGamma(l)$ is a 
Hilbert space of finite dimension $\vartheta_\varGamma(l)$ 
equipped with the inner product
$$
\big\langle{\varpi_1,\varpi_2}\big\rangle_l
=\int_{\eusm{F}}\varpi_1(z)\overline{\varpi_2(z)}\,
y^{2l}d\mu(z),\eqno(21.7)
$$
where $\e{F}$ is as in $(1.11)$. The dimension formula
$$
\vartheta_\varGamma(l)=\cases{\hfil[l/6]&
$l\not\equiv1\bmod6$,\cr
[l/6]-1& $l\equiv1\bmod6,\, l\ge7$,}
\eqno(21.8)
$$ 
is well-known. 
\par
Let $\big\{(e^{2il\theta} y^l)^{-1}
\lambda_V^{(l)}(\r{g})\big\}$ be an
orthonormal basis of $\e{C}_\varGamma(l)$,
with $V$ running over $\vartheta_\varGamma(l)$ labels.
In particular, we have
$$
\big\langle{\lambda_V^{(l)},\lambda_{V'}^{(l)}}
\big\rangle=\delta_{V,V'};\eqno(21.9)
$$
that is, the
inner product $(21.7)$ on $\e{C}_\varGamma(l)$ is
translated into $(9.5)$ on $L^2(\varGamma\backslash
\r{G})$.
 We write
the Fourier expansion of $\lambda_V^{(l)}$ as
$$
\displaystyle{\lambda_V^{(l)}(\r{g})=\pi^{1/2-l}
\Gamma(2l)^{1/2}
\sum_{n=1}^\infty{\varrho_V(n)\over\sqrt{n}}
\e{A}^+\!\phi_l\big(\r{a}[n]\r{g},
l-\txt{1\over2}\big),}\atop{
\Omega\lambda_V^{(l)}=
\big(\txt{1\over4}-\nu_V^2\big)
\lambda_V^{(l)},\quad \nu_V=l-\txt{1\over2},
}\eqno(21.10)
$$
which converges absolutely, and 
is to be compared with $(17.3)$.
This notation may appear misleading but its
employment will be justified in the next section.
In view of the lower line of $(15.4)$, 
the expansion $(21.10)$ means that we have 
re-normalised the Petersson--Fourier
coefficients [22,  $(2.2.3)$] as
$$
\rho_{j,l}(n)=(-1)^l{2^{2l}\pi^{l+1/2}\over
\Gamma(2l)^{1/2}}\varrho_V(n).\eqno(21.11)
$$
Corresponding to $(17.11)$, we have,
with $V$ as in $(21.10)$,
$$
\sum_V|\varrho_V(n)|^2\ll l+n^{4/5},\eqno(21.12)
$$
where the implied constant is absolute,
which contains 
a simple  improvement upon [22, (2.2.10)] 
concerning the dependency
on weights; see [27, Vol.\ 2, (4.1.21)]. 
\medskip
\noindent
{\csc Notes:} For the proof of $(21.8)$ see, e.g., 
Maass [20, Chapter 2]; Petersson's theory 
of Poincar\'e series is applied. In the notes to Section 38
is an indication of a proof via Selberg's trace formula.
We should not miss mentioning, even though this historical fact is
very much well-known, that $(21.7)$ is the metric introduced by
Petersson [29], by which he initiated the modern development of
the theory of automorphic functions and automorphic 
representations such as the spectral decomposition of
$L^2(\varGamma\backslash\r{G})$ stated in Section 23.
\medskip
\noindent
{\bf 22. Hilbert space $L^2_\ell(\varGamma\backslash
\r{G})$.} We sum up the above 
discussion: The expansion $(19.6)$ together
with $(21.1)$ makes precise the spectral
structure of $L^2_\ell(\varGamma
\backslash\r{G})$, $\ell>0$.
First, in the context of Section 17,
 we put
$$
\lambda_V^{(\ell)}(\r{g})
={\Gamma\big({1\over2}+\nu_V\big)
\over 2^\ell\Gamma\big(\txt{1\over2}+\nu_V+\ell\big)}
({\bf e}^+)^\ell\lambda_V^{(0)}(\r{g}).\eqno(22.1)
$$
By $(17.3)$ and $(18.5)$ we have
$$
\lambda_V^{(\ell)}(\r{g})=\sum_{\scr{n=-\infty}\atop
\scr{n\ne0}}^\infty{\varrho_V(n)\over\sqrt{|n|}}
\e{A}^{\sgn(n)}\!\phi_\ell(\r{a}[|n|]\r{g},\nu_V),
\quad \nu_V\in i\B{R}.
\eqno(22.2)
$$
Second, in the context of the previous section,
we put, for $\ell\ge l$,
$$
\lambda_V^{(\ell)}(\r{g})=2^{l-\ell}
\left({\Gamma(2l)\over\Gamma(\ell+l)
\Gamma(\ell-l+1)}\right)^{1/2}({\bf
e}^+)^{\ell-l}\lambda_V^{(l)}(\r{g}).\eqno(22.3)
$$
By $(21.10)$ we have
$$
{\displaystyle
\lambda_V^{(\ell)}(\r{g})
=\pi^{1/2-l}\left({\Gamma(\ell+l)\over
\Gamma(\ell-l+1)}\right)^{1/2}
\sum_{n=1}^\infty{\varrho_V(n)\over\sqrt{n}}
\e{A}^+\!\phi_\ell\big(\r{a}[n]\r{g},\nu_V\big),}
\atop{\displaystyle\nu_V=l-{1\over2},\quad l\in\B{N}.}
\eqno(22.4)
$$
\par
With this, we have
$$
\big\langle{\lambda_V^{(\ell)},
\lambda_{V'}^{(\ell)}}\big\rangle
=\delta_{V,V'}.\eqno(22.5)
$$
Namely, the system consisting of 
functions $(22.2)$ and $(22.4)$
is orthonormal in $L^2_\ell(\varGamma\backslash
\r{G})$. Thus, if $\lambda_V^{(\ell)}$ and
$\lambda_{V'}^{(\ell)}$ are both in the category
$(22.1)$, then $(17.2)$ and $(19.5)$ give the assertion;
and the case with the category $(22.3)$ is treated via
the identities $(19.3)$ and $(21.9)$. Further, if these
functions belong to different categories, then they are
orthogonal, since their $\Omega$-eigenvalues are
different and we may appeal to $(11.4)$; or rather
the fact comes down to $(\b{e}^-)^\ell\lambda_V^{(\ell)}\equiv0$,
with $\lambda_V^{(\ell)}$ defined by $(22.3)$,
which may also be used to show the orthogonality of these
$\lambda_V^{(\ell)}$ against
the integrated part of $(19.6)$.
As a consequence, we can rewrite $(21.1)$ as
$$
\skew3\tilde g(\r{g})=\sum_{l=1}^\ell\sum_{\scr{V}
\atop\scr{\nu_V=l-{1\over2}}}
\big\langle{g,\lambda_V^{(\ell)}}\big\rangle 
\lambda_V^{(\ell)}(\r{g}).\eqno(22.6)
$$
In fact, it suffices to plug $(19.6)$ into the right side.
\par
Hence we have
$$
L_\ell^2(\varGamma\backslash\r{G})=
{}^0\!L_\ell^2(\varGamma\backslash\r{G})
\oplus {}^e\!L_\ell^2(\varGamma\backslash\r{G}),\quad
\ell>0,
\eqno(22.7)
$$
with
$$
{}^0\!L_\ell^2(\varGamma\backslash\r{G})
={\mathop{\oplus}_V}^{(\ell)}\,\B{C}\cdot\lambda_V^{(\ell)}
\eqno(22.8)
$$
and
$$
{}^e\!L_\ell^2(\varGamma\backslash\r{G})
=\left\{\int_{(0)} \eta(\nu)E_\ell(\r{g},\nu)
d\nu:\eta\in L^2(i\B{R})\right\}.\eqno(22.9)
$$
Or equivalently, we have, 
for any $f\in L_\ell^2(\varGamma\backslash\r{G})$,
$$
f(\r{g})={\sum_V}^{(\ell)}
\big\langle f,\lambda_V^{(\ell)}\big\rangle
\lambda_V^{(\ell)}(\r{g})
+{1\over4\pi i}\int_{(0)}\e{E}_\ell(f,\nu)
E_\ell(\r{g},\nu)d\nu,\eqno(22.10)
$$
with
$$
\e{E}_\ell(f,\nu)=\int_{\varGamma\backslash\r{G}}f(\r{g})
\overline{E_\ell(\r{g},\nu)}d\r{g}
\quad\hbox{in $\; L^2(i\B{R})$}.\eqno(22.11)
$$
Here both 
$\oplus^{(\ell)}$ and $\sum^{(\ell)}$ denote that
$\nu_V$ is either pure imaginary or equal to
$l-{1\over2}$, $l\in\B{N}$, $l\le\ell$. The right side of $(22.8)$
is understood similarly to that of $(17.2)$.
The integrated part of
$(22.10)$ is a consequence of a combination
of $(18.1)$, $(18.5)$, $(19.2)$ and the corresponding part of
$(19.6)$. The identity $(22.10)$ holds in the same sense
as $(17.7)$, with $(22.11)$ being an analogue of $(17.8)$.
\par
The case $\ell<0$ is analogous; it suffices to
apply the map $J$ to $(22.7)$. Since
$J\r{g}_1\r{g}_2=J\r{g}_1 J\r{g}_2$ and
$d\r{g}=dJ\r{g}$, we have
$\lambda_V^{(|\ell|)}(J\r{g})\in L^2_\ell
(\varGamma\backslash\r{G})$. Also,
$\exp(2li\theta)y^{-l}\lambda_V^{(l)}(J\r{g})$
with $\lambda_V^{(l)}$ as in $(21.10)$ is
an antiholomorphic cusp form of weight $2l$.
\medskip
\noindent
{\csc Notes:} The spectral decomposition $(22.7)$--$(22.11)$
can be approached via the the Green function for $\Omega_\ell$.
Salient points of the argument will be indicated in Sections 32--37.
\medskip
\noindent
{\bf 23. Spectral decomposition of 
$L^2(\varGamma\backslash\r{G})$.} Combining
the assertions $(12.2)$, $(17.1)$--$(17.2)$, $(17.6)$,
and $(22.7)$--$(22.9)$ we obtain our main result. This
is essentially the same as the implication of
$(17.7)$ and $(22.10)$. Or more drastically rendering,
it is the same as
the rearrangement of the result so far established by
exchanging the order of the indices $V$ and $\ell$:
\smallskip
\noindent
{\csc Theorem.} {\it Let
${}^0\!L^2(\varGamma\backslash\r{G})$ be the 
cuspidal subspace
of $L^2(\varGamma\backslash\r{G})$ which is
spanned by all vectors whose constant terms in the
Fourier expansion with respect to the left action of
$\r{N}$ vanish, 
and let ${}^e\!L^2(\varGamma\backslash\r{G})$ 
be the subspace generated by integrals of all
Eisenstein series $E_\ell$, $\ell\in\B{Z}$,
as indicted by $(22.9)$. Then we have
$$
\displaystyle{L^2(\varGamma\backslash\r{G})
={\Bbb C}\cdot1\oplus
{}^0\!L^2(\varGamma\backslash\r{G})
\oplus{}^e\!L^2(\varGamma\backslash\r{G})}\atop
{\displaystyle{{}^0\!L^2(\varGamma\backslash\r{G})
=\oplus V,\quad V=\mathop\oplus_\ell 
{\Bbb C}\lambda_V^{(\ell)},}
\atop \displaystyle{{}^e\!L^2(\varGamma\backslash\r{G})
=\mathop\oplus_{\ell=-\infty}^\infty
{}^e\!L^2_\ell(\varGamma\backslash\r{G}).
}}\eqno(23.1)
$$
Here $V$'s are all irreducible subspaces with
$\Omega\lambda_V^{(\ell)}
=\big(\txt{1\over4}-\nu_V^2\big)\lambda_V^{(\ell)}$.
If $V$ is generated by a real analytic cusp 
form via $(22.1)$, that is, $\nu_V\in i\B{R}$,
then the corresponding index $\ell$ runs over all integers.
Otherwise, either $\ell\ge l$ or $\ell\le-l$ according as
$V$ is generated by a holomorphic cusp form of positive even 
integral weight $2l$ via $(22.3)$ or via an
antiholomorphic cusp form analogously, that is, $\nu_V=
l-{1\over2}$.
Further, we have, for any
$f\in L^2(\varGamma\backslash\r{G})$,
$$
\eqalignno{
f(\r{g})&={3\over\pi}\big\langle f,1\big\rangle+
\sum_V\sum_\ell\big\langle f,\lambda_V^{(\ell)}\big\rangle
\lambda_V^{(\ell)}(\r{g})\cr
&+\sum_{\ell=-\infty}^\infty
{1\over4\pi i}\int_{(0)}\e{E}_\ell(f,\nu)
E_\ell(\r{g},\nu)d\nu,&(23.2)
}
$$
where the range of $\ell$ in the 
first line is the same as in the corresponding part of
$(23.1)$, and $\e{E}_\ell(f,\nu)$ is defined by $(22.11)$,
although $f$ is not restricted to
$L^2_\ell(\varGamma\backslash\r{G})$.
}
\medskip
\noindent
{\csc Notes:} More precisely, each
$V$ stands for the closure of the sum $\mathop\oplus_\ell 
{\Bbb C}\lambda_V^{(\ell)}$; and
${}^0\!L^2(\varGamma\backslash\r{G})$, 
and ${}^e\!L^2(\varGamma\backslash\r{G})$ are understood
similarly. Also the identity $(23.2)$ holds, with 
sums and integrals converging in the mean.
\medskip
\noindent
{\bf 24. Series of irreducible
representations.} The invariance and the
irreducibility of each $V$ remains to be established;
we shall achieve it in due course. Thus the
title of the present section
might be premature. We shall proceed 
with caution so that the use of this notion 
should not cause any confusion.
\par
We introduce first the classification of $V$'s:
$$
\eqalign{
&\hbox{Unitary principal series $\Leftarrow$ real
analytic cusp forms,}\cr
&\hbox{Discrete series $\Leftarrow$
either holomorphic or antiholomorphic cusp forms,}
}\eqno(24.1)
$$
where the arrows mean that each series are generated
by respective variety of cusp forms on $\B{H}^2$;
thus the latter splits into the holomorphic and
antiholomorphic discrete series.
This is the same as
$$
\eqalign{
&\hbox{$V$ in the unitary principal series $\Leftrightarrow$
$\nu_V\in i\B{R}$,}\cr
&\hbox{$V$ in the discrete series $\Leftrightarrow$
$\nu_V\in \B{N}-{1\over2}$,}
}\quad\quad\hbox{with $\Omega\vert_V=
\big(\txt{1\over4}-\nu_V^2\big)\cdot 1\,$.}
 \eqno(24.2)
$$
In general, there possibly exists
the complementary series, whose constituents are
generated by real analytic 
cusp forms associated with exceptional
eigenvalues of $\Omega_0$. With $\varGamma=
\r{PSL}(2,\B{Z})$ such forms do not exist as $(17.7)$ indicates;
see [22, Lemma 1.4]. 
\par
We observe that the Fourier expansions $(22.2)$ and
$(22.4)$ of the basis vectors $\lambda_V^{(\ell)}$ 
as well as their $J$-images 
are conveniently expressed as
$$
\lambda_V^{(\ell)}(\r{g})=
\bigg|\pi^{-2\nu_V}{\Gamma\big(|\ell|+\nu_V
+{1\over2}\big)
\over\Gamma\big(|\ell|-\nu_V+{1\over2}\big)}
\bigg|^{1/2}
\sum_{n\ne0}{\varrho_V(n)\over\sqrt{|n|}}
\e{A}^{\sgn\!(n)}\!\phi_\ell(\r{a}[|n|]\r{g},\nu_V),
\eqno(24.3)
$$
where $(15.4)$ is to be invoked to see that $(22.4)$ is
well included. With this,
$$
\hbox{the sequence $\{\varrho_V(n): |n|\in\B{N}\}$ is 
dependent solely on the space $V$,}\eqno(24.4)
$$
except for an arbitrary multiplier of unit absolute value. These
coefficients do not depend on weights.
Hence, it is appropriate to designate them 
as the Fourier coefficients 
of an irreducible subspace or representation $V$. 
It should be stressed that this normalisation of
Fourier coefficients of cusp forms 
has been made possible by the use of
the Jacquet operator $\e{A}^\delta$
and the weight function $\phi_\ell(\r{g},\nu)$.
\medskip
\noindent
{\bf 25. Local structure.} In order to illuminate
the grand assertion $(23.1)$,
we fantasise that each subspace
$L^2_\ell(\varGamma\backslash\r{G})$ is the
galaxy $\ell$. Then we observe that all the galaxies revolve
at respective angular velocities under the right
action of $\r{K}$. They
are, however, not independent
of each other. Maass operators transport vectors
from a galaxy to others
along light-ladders $V$, most of which start in
the fount $L^2_0(\varGamma\backslash\r{G})\equiv
L^2(\varGamma\backslash\B{H}^2)$ 
and extend to infinity in both directions, while
others emerge spontaneously and discretely
in pairs above and below the fount and 
extend to respective infinity. 
The principal part of each galaxy is a slice
of the sum of $V$'s in the unitary principal series, which
is a unitary image of
the fount. The rest is a slice of the sum of $V$'s in
the discrete series.
To see how the latter expands
with $|\ell|$, one should observe that a particular $V$
in the discrete series can emerge in the galaxy $\ell$
only when $\nu_V=|\ell|-{1\over2}$.
\par
This is a view which although fairly beautiful is
nevertheless not of much use, especially in applications to
problems in analytic number theory. What matters 
in practice is nothing else but 
to have precise analytical structures 
of individual light-ladders $V$ so that
the projection to $V$, i.e., the relevant sum 
over the weights in $(23.2)$, 
of a given vector can be computed explicitly and
effectively, without recourse to deeper natures
of the Fourier coefficients $\{\varrho_V(n)\}$,
preferably independently of them. We stress in this context
that the construction of the
subspace ${}^e\!L^2(\varGamma\backslash\r{G})$ 
is fairly representational and visible as Eisenstein series
$E_\ell$ are defined by $(18.1)$ and have expansions
$(18.2)$, whereas basis vectors $(24.3)$ of
the cuspidal subspace are far more
abstract objects as the nature of the numbers
$\{\varrho_V(n)\}$ remains largely mysterious.
\par
Hence, our next task is to cast light on the
structure of each $V$ in the cuspidal
subspace, an issue which is local in the spectral context 
but central for practical purposes. 
\medskip
\noindent
{\csc Notes:} The Fourier coefficients 
$\sigma_{2\nu}(n)$ correspond
to $\varrho_V(n)$. The former has
a visible inner structure as being expressed in terms of a sum
over divisors of $n$, which is not well shared by the
latter, even though the Hecke operators disclose
similarities between these Fourier coefficients. However,
this ostensive nature of 
divisor functions is misleading.
They are in fact equally mysterious, as is demonstrated,
for instance, by
the additive divisor problem which concerns the sum
$$
\sum_{n=1}^\infty\sigma_\alpha(n)\sigma_\beta(n+f)
W(n/f),\eqno(25.1)
$$
where $\alpha,\beta\in\B{C}$, 
$f\in\B{N}$, and the weight function $W$
is supposed to be sufficiently smooth and of rapid decay; see [21]
and [27, Section 6.4].
\medskip
\noindent
{\bf 26. Invariance.} We shall first prove that the $V$'s in $(23.1)$
are all invariant with respect to the right action of $\r{G}$. 
Dealing with a $V$ in the unitary principal 
series, i.e., $\nu_V\in i\B{R}$, we put
$$
 V^\infty=\left\{\hbox{$\displaystyle\lambda(\r{g})=
\sum_{\ell=\infty}^\infty c_\ell
\lambda_V^{(\ell)}(\r{g}):\,c_\ell\ll(|\ell|+1)^{-M}$
with any $M>0$}\right\},\eqno(26.1)
$$
where the implied constant may depend on $M$. 
The absolute convergence of the sum can be
confirmed by means of $(15.5)$ and $(17.11)$;
we plug $(24.3)$ into the definition of $\lambda(\r{g})$, and
exchange the order of summation, which is justified
in much the same way as in $(20.2)$, although this
time we sum over $\ell$ instead of $V$. Then we note that
$\e{A}^\delta$ can be applied term-wise to
$$
\displaystyle{\phi(\r{g})=\sum_{\ell=-\infty}^\infty
c_\ell\phi_\ell(\r{g},\nu_V)=y^{1/2+\nu_V}\Phi(\theta),}
\atop\displaystyle{\Phi(\theta)=\sum_{\ell=-\infty}^\infty 
c_\ell e^{2\ell i\theta}\in C^\infty\big(\B{R}/\pi\B{Z}\big).}
\eqno(26.2)
$$
In fact we apply integration by parts to $(15.2)$, for
$\nu=\nu_V$; then,
after exchanging the order of integration and
summation, we undo the
integration by parts. In this way we get
$$
\lambda(\r{g})=\sum_{n\ne0}{\varrho_V(n)
\over\sqrt{|n|}}\e{A}^{\sgn(n)}\phi(\r{a}[|n|]\r{g}).
\eqno(26.3)
$$
Before applying right translations,
we invoke $(5.4)$ and have that for $\r{h}=\r{n}[\alpha]
\r{a}[\beta]\r{k}[\tau]$
$$
\eqalignno{
r_\r{h}\phi(\r{g})
&=y^{1/2+\nu_V}{\beta^{1/2+\nu_V}\Phi(\tau+
\vartheta(\theta))\hfil\over
\big((\cos\theta-\alpha\sin\theta)^2
+(\beta\sin\theta)^2\big)^{1/2+\nu_V}}\cr
&=\sum_{\ell=-\infty}^\infty
c_\ell^\r{h}\phi_\ell(\r{g},\nu_V),\quad
c_\ell^\r{h}\ll (|\ell|+1)^{-M},&(26.4)
}
$$
for any $M>0$, since the first line belongs to
$C^\infty(\B{R}/\pi\B{Z})$ as a function of $\theta$.
Reversing the order of reasoning,
we find that $r_\r{h}\lambda(\r{g})\in V^\infty$,
i.e.,  $r_\r{h} V^\infty\subseteq V^\infty$. As $V^\infty$
is dense in $V$ and $r_\r{h}$ is unitary, 
we conclude that 
$$
\hbox{the $V$'s in the unitary principal series are all invariant.}
\eqno(26.5)
$$
\par
We next consider a $V$ in the holomorphic discrete series;
thus $\nu_V=l-{1\over2}$, $l\in\B{N}$. 
In place of $(26.1)$--$(26.3)$, we put
$$
 V^\infty=\left\{\hbox{$\displaystyle\lambda(\r{g})=
\sum_{\ell= l}^\infty c_\ell
\lambda_V^{(\ell)}(\r{g}):\,c_\ell\ll(\ell+1)^{-M}$
with any $M>0$}\right\},\eqno(26.6)
$$
$$
\phi(\r{g})=\pi^{1/2-l}y^l\sum_{\ell=l}^\infty 
c_\ell\left({\Gamma(\ell+l)
\over\Gamma(\ell-l+1)}\right)^{1/2}
 e^{2\ell i\theta},\eqno(26.7)
$$
$$
\lambda(\r{g})=\sum_{n=1}^\infty{\varrho_V(n)
\over\sqrt{|n|}}\e{A}^+\phi(\r{a}[|n|]\r{g}).
\eqno(26.8)
$$
The verification of the last identity may be skipped, as it is
analogous to that of $(26.3)$. However, in order to 
show the analogue of $(26.4)$ 
we need a minor contrivance: Via $(5.4)$ we have,
with $z=e^{i\theta}$,
$$
\eqalignno{
r_\r{h}\phi(\r{g})=\,&\pi^{1/2-l}
(4\beta y)^l z^{2l}\sum_{\ell=l}^\infty
c_\ell e^{2\ell i\tau}\left({\Gamma(\ell+l)
\over\Gamma(\ell-l+1)}\right)^{1/2}\cr
&\times{\big((1+\beta+\alpha i)z^2
+1-\beta-\alpha i\big)^{\ell-l}
\over\big((1-\beta+\alpha i)z^2
+1+\beta-\alpha i\big)^{\ell+l}}.&(26.9)
}
$$
This sum is an even regular function for $|z|<1$ and
also belongs to $C^\infty(\B{R}/\pi\B{Z})$ as
a function of $\theta$. We conclude that
$r_\r{h}V^\infty\subseteq V^\infty$.
Thus, together with an application of the
map $(13.2)$, we conclude that
$$
\hbox{the $V$'s in the discrete series are all invariant.}
\eqno(26.10)
$$
\par
We may now use safely the term `a representation $V$'. The
irreducibility is to be proved in Section 30.
\medskip
\noindent
{\csc Notes:} The invariance assertion is
usually confirmed via the Lie algebra $\f{g}$. Our
argument may appear to be unconventional. However,
the discussion based on explicit group actions 
should also be worth reporting. 
\medskip
\noindent
{\bf 27. Kirillov map.} We shall deal with
the task set out in Section 25. To this end, we 
introduce the map:
$$
\e{K}\phi(u)=\e{A}^{\sgn(u)}\phi(\r{a}[|u|]),
\eqno(27.1)
$$
following Kirillov [16], where $\phi$ can be any function 
on $\r{G}$ as
far as the relevant integral converges in the same sense
as in $(14.1)$. In what follows, 
we are concerned mainly with the specialisation
$\phi(\r{g})=\phi_\ell(\r{g},\nu)$, i.e.,
$\e{K}\phi(u)=\e{K}\phi_\ell(u,\nu)$.
We shall show that
$$
\hbox{for each $\nu\in i\B{R}$, 
the set $\big\{\e{K}\phi_\ell(u,\nu):\ell\in\B{Z}\big\}$ is}\atop
\hbox{a complete orthonormal system of
$L^2(\B{R}^\times,d^\times\!/\pi)$,}\eqno(27.2)
$$
as well as that
$$
\hbox{for each $l\in \B{N}$, 
the set $\big\{
\pi^{1/2-l}(\Gamma(\ell+l)/\Gamma(\ell-l+1))^{1/2}
\e{K}\phi_\ell\big(u,l-{1\over2}\big):
\ell\ge l\big\}$ is}\atop
\hbox{a complete orthonormal system of
$L^2(\B{R}_+^\times,d^\times\!/\pi)$}.\eqno(27.3)
$$
Here $\B{R}^\times=\B{R}\!\setminus\!\{0\}$,
$d^\times\! u=du/|u|$; and $\B{R}_+^\times$ is the 
set of positive real numbers.
\par
By $(16.1)$ we have, for 
$\nu\in i{\Bbb R}$ and $\ell,\ell'\in{\Bbb Z}$,
$$
\eqalignno{
\big\langle\e{K}\phi_\ell(\cdot,\nu),&\,
\e{K}\phi_{\ell'}(\cdot,\nu)\big\rangle
={(-1)^{\ell+\ell'}\over
\Gamma\big(\ell+\nu+{1\over2}\big)
\Gamma\big(\ell'-\nu+{1\over2}\big)}
\int_0^\infty W_{\ell,\nu}(y)
W_{\ell',\nu}(y){dy\over y}\qquad\cr
&+{(-1)^{\ell+\ell'}\over\Gamma\big(-\ell+\nu+{1\over2}\big)
\Gamma\big(-\ell'-\nu+{1\over2}\big)}
\int_0^\infty W_{-\ell,\nu}(y)
W_{-\ell',\nu}(y){dy\over y},&(27.4)
}
$$
where the inner product is taken in $L^2(\B{R}^\times,
d^\times\!/\pi)$ and
we have used the fact that
$W_{\ell,\nu}(y)$ is real, as $(16.3)$ implies. Replacing
$\ell,\ell'$ by $\alpha,\beta\in\B{C}$ 
and applying $(16.8)$, the right side equals
$$
\eqalignno{
{1\over\pi(\alpha-\beta)\sin(2\pi\nu)}
\Big\{&\sin\big(\pi(\txt{1\over2}-\nu+\alpha)\big)
\sin\big(\pi(\txt{1\over2}+\nu+\beta)\big)\cr
-&\,\sin\big(\pi(\txt{1\over2}+\nu+\alpha)\big)
\sin\big(\pi(\txt{1\over2}-\nu+\beta)\big)\Big\}.&(27.5)
}
$$
On the assumption that $\alpha,\beta$ are
unequal integers, this
vanishes; and when $\ell=\ell'$, we take a
limit in $(27.5)$, getting 
$\Vert{\e{K}\phi_\ell(\cdot,\nu)}\Vert=1$.
\par
As for the completeness assertion in $(27.2)$, let
$g$ be a smooth function, compactly supported in
 ${\Bbb R}^\times$, which is orthogonal
to all $\e{K}\phi_\ell(\cdot,\nu)$, $\ell\in
\B{Z}$, in the space $L^2\big(\B{R}^\times,
d^\times\!/\pi\big)$. We apply integration by parts to
$(15.2)$, for $\r{g}=\r{a}[|u|]$, 
$\delta=\sgn(u)$, and take complex conjugate. We multiply
both sides by the factor $g(u)$
and integrate over $\B{R}^\times$ against $d^\times\! u$,
obtaining a double integral.
Because of absolute convergence we may exchange the
order of integrals; and in the new outer integral we
undo the integration by parts. In this way, we have
$$
0=\int_{-\infty}^\infty
{1\over(\xi^2+1)^{{1/2}-\nu}}
\left({\xi-i\over\xi+i}\right)^\ell\int_{-\infty}^\infty 
g(u) |u|^{-{1/2}+\nu}e(-u\xi)du\,d\xi. \eqno(27.6)
$$ 
Then we invoke that the set
$\left\{\big((\xi-i)/(\xi+i)\big)^\ell:
\,\ell\in{\Bbb Z}\right\}$ is a complete orthonormal
system of the space $L^2\!\left({\Bbb R},
(\pi(\xi^2+1))^{-1}d\xi\right)$, as
can be readily seen by the change of
variable $\xi\mapsto\tan\vartheta$.
Hence the Fourier transform of 
$g(u)|u|^{-{1/2}+\nu}$ vanishes, which ends the proof
of $(27.2)$.
\par
We next proceed to the proof of $(27.3)$.
By $(15.2)$ and $(16.1)$, the function
$W_{\ell,l-{1\over2}}(u)$, for $\ell\ge l$, equals
$u^l\exp(-u/2)$ multiplied by a polynomial factor
of degree $\ell-l$; one may apply residue calculus to 
$(15.2)$ under the present specification. In particular,
although the condition on $\Re\nu$ does not
hold, the argument $(16.10)$ extends nevertheless
to the product $W_{\ell,l-{1\over2}}(u)
W_{\ell',l-{1\over2}}(u)$, $\ell,\,\ell'\in\B{Z}$;
and we get the orthogonality for $\ell\ne\ell'$. On
the other hand, dealing with the case $\ell=\ell'$, 
we argue as follows: We let
$\omega^\pm$ be as in $(16.7)$. Then we have,
for $\ell>l$,
$$
\eqalignno{
\int_0^\infty \Big(W_{\ell,l-{1\over2}}(y)\Big)^2
{dy\over y}&\,=-\int_0^\infty \big[\omega^++(\ell-1)\big]
W_{\ell-1,l-{1\over2}}(y)
W_{\ell,l-{1\over2}}(y){dy\over y}\cr
&\,=\int_0^\infty W_{\ell-1,l-{1\over2}}(y)
\big[\omega^--\ell\big]W_{\ell,l-{1\over2}}(y){dy\over y}\cr
&\,=(\ell-l)(\ell+l-1)
\int_0^\infty \Big(W_{\ell-1,l-{1\over2}}(y)\Big)^2
{dy\over y}\cr
&\,=\Gamma(\ell-l+1)\Gamma(\ell+l).&(27.7)
}
$$
We have applied integration by parts as well as
the orthogonality of $W_{\ell-1,l-{1\over2}}$ and 
$W_{\ell,l-{1\over2}}$; and the last line is due to
$W_{l,l-{1\over2}}(y)
= y^l\exp\big(-{1\over2}y\big)$ as is implied by $(15.4)$
and $(16.1)$.
We have verified the orthonormality assertion in $(27.3)$.
On the other hand, let $g$ be a smooth 
function compactly supported in
$\B{R}^\times_+$ which is orthogonal to all 
$\e{K}\phi_\ell(\cdot,l-{1\over2})$,
$\ell\ge l$, in the space $L^2\big(\B{R}_+^\times,
d^\times\!/\pi\big)$. Then, because of the construction
of $W_{\ell,l-{1\over2}}(u)$ mentioned above, we have
$$
\int_0^\infty g(u)\exp(-2\pi u)u^{\ell-1}du=0,
\quad \ell\ge l.\eqno(27.8)
$$ 
Hence the Fourier transform of $g(u)\exp(-2\pi u)u^{l-1}$
vanishes, as it follows via the Taylor expansion of the factor
$\exp(i\xi u)$ to be multiplied. We end the proof of $(27.3)$.
\medskip
\noindent
{\csc Notes:} The argument of this section is
taken from [23, Part XII] and [5]. 
The orthonormality assertion in $(27.2)$ can
be proved in a smarter way as follows:
By $(14.2)$ we have, for $\phi(\r{g})=y^{1/2+\nu}
\Phi(\theta)$, $\nu\in i\B{R}$,
$$
\e{K}\phi(u)=|u|^{1/2-\nu}\int_{-\infty}^\infty
{e(-u\xi)\over(1+\xi^2)^{1/2+\nu}}
\Phi(\varkappa(\xi))d\xi,\quad
\exp(2i\varkappa(\xi))={\xi-i\over\xi+i}.\eqno(27.9)
$$
The Parseval formula for Fourier integrals gives
$$
\big\langle\e{K}\phi_1,\e{K}\phi_2\big\rangle
=\int_{-\infty}^\infty
{\Phi_1(\varkappa(\xi))\overline{\Phi_2(\varkappa(\xi))}
\over1+\xi^2}{d\xi\over\pi}=
\int_{\B{R}/\pi\B{Z}}\Phi_1(\theta)
\overline{\Phi_2(\theta)}\,
{d\theta\over\pi},\eqno(27.10)
$$
with an obvious correspondence,
where the inner product is taken 
in $L^2(\B{R}^\times,d^\times\!/\pi)$, and
the change of variable $\xi\mapsto\tan\theta$ has
been applied. This yields the assertion. However,
the argument does not seem to readily extend to the
discrete series, which is the reason why we have employed
an argument that might appear verbose. Naturally,
one may try to exploit the action such as $(18.5)$ 
of Maass operators in conjunction with $(14.3)$,
which is but essentially the same as the use of $(16.7)$,
duly indicated in the relevant notes. The argument of the 
present section
extends to $\r{PSL}(2,\B{C})$; see [24].
\medskip
\noindent
{\bf 28. Representations realised.} Now, we recapitulate
the most salient points of our discussion 
in the last two sections: We pick up a $V$
in the unitary principal series; thus $\nu_V\in i\B{R}$. 
The assertion $(27.2)$ means that
we have the unitary and surjective map
$$
{\e{L}: V\mapsto L^2\big(\B{R}^\times,d^\times\!/\pi\big),}
\atop{\e{L}\lambda_V^{(\ell)}(u)
=\e{K}\phi_\ell\big(u,\nu_V\big).}\eqno(28.1)
$$
Combined with $(26.5)$, this yields a realisation of
the representation $V$:
$$
{\hbox{The map $r^V: \r{h}\mapsto \e{L}
r_\r{h}\e{L}^{-1}$}\atop 
\hbox{is a unitary representation of $\r{G}$ 
on $L^2\big(\B{R}^\times,d^\times\!/\pi\big)$,}}
\atop\hbox{which is equivalent to the representation
$V$.}\eqno(28.2)
$$
That is, any right action of $\r{G}$ 
in the space $V$ is realised faithfully in
$L^2\big(\B{R}^\times,d^\times\!/\pi\big)$. 
We have seen at $(26.4)$ that inside $V^\infty$
the mode of the translation $\lambda\mapsto
r_\r{h}\lambda$ is exactly the same as that
of $\phi\mapsto r_\r{h}\phi$; namely, $\e{L}r_\r{h}\lambda
=\e{K}r_\r{h}\phi$ over $V^\infty$, with $\phi$ as in
$(26.2)$. In other words, we have 
$r^V_\r{h}\e{K}\phi=\e{K}r_\r{h}\phi$. Since
$\e{K}r_{\r{n}[\alpha]\r{a}[\beta]}\phi
(u)=e(\alpha u)\e{K}\phi(\beta u)$ and
$\big\{\e{K}\phi\big\}$ is dense in
$L^2\big(\B{R}^\times,d^\times\!/\pi\big)$, we have
$$
r^V_{\r{n}[\alpha]\r{a}[\beta]}f(u)=e(\alpha u)
f(\beta u),\quad \forall 
f\in L^2\big(\B{R}^\times,d^\times\!/\pi\big).\eqno(28.3)
$$
\par
For a $V$ in the holomorphic discrete series, we need only 
to replace $L^2(\B{R}^\times, d^\times\!/\pi)$ and
$\e{K}\phi_\ell(u,\nu_V)$ in $(28.1)$--$(28.3)$ by 
$L^2\big(\B{R}_+^\times, d^\times\!/\pi\big)$ and
$\pi^{1/2-l}\big(\Gamma(\ell+l)/
\Gamma(\ell-l+1)\big)^{1/2}\e{K}\phi_\ell(u,\nu_V)$,
$\ell\ge l$, respectively. The antiholomorphic discrete series
is analogous; it suffices to apply
the involution $J$ defined by $(13.2)$.
\medskip
\noindent
{\csc Notes:} In literature, the fact $(28.2)$ 
is often termed the Kirillov model. 
This was recently employed by the present author 
[23, Parts XIV and XV] in his
resolution of Selberg's decades-old problem [33] to find
a complete spectral decomposition, within 
the structure of $L^2(\varGamma
\backslash\r{G})$, of the
shifted convolution of Fourier coefficients of cusp forms,
i.e., sums
$$
\sum_{n=1}^\infty
\varrho_V(n+f)\overline{\varrho_V(n)}\,W(n/f),\eqno(28.4)
$$
with generic cuspidal representations $V$ 
as is defined by $(24.3)$, which is obviously
an analogue of $(25.1)$. In the solution
the use of the right action of $\r{K}$ is
essential; that is, the spectral theory of cusp forms of
weight zero is inadequate, and representation theory comes into
play in an indispensable fashion. See also 
[27, Vol.\ 2, Section 7.2].
\medskip
\noindent
{\bf 29. Bessel functions of representation.} We point out an 
obvious incompleteness of the formula $(28.3)$: 
It lacks the description of the right action of $\r{K}$.
This is because of the fact that when considering 
$r^V_{\r{k}[\tau]}f$ with the general
combination of $\r{k}[\tau]$ and $f$ 
we need first to spectrally
decompose $f$ in terms of the system 
$\big\{\e{K}\phi_\ell(u,\nu_V)\big\}$ 
if we follow the argument so far developed. 
Thus one may
ponder whether it is possible or not 
to avoid this tedious procedure. It should be observed here that
we actually do not need to consider all 
of $\{\r{k}[\tau]\}$, since we have the Bruhat decomposition
$(14.4)$. Namely, what really matters is
the action of the Weyl element 
$\r{w}=\r{k}\big[{1\over2}\pi\big]$. 
We need to express $r^V_\r{w}f$
in terms of $f$. The answer is 
an integral transform, and its kernel is the
Bessel function of representation in the context
of $\r{PSL}(2,\B{R})$, a fundamental concept
due to Gel'fand--Graev--Pyatetskii-Shapiro [8]:
In a slightly
generalised form the kernel is defined by
$$
j_\nu(\lambda)=\pi{\sqrt{|\lambda|}\over\sin\pi\nu}
\left(J_{-2\nu}^{\sgn(\lambda)}
\big(4\pi\sqrt{|\lambda|}\big)-J^{\sgn(\lambda)}_{2\nu}
\big(4\pi\sqrt{|\lambda|}\big)\right),
\quad\lambda\in\B{R}^\times,\;
|\Re\nu|<{1\over2},\eqno(29.1)
$$ 
where $J^+_\nu=J_\nu$ and $J^-_\nu=I_\nu$ 
with the usual notation for Bessel functions. 
We have that
$$
{\hbox{if $f\in L^2\big(\B{R}^\times, 
d^\times\!/\pi\big)$ is compactly 
supported in $\B{R}^\times$, then}\atop
\hbox{$\displaystyle{r^V_\r{w}f(u)=\int_{\B{R}^\times}
j_{\nu_V}(u\lambda)f(\lambda)
d^\times\!\lambda\,}$ in
$L^2\big(\B{R}^\times, 
d^\times\!/\pi\big)$,}}\atop
\hbox{for any $V$ in the unitary principal series.}
\eqno(29.2)
$$
It is meant that the integral transform maps a dense subset of
$L^2\big(\B{R}^\times, d^\times\!/\pi\big)$
into the space unitarily. Together with $(14.4)$ and 
$(28.3)$, this describes explicitly the
action of $\r{G}$ over $L^2(\B{R}^\times, 
d^\times\!/\pi)$ via the invariant subspace $V$ of
$L^2(\varGamma\backslash\r{G})$.
\par
For the proof of $(29.2)$ we consider first the Mellin transform
$$
\Gamma_\ell(s,\nu)=\int_0^\infty
\e{A}^+\!\phi_\ell(\r{a}[y],\nu)
y^{s-{3/2}}dy.\eqno(29.3)
$$ 
We shall show that $\Gamma_\ell(s,\nu)$ exists in
the domain $|\Re s-{1\over2}|<\Re\nu+{1\over2}$ and
satisfies the local functional equation of 
Jacquet--Langlands [14, p.\ 196]
$$
\eqalignno{ 
(-1)^\ell&\Gamma_\ell(s,\nu) 
=2(2\pi)^{-2s}\Gamma(s+\nu)
\Gamma(s-\nu)\cr &\times\Big(\cos(\pi s)
\,\Gamma_\ell(1-s,\nu)
+\cos(\pi\nu)\,\Gamma_{-\ell}(1-s,\nu)\Big).&(29.4)
}
$$ 
To this end, in $(15.2)$, for $\r{g}=\r{a}[y]$, 
we shift the contour to
$\Im\xi=-{1\over2}$, and insert the result into $(29.3)$.
The double integral converges absolutely in the domain
$0<\Re\nu<\Re s$. After exchange, the inner integral
is seen to be $(2\pi i\xi)^{\nu-s}\Gamma(s-\nu)$,
$|\arg(i\xi)|<{1\over2}\pi$. We shift the $\xi$-contour 
back to the original, getting
$$
\eqalignno{
\Gamma_\ell(s,\nu)&=(2\pi)^{\nu-s}\Gamma(s-\nu)\cr
\times&\Big[\exp\left(\txt{1\over2}
\pi i(s-\nu)\right)\r{Y}_{-\ell}(s,\nu)
+\exp\left(-\txt{1\over2}\pi i(s-\nu)\right)
\r{Y}_\ell(s,\nu)\Big],\qquad&(29.5)\cr
\r{Y}_\ell(s,\nu)&=\int_0^\infty{\xi^{-s+\nu}
\over(\xi^2+1)^{\nu+{1/2}}}
\left({\xi-i\over\xi+i}\right)^\ell d\xi.&(29.6)
}
$$
This integral converges absolutely and uniformly
for $|\Re s-{1\over2}|<\Re\nu+{1\over2}$.
We apply the change
of variable $\xi\mapsto\xi^{-1}$, and find that
in the same domain
$$
(-1)^\ell\r{Y}_\ell(s,\nu)=\r{Y}_{-\ell}(1-s,\nu).\eqno(29.7)
$$ 
We then apply, to $(29.5)$, the transformations
$s\mapsto 1-s$,
$\ell\mapsto\pm\ell$ as well as $(29.7)$,
which yields, after an elimination,
$$
\eqalignno{
(-1)^\ell Y_\ell&(s,\nu)=(2\pi)^{-s-\nu}\Gamma(s+\nu)\cr
&\times\Big(\Gamma_\ell(1-s,\nu)
\exp\big(\txt{1\over2}(s+\nu)\big)+
\Gamma_{-\ell}(1-s,\nu)
\exp\big(-\txt{1\over2}(s+\nu)\big)\Big).&(29.8)
}
$$
Inserting this back to $(29.5)$ we obtain $(29.4)$.
\par
On the other hand, we have, for
 $|\Re\nu|<\Re s<{1\over4}$,
$$
\int_0^\infty j_\nu(\lambda)\lambda^{s-1/2}d^\times\!\lambda
=2(2\pi)^{-2s}
\cos(\pi s)\Gamma\big(s+\nu\big)
\Gamma\big(s-\nu\big),
\eqno(29.9)
$$
and, for $|\Re\nu|<\Re s$,
$$
\int_{-\infty}^0 j_\nu(\lambda)|\lambda|^{s-1/2}
d^\times\!\lambda=2(2\pi)^{-2s}
\cos(\pi\nu)\Gamma\big(s+\nu\big)
\Gamma\big(s-\nu\big),\eqno(29.10)
$$
which are consequences of the well-known integral
formulas
$$
{{\displaystyle\int_0^\infty J_\nu(y)y^{s-1}dy
=2^{s-1}{\Gamma\big(\txt{1\over2}(s+\nu)\big)\over
\Gamma\big(1-\txt{1\over2}(s-\nu)\big)},\quad
-\Re\nu<\Re s<{1\over2}\,,}\atop
{\displaystyle\int_0^\infty K_{\nu}(y)y^{s-1}dy=2^{s-2}
\Gamma\left(\txt{1\over2}(s+\nu)\right)
\Gamma\left(\txt{1\over2}(s-\nu)\right),\quad
|\Re\nu|<\Re s.}}\eqno(29.11)
$$
By  $(29.9)$--$(29.10)$, we may rewrite $(29.4)$ as
$$
(-1)^\ell\Gamma_\ell(s,\nu)
=\int_{{\Bbb R}^\times}j_\nu(\lambda)|\lambda|^{s-{1/2}}
\Gamma_{\sgn(\lambda)\ell}(1-s,\nu)
d^\times\!\lambda,\eqno(29.12)
$$
although the range of $(\nu,s)$ is to be restricted to have $(29.9)$.
We replace $\ell$ by $\sgn(u)\ell$, multiply both sides
by the factor $|u|^{{1/2}-s}/2\pi i\ne0$, and integrate 
along $\Re s=\beta$, $|\Re\nu|<\beta<{1\over4}$.
We get
$$
\eqalignno{
&{(-1)^\ell\over2\pi i}\int_{(\beta)}
\Gamma_{\sgn(u)\ell}(s,\nu)|u|^{{1/2}-s}ds\cr
&=\int_{{\Bbb R}^\times}j_\nu(\lambda)
\left\{{1\over2\pi i}\int_{(\beta)}\Gamma_{\sgn(\lambda v)\ell}
(1-s,\nu)|\lambda/u|^{s-{1/2}}ds\right\}
d^\times\!\lambda. &(29.13)
}
$$
This exchange is legitimate, since
the function $\Gamma_\ell (s,\nu)$ is of rapid decay, as
can be seen by turning the line of integration in $(29.6)$
through a small angle round the origin. Being a Mellin inversion of
$(29.3)$, the left side equals
$(-1)^\ell\e{A}^+\!\phi_{\sgn(u)\ell}(\r{a}[|u|],\nu)
={\eusm K}r_\r{w}\phi_\ell(u,\nu)$, while the
inner integral is
$\e{A}^+\!\phi_{\sgn(\lambda u)\ell}(\r{a}[|\lambda/u|],\nu)
={\eusm K}\phi_\ell(\lambda/u,\nu)$. Namely, we have
obtained the following point-wise
identity, but only for $|\Re\nu|<{1\over4}$:
$$
\e{K}r_\r{w}\phi_\ell(u,\nu)=\int_{\B{R}^\times}
j_\nu(u\lambda)\e{K}
\phi_\ell(\lambda,\nu)d^\times\!\lambda,
\quad |\Re\nu|<{1\over2}.\eqno(29.14)
$$
The extension of the range of $\nu$ can be
attained by $(15.5)$ and analytic continuation.
\par
Now, let $f$ be as in $(29.2)$.
We consider the double integral
$$
\int_{\B{R}^\times}\left(\int_{\B{R}^\times}j_{\nu_V}(u\lambda)
f(\lambda)d^\times\!\lambda\right)\e{K}\phi_\ell(u,\nu_V)
{d^\times\!u\over\pi}.\eqno(29.15)
$$
Invoking $(15.5)$ again, we see that this converges
absolutely, and find that $(29.14)$ implies
that it equals
$$
\big\langle{f,\e{K}r_\r{w}\phi_\ell(\cdot,\nu_V)}\big\rangle
=\big\langle{f,r^V_\r{w}\e{K}\phi_\ell(\cdot,\nu_V)}\big\rangle
=\big\langle{r_\r{w}^Vf,\e{K}\phi_\ell(\cdot,\nu_V)}\big\rangle
\eqno(29.16)
$$
in $L^2(\B{R}^\times,d^\times\!/\pi)$, as $r_\r{w}^V$ is
an involution.
Via $(27.2)$ one might conclude that 
$(29.2)$ has been confirmed. 
However, it remains for us to show that
$$
\hbox{$\displaystyle
F(u)=\int_{\B{R}^\times}j_{\nu_V}(u\lambda)
f(\lambda)d^\times\!\lambda\,$ is in  
$L^2(\B{R}^\times,d^\times\!/\pi).$
}\eqno(29.17)
$$
When $|u|\le1$, the definition $(29.1)$ implies that
$j_{\nu_V}(u\lambda)\ll \sqrt{|u\lambda|}\,$; and
$F(u)\ll\sqrt{|u|}$, which means that
the integral of $|F(u)|^2$ against $d^\times\!u$ 
over $|u|\le1$ is finite. 
When $|u|\ge1$ and $u\lambda<0$,
the asymptotic expansion for $K$-Bessel functions gives
$j_{\nu_V}(u\lambda)\ll \exp(-\sqrt{|u\lambda|})$; that is,
the corresponding part of $F(u)$ can be ignored. 
Hence it suffices to show that 
$$
\int_{U_1}^{U_2}\left|\int_0^\infty j_{\nu_V}(u\lambda)
f(\lambda)d^\times\!\lambda\right|^2d^\times\!u\ll1
\;\hbox{uniformly for $1\le U_1<U_2$,}
\eqno(29.18)
$$
since the part with $u,\lambda<0$ is analogous.
We then invoke the asymptotic expansion for $J$-Bessel functions;
and the discussion is reduced to that on the expression
$$
\eqalignno{
&\int_{U_1}^{U_2}\left|\int_0^\infty(u\lambda)^{1/4}
\exp\big(4\pi i(u\lambda)^{1/2}\big)
f(\lambda)d^\times\!\lambda\right|^2d^\times\!u\cr
=&\,8\int_{\sqrt{U_1}}^{\sqrt{U_2}}
\left|\int_0^\infty\mu^{-1/2}
\exp\big(4\pi i v\mu\big)
f(\mu^2)d\mu\right|^2dv,&(29.19)
}
$$
in which the
change of variables $(u,\lambda)\mapsto(v^2,\mu^2)$
has been applied. Hence, by the Parseval 
formula for Fourier integrals, we
obtain $(29.17)$. This ends the proof of $(29.2)$. 
\par
The discussion on the discrete series is skipped,
as it is fairly analogous to the above. We remark only that
for any $l\in\B{N}$
$$
j_{l-{1\over2}}(u)=\cases{\hfil 0& $u<0$,\cr
2\pi(-1)^l\sqrt{u}\, J_{2l-1}(4\pi\sqrt{u})& $u>0$.
}\eqno(29.17)
$$
In passing, we remark that the complementary series of
irreducible representations can be discussed in much the
same way; then the range of $\nu$ in $(29.14)$
becomes relevant, 
although this is immaterial in our present context.
\medskip
\noindent
{\csc Notes:} As to the formulas $(29.11)$  
as well as the asymptotic
expansion of Bessel functions, see
Watson [37]; a concise proof of these
can be found in [27, Vol.\ 2, Section 1.4]. 
In conjunction with the Kirillov map, the
Bessel function of representation $(29.1)$ 
has played a fundamental
r\^ole in the modern developments of analytic number theory.
As Cogdell and Piatetskii-Shapiro pointed out in their inspiring
monograph [6], it arises in the spectral decomposition
of Poincar\'e series on $\r{G}$ in general, typical instances
of which occur in the theory of
sums of Kloosterman sums 
due to Bruggeman [2] and Kuznetsov [17]
and in the theory of the fourth power
moment of the Riemann zeta-function due to the
present author [22]. See also [5], 
[23, Parts XIV and XV] and [26], for instance. 
The above proof of the fundamental transformation
formula in $(29.2)$ is an adaptation of 
a part of [23, Part XII] where
smooth $f$, i.e., those $\e{K}\phi$ with $\phi$ as in $(26.2)$,
is in fact dealt with; and this time we have
extended the assertion to compactly supported vectors by means of
the augmentation $(29.15)$--$(29.19)$. Naturally, $(15.5)$
and $(29.14)$ readily imply the identity
in $(29.2)$ for any smooth $f$.
The above proof of the very basic functional equation
$(29.4)$ comes also from [23, Part XII].
An alternative
and independent proof of the identity in $(29.2)$, 
but only for smooth $f$, is given 
by Baruch and Mao [1, Section 6 and Appendix 2], which is in fact a
verification of Vilenkin's claim made at 
the end of Chapter VII of the older
edition of [36]; see also [36, p.\ 454]. It should
be stressed that for any smooth $f$ the identity in $(29.2)$
holds point-wise.
We add that the statement $(29.2)$ can readily be extended to
any vector $f$ in terms of the mean convergence; see
the first line of $(30.4)$ below.
Further, we remark that the argument of this
section extends to $\r{PSL}(2,\B{C})$; necessary means
may be found in [24].
Incidentally, it will be worth remarking that
the discussion in [6] is incomplete in the sense that it lacks
an inversion procedure which is necessary in
stating two versions of the spectral expansion of
sums of Kloosterman sums, due originally
to Kuznetsov [17]. This task of reversing Bessel 
transforms containing the
kernel $(29.1)$ is resolved in [25]; see also 
Bruggeman [3] for an $L^2$-approach. 
The notes to Section 34 below contains another 
aspect of $j_\nu$.
\medskip
\noindent
{\bf 30. Irreducibility.} We now prove that
any subspace $V$ in the unitary principal
series is irreducible. Because of $(28.2)$ it suffices
to show that $r^V$ is an irreducible representation.
Thus, let $U_1$ be an invariant subspace of
$L^2\big({\Bbb R}^\times,d^\times\!/\pi\big)$ and $U_2$ be
its orthogonal complement. For each $f_1\in U_1$ 
we have $r^V_{\r{n}[\alpha]\r{a}[\beta]}f_1\in U_1$,
$\forall\alpha\in\B{R}$ and $\forall\beta>0$; 
and  by $(28.3)$ the Fourier
transform of $f_1(\beta u)\overline{f_2(u)}/|u|$
vanishes identically for any $f_2\in U_2$; that is, 
$$
\int_{\B{R}^\times} |f_1(\beta u)f_2(u)|d^\times\!u=0.
\eqno(30.1)
$$
Integrating this over the positive real axis 
against the measure $d^\times\!\beta$, we find
via Fubini's theorem  that
$$
\eqalignno{
&\left(\int_0^\infty|f_1(u)|d^\times\!u\right)
\left(\int_0^\infty|f_2(u)|d^\times\!u\right)=0,&(30.2)\cr
&\left(\int_{-\infty}^0|f_1(u)|d^\times\!u\right)
\left(\int_{-\infty}^0|f_2(u)|d^\times\!u\right)=0.&(30.3)
}
$$
We are, however, unable to assert that any 
combination of $(30.2)$ and
$(30.3)$ implies that one of $f_1$ and $f_2$ vanishes
almost everywhere in $\B{R}^\times$. Overcoming this
difficulty, we argue as follows: According to $(30.3)$, one of the
sets $\{u<0: f_1(u)\ne0\}$ and
$\{u<0: f_2(u)\ne0\}$ has Lebesgue measure zero.
We assume that the former holds
for all elements in $U_1$;
otherwise we may exchange $U_1$ and $U_2$. 
We then apply the assertion $(29.2)$ to $f_L$
the restriction of an arbitrary $f\in U_1$ to $[1/L,L]$, and have
$$
\displaystyle{\lim_{L\to\infty}
\int_{\B{R}^\times}\left|r^V_\r{w}f(u)-\int_{\B{R}^\times}
j_{\nu_V}(u\lambda)f_L(\lambda)d^\times\!\lambda\right|^2
{d^\times\! u}=0,}\atop
\hbox{$r^V_\r{w}f(u)=0$ for almost all $u<0$,}\eqno(30.4)
$$
which implies that
$$
\lim_{L\to\infty}\int_{-\infty}^0\left|\int_0^\infty
j_{\nu_V}(u\lambda)f_L(\lambda)d^\times\!\lambda\right|^2
{d^\times\! u}=0,\eqno(30.5)
$$
The assertion $(30.4)$ transfers the
fact on the negative real axis to $(30.5)$ which
concerns the values of $f$ on the positive real axis.
In $(30.5)$, $f_L(\lambda)$ has the Mellin transform $f_L^*(s)$
for any $s\in \B{C}$; and $(29.10)$ gives that of
$j_{\nu_V}(u\lambda)\lambda^{-1}$
for any $s$ with $\Re s>{1\over2}$. Applying the Mellin
inversion to the latter on $\Re s=1$ and an exchange to the
resulting double integral,
we see that $(30.5)$ is equivalent to
$$
\lim_{L\to\infty}\int_\B{R}\left|\int_\B{R}\Gamma
\big(\txt{1\over2}-it+\nu_V\big)
\Gamma\big(\txt{1\over2}-it-\nu_V\big)f_L^*(it)
e^{-i\xi t}dt\right|^2d\xi=0.
\eqno(30.6)
$$
Thus we have
$$
\lim_{L\to\infty}\int_\B{R}\left|\Gamma
\big(\txt{1\over2}-it+\nu_V\big)
\Gamma\big(\txt{1\over2}-it-\nu_V\big)f_L^*(it)
\right|^2dt=0.\eqno(30.7)
$$
On the other hand, according to the $L^2$-theory of
Mellin transforms, $f^*_L(it)$ converges in the mean to an
$f^*(it)\in L^2(i\B{R})$; and  
$\Vert{f^*}\Vert=\Vert{f}\Vert$ with an obvious specification of
the norms. In particular, we have
$$
\eqalignno{
\lim_{L\to\infty}&\int_\B{R}\left|\Gamma
\big(\txt{1\over2}-it+\nu_V\big)
\Gamma\big(\txt{1\over2}-it-\nu_V\big)
(f^*(it)-f_L^*(it))\right|^2dt\cr
&\ll\lim_{L\to\infty}\int_{\B{R}}|f^*(it)-f_L^*(it)|^2dt=0.
&(30.8)
}
$$
It follows that
$$
\int_\B{R}\left|\Gamma
\big(\txt{1\over2}-it+\nu_V\big)
\Gamma\big(\txt{1\over2}-it-\nu_V\big)f^*(it)\right|^2
dt=0,\eqno(30.9)
$$
which implies $\Vert{f^*}\Vert=0$. Hence, we have
$\Vert f\Vert =0$. This ends the treatment of the
unitary principal series. 
\par
The discussion of the
discrete series is much simpler; the assertion corresponding
$(30.2)$--$(30.3)$ suffices. As before we remark that the
above extends to the complementary series,
although this is irrelevant in our present context.
\medskip
\noindent
{\csc Notes:} The irreducibility of $V$'s is
usually confirmed via the Lie algebra $\f{g}$. Our
argument may appear to be unconventional. However,
the discussion based on explicit group actions $(28.3)$ and
$(29.2)$ should also be worth reporting. 
As for the $L^2$-theory of
Mellin transforms, see Titchmarsh [35, Section 3.17].
\medskip
\noindent
{\bf 31. Functional equations.} This is kind of
a rumination. The r\^ole
played either visibly or invisibly
by the Weyl element in functional equations for
automorphic $L$-functions 
is condensed in $(29.4)$; in other words those equations can
be all traced back to respective versions of $(29.4)$. 
Although only remotely relevant to the
motivation of the present article,
it might be worthwhile to ponder upon the nature of
functional equations for general $L$-functions. 
In fact,
we have set out already the following view 
point in [27, Vol.\ 2],  starting with the Poisson and 
the Vorono{\"\i} sum formulas, which are
equivalent to the functional equations for
the Riemann zeta-function and for the
product of two zeta-values, respectively: After
observing that these classical 
sum formulas yield a variety of significantly non-trivial
bounds, that is, massive cancellations among number theoretical
terms, we assert that
an art of counting, which is at the core of mathematics,
is to surmise that those discrete objects
to be counted are floating on waves, i.e., harmonic
extensions of discrete existence, whence
one may see the very possibility of
applying analysis of continuous objects to
discrete objects. Indeed, the zeta and $L$-functions are
additive collections of the zeta-waves $\{n^{it}:\,
n\in\B{N}\}$. What counts most is thus to try to 
express counting functions as precisely as possible in
terms of these analytic functions. However, 
we face a dilemma. Counting functions are 
genetically discontinuous, and we need infinitely 
many waves to exactly describe them, that is, limiting
procedure is inevitable. There are thus gaps
in between discrete and continuous existence; and the
speculation on the best possible 
bounding of the gaps leads us to
conjectures. Most tantalising 
problems in analytic number
theory such as the distribution of prime numbers concern certain
error terms in approximations by means of
continuous main terms to counting
functions. With this, one may wonder how to
detect cancellations among waves. Thence we come
to the view that functional equations should be 
the evidence of the existence of cancellations, and they
must be unable to hold without
number theoretical peculiarities of those coefficients attached to
zeta-waves. 
The contrary is indeed hard to imagine. 
Therefore it is of the foremost importance for us to
know the origin of those functional equations. It is in
the action of the Weyl element whose existence
defines $\r{PSL}(2,\B{R})$ as a matrix
group. Namely, various cancellations among zeta-waves
are the work of the group structure of $\r{PSL}(2,\B{R})$.
We are led to the vision that
there should be unified approaches to quantitative problems
in the theory of $L$-functions. See [15] and [23, Part XIV and XV]
for investigations along this line of thoughts.
\medskip
\noindent
{\csc Notes:} The statistical detection of cancellations among Fourier
coefficients of cusp forms, which is well in the category of
the above thoughts, can sometimes 
yield assertions far deeper than those
resulting on the Ramanujan conjecture. 
A typical instance is the spectral analysis of
sums of Kloosterman sums which leads us 
to a region far beyond that the Weil bound does.
Here we see a structure similar to that witnessed in the
distribution of prime numbers via the method of the
large sieve. 
\medskip
\noindent
{\bf 32. A weight stratum.} The rest is to be
regarded as an appendix to the foregoing discussion, although
it is not directly related to automorphic representations
and moreover disproportionally long. 
\par
We shall stay on a single stratum rather than view the whole 
of all strata.
We shall indicate how to approach to the assertion of Section 22 
by extending the argument of [22, Chapter 1] to 
$L_\ell^2(\varGamma\backslash\r{G})$; however, the spectral 
decomposition $(22.7)$--$(22.11)$ proper
is not treated, because of an obvious redundancy.
The task is, in essence, to find the Green function for 
$\Omega_\ell$ and investigate its analytic nature, which will occupy
the next five sections.
Then in Section 38 we shall turn to the Selberg trace formula.
We shall restrict ourselves to
the situation of weight zero, i.e., the strictly real analytic
environment, which will not cause any significant loss 
of generality, as is to be explained.
\par
To begin with, we shall view the invariance 
asserted in the first section from the Lie group $\r{G}$. Thus, we
shall show first that the Maass operators $\b{e}^\pm$ can be
regarded as extensions of the hyperbolic
outer-normal differential: We let $f$ be a smooth
function on $\B{H}^2$ and put $F(\r{n}[x]\r{a}[y]
\r{k}[\theta])=f(z)e^{2\ell i\theta}$ and
$l_\r{h}F(\r{g})=f^\r{h}(z)e^{2\ell i\theta}$, with $z=x+iy$ and
$\r{h}\in\r{G}$. We have $\b{e}^\pm F=2(\r{e}^\pm_\ell f)
e^{2(\ell\pm1) i\theta}$,
with
$$
\r{e}^\pm_\ell=\pm iy{\partial\over\partial x}+
y{\partial\over\partial y}\pm\ell.\eqno(32.1)
$$
The left invariance of $\b{e}^\pm$ implies that
$$
2(\r{e}_\ell^\pm f^\r{h})e^{2(\ell\pm1)i\theta}
=\b{e}^\pm(f^\r{h}(z)e^{2\ell i\theta})=\b{e}^\pm l_\r{h}F=
l_\r{h}\b{e}^\pm F=2(\r{e}_\ell^\pm f)^\r{h}
e^{2(\ell\pm1) i\theta},\eqno(32.2)
$$
which in view of $(3.4)$ is equivalent to 
$$
\r{e}^\pm_\ell\left(f(\r{h}(z))\left({\jmath(\r{h},z)
\over\jmath(\r{h},\overline{z})}\right)^{-\ell}\right)=
(\r{e}^\pm_\ell f)(\r{h}(z))\left({\jmath(\r{h},z)
\over\jmath(\r{h},\overline{z})}\right)^{-(\ell\pm1)}.
\eqno(32.3)
$$
When $\ell=0$, this is the same as
$$
y\left(\pm i{\partial\over\partial x}+
{\partial\over\partial y}\right)f(w){dx\pm idy\over|dz|}
=v\left(\pm i{\partial\over\partial u}+
{\partial\over\partial v}\right)f(w)
{du\pm idv\over|dw|},\eqno(32.4)
$$
since $\jmath(\r{h},z)/\jmath(\r{h},
\overline{z})=(dw/dz)/|dw/dz|$, for $\r{h}(z)=u+iv=w$.
A rearrangement gives $(1.7)$. 
\par
Analogously the left invariance of $\Omega$ implies that
$$
\Omega_\ell \left(
f(\r{h}(z))\left({\jmath(\r{h},z)\over
\jmath(\r{h},\overline{z})}\right)^{-\ell}\right)
=(\Omega_\ell f)(\r{h}(z))
\left({\jmath(\r{h},z)\over
\jmath(\r{h},\overline{z})}\right)^{-\ell},
\quad \forall\r{h}\in\r{G},\eqno(32.5)
$$
which is an extension of $(1.8)$. Alternatively one may use
the relation
$$
\Omega_\ell=-\r{e}_{\ell-1}^+\r{e}_\ell^--\ell(\ell-1),
\eqno(32.6)
$$
which is an interpretation of the second line of $(8.3)$.
This yields also, for any smooth functions 
$f,\,g$ on $\B{H}^2$,
$$
\int_D\big((\Omega_\ell+\ell^2-\ell)f\big)\cdot g\,d\mu(z)
-\int_D (\r{e}_\ell^-f)\cdot(\r{e}_{-\ell}^+g)d\mu(z)
=\int_{\partial D}(\r{e}_\ell^-f)\cdot g{\overline{dz}\over y},
\eqno(32.7)
$$
where $D$ and $\partial D$ are as in $(1.9)$. 
It suffices to apply integration by parts once; indeed, such an effect is
a merit of introducing the Maass operators.
The assertion $(1.9)$ follows from $(32.7)$ for $\ell=0$.
\par
The relations $(32.3)$ and $(32.5)$, 
for $\r{h}\in\varGamma$, mean that
$$
\hbox{$\r{e}^\pm_\ell$ maps 
$B^\infty_\ell(\varGamma\backslash\B{H}^2)$ into
$ B^\infty_{\ell\pm1}(\varGamma\backslash\B{H}^2)$,}
\atop\hbox{$\Omega_\ell$ maps 
$B^\infty_\ell(\varGamma\backslash\B{H}^2)$ into itself,}
\eqno(32.8)
$$
where $B^\infty_\ell(\varGamma\backslash\B{H}^2)
\exp(2i\ell\theta)=B^\infty_\ell(\varGamma\backslash\r{G})$;
see $(13.1)$.
The second line is, however, trivial. Also we have the symmetry: 
$$
\big\langle{\Omega_\ell f_1,f_2}\big\rangle
=\big\langle{f_1,\Omega_\ell f_2}\big\rangle,
\quad  f_1,\,f_2\in 
B^\infty_\ell(\varGamma\backslash\B{H}^2),\eqno(32.9)
$$
where the inner product is taken in 
$L^2_\ell(\varGamma\backslash\B{H}^2)$.  This
is an immediate consequence of $(11.4)$. In the present
context, it should, however, be
more expedient to derive $(32.9)$ from $(32.7)$: 
We have
$$
\int_\e{F}\big((\Omega_\ell+\ell^2-\ell)f_1\big)
\cdot\overline{f_2}\,d\mu(z)
=\int_\e{F} (\r{e}_\ell^-f_1)
\cdot\overline{(\r{e}_\ell^-f_2)}\,d\mu(z),
\eqno(32.10)
$$
since
$$
\int_{\partial\e{F}}(\r{e}_\ell^-f_1)\cdot \overline{f_2}\,
{\overline{dz}\over y}=0.\eqno(32.11)
$$
In fact $(32.3)$ implies that $(\r{e}_\ell^-f_1)\cdot 
\overline{f_2}\cdot\overline{dz}/|dz|$ is a 
$\varGamma$-automorphic quantity; 
and thus each part of the integral $(32.11)$
corresponding to $\varGamma$-congruent
pairs of the sides of $\e{F}$
vanishes, as the point $z$ on the
contour proceeds along
the sides in the direction opposite to that implied by
the $\varGamma$-maps. We get $(32.9)$ from $(32.10)$.
In particular we have
$$
\big\langle{\Omega_\ell f,f}\big\rangle\ge -\ell(\ell-1)\Vert{f}
\Vert^2,\quad\forall f\in 
B^\infty_\ell(\varGamma\backslash\B{H}^2),\eqno(32.12)
$$
which is equivalent to the restriction on the range of summands
in $(22.8)$.
\medskip
\noindent
{\csc Notes:} The assertions of this section 
are stated in a generalised form by Maass [20,
Chapter 4]. Our discussion is to be compared with 
Roelcke [31, Teil I]. 
\medskip
\noindent
{\bf 33. Green's function to find.} Now, we 
suppose that $g_{\alpha,\ell}(z,w)$ stands for
the free-space Green function for the operator
$\Omega_\ell+\alpha(\alpha-1)$ acting 
over $\B{H}^2$; the term `free-space' is to
indicate that the action of $\varGamma$ is not taken into
account yet. This is to satisfy
$$
\big(\Omega_\ell+\alpha(\alpha-1)\big)_{\!z}\,
g_{\alpha,\ell}(z,w)=0,\quad z\ne w,\eqno(33.1)
$$
and have logarithmic singularities along the diagonal $z=w$. 
We shall assume that
$$
\Re\alpha>\ell+1,\quad \ell\ge0.
\eqno(33.2)
$$
It will turn out that the cases $\ell\le0$ and $\ell\ge0$
are essentially the same; see $(34.2)$ below.
The requirement about the singularities of $g_{\alpha,\ell}$ 
is suggested by the potential theory on the Euclidean plane.
As a matter of fact we may employ the definition
$(34.1)$ below in an a priori manner; however, it should be more
beneficial to know how to reach $(34.1)$. 
\par
We begin with the postulation that
there exists a function $p_{\alpha,\ell}$ on positive reals such that
$$
g_{\alpha,\ell}(z,w)=p_{\alpha,\ell}\big(\varrho(z,w)\big)
H_\ell(z,w),
\quad H_\ell(z,w)=\left({\overline{z}-w\over 
\overline{w}-z}\right)^\ell,\eqno(33.3)
$$
where $\varrho(z,w)$ is as in $(2.1)$. The factor $H_\ell(z,w)$ is 
related to the Cartan decomposition of $\r{G}$. It
might, in fact, be easier to work with the unit disk model than
with $\B{H}^2$ of the 2-dimensional hyperbolic geometry; that is,
discussion will become somewhat transparent if viewed via
the map $(2.3)$. 
Namely, $(33.3)$ is an expression of the natural requirement that 
$g_{\alpha,\ell}$ be a point-pair invariant $h(\r{g}_1^{-1}
\r{g}_2)$ if viewed from $\r{G}\times\r{G}$ 
with which the left and right actions of $\r{K}$
induce characters. More precisely, according to $(4.5)$ we should
have
$$
{g_{\alpha,\ell}(z,w)=h\big((\r{n}[x]\r{a}[y])^{-1}
\r{n}[u]\r{a}[v]\big),}\atop {h\big(\r{k}[\tau_1]\r{g}
\r{k}[\tau_2]\big)=\exp\big(-2\ell i(\tau_1+\tau_2)\big)
h(\r{g})},\eqno(33.4)
$$
where $z=x+iy,\, w=u+vi$.
Also, the construction $(8.3)$
of the Casimir operator in terms of the Maass operators could be
utilised in the computation below, which also yields simplifications.
We shall, nevertheless, work with $\B{H}^2$ directly.
\smallskip
We have, with $\varrho=\varrho(z,w)$,
$$
\matrix{
\varrho_x=(x-u)/(2yv),\hfill&
\varrho_y=\left(y^2-v^2-(x-u)^2\right)/(4y^2v),\hfill\cr\cr
\varrho_{xx}=1/(2yv),\hfill
&\varrho_{yy}=\left((x-u)^2+v^2\right)/(2y^3v),\hfill\cr
}\eqno(33.5)
$$ 
and, with $p=p_{\alpha,\ell}$,
$$
\eqalignno{
{\partial\over\partial x}g_{\alpha,\ell}(z,w)
&=p'(\varrho)\varrho_xH_\ell+\ell
p(\varrho)H_{\ell-1}{\overline{w}-z+\overline{z}-w
\over (\overline{w}-z)^2},\cr
{\partial\over\partial y}g_{\alpha,\ell}(z,w)
&=p'(\varrho)\varrho_yH_\ell-i\ell
p(\varrho)H_{\ell-1}{\overline{w}-z-\overline{z}+w
\over (\overline{w}-z)^2},\cr
\left({\partial\over\partial x}\right)^2
g_{\alpha,\ell}(z,w)&=p''(\varrho)\varrho_x^2
H_\ell+p'(\varrho)\varrho_{xx}H_\ell
+2\ell p'(\varrho)\varrho_xH_{\ell-1}
{\overline{w}-z+\overline{z}-w\over (\overline{w}-z)^2}\cr
+\ell(\ell-1)&p(\varrho)H_{\ell-2}
{(\overline{w}-z+\overline{z}-w)^2
\over (\overline{w}-z)^4}
+2\ell p(\varrho)H_{\ell-1}{\overline{w}-z+\overline{z}-w
\over (\overline{w}-z)^3},\cr
\left({\partial\over\partial y}\right)^2
g_{\alpha,\ell}(z,w)&=p''(\varrho)\varrho_y^2H_\ell+
p'(\varrho)\varrho_{yy}H_\ell
-2i\ell p'(\varrho)\varrho_yH_{\ell-1}
{\overline{w}-z-\overline{z}+w
\over (\overline{w}-z)^2}\cr
-\ell(\ell-1)&p(\varrho)H_{\ell-2}
{(\overline{w}-z-\overline{z}+w)^2
\over (\overline{w}-z)^4}+2\ell p(\varrho)H_{\ell-1}
{\overline{w}-z-\overline{z}+w
\over (\overline{w}-z)^3}.&(33.6)
}
$$
Thus
$$
\eqalignno{
(\Omega_\ell)_z\,& g_{\alpha,\ell}(z,w)
=-y^2\{(\varrho_x^2+\varrho_y^2)p''(\varrho)+
(\varrho_{xx}+\varrho_{yy})p'(\varrho)\}H_\ell\cr
&-2\ell y^2 p'(\varrho)H_{\ell-1}
{\varrho_x(\overline{w}-z+\overline{z}-w)-i\varrho_y
(\overline{w}-z-\overline{z}+w)
\over (\overline{w}-z)^2}\cr
&-\ell(\ell-1)y^2p(\varrho)H_{\ell-2}
{(\overline{w}-z+\overline{z}-w)^2-
(\overline{w}-z-\overline{z}+w)^2
\over(\overline{w}-z)^4}\cr
&-4\ell y^2 p(\varrho)H_{\ell-1}
{1\over (\overline{w}-z)^2}+2i\ell yp'(\varrho)\varrho_xH_\ell
+2i\ell^2 yp(\varrho)H_{\ell-1}
{\overline{w}-z+\overline{z}-w\over (\overline{w}-z)^2}\cr
=\,&-\left\{y^2\big\{(\varrho_x^2+\varrho_y^2)p''(\varrho)+
(\varrho_{xx}+\varrho_{yy})p'(\varrho)\big\}
+{4\ell^2yv\over |\overline{w}-z|^2}p(\varrho)\right\}
H_\ell\cr
=\,&-\left((\varrho^2+\varrho)
\left({d\over d\varrho}\right)^2
+(2\varrho+1){d\over d\varrho}
+{\ell^2\over\varrho+1}\right) p(\varrho)\cdot H_\ell.&(33.7)
}
$$
In order to relate the last line
with the hypergeometric
differential equation of Gauss, we put
$$
q(\xi)=\left({\varrho\over1+\varrho}\right)^\ell
\varrho^\alpha p(\varrho),\quad
\xi=-{1\over\varrho}\,.\eqno(33.8)
$$
Then we have
$$
\xi^2 q'(\xi)=
{\ell\varrho^{\ell-1}\over(1+\varrho)^{\ell+1}}
\varrho^\alpha
p(\varrho)+\alpha\left({\varrho\over1+\varrho}\right)^\ell
\varrho^{\alpha-1}
p(\varrho)+\left({\varrho\over1+\varrho}\right)^\ell
\varrho^\alpha p'(\varrho);\eqno(33.9)
$$
thus
$$
\left({\varrho\over1+\varrho}\right)^\ell
\varrho^\alpha p'(\varrho)=\xi^2q'(\xi)
+\xi\left({\ell \xi\over \xi-1}+\alpha \right)q(\xi).\eqno(33.10)
$$
Further,
$$
\eqalignno{
&\xi^2{d\over d\xi}\left(\xi^2 q'(\xi)
+\xi \left({\ell \xi\over \xi-1}+\alpha\right)q(\xi)\right)\cr
=&\,{\ell\varrho^{\ell-1}\over(1+\varrho)^{\ell+1}}
\varrho^\alpha
p'(\varrho)+\alpha\left({\varrho\over1+\varrho}\right)^\ell
\varrho^{\alpha-1}
p'(\varrho)+\left({\varrho\over1+\varrho}\right)^\ell
\varrho^\alpha p''(\varrho)\cr
=&\, -\xi\left({\ell\xi\over\xi-1}+\alpha\right)
\left({\varrho\over1+\varrho}\right)^\ell
\varrho^\alpha p'(\varrho)
+\left({\varrho\over1+\varrho}\right)^\ell
\varrho^\alpha p''(\varrho),&(33.11)
}
$$
which means that
$$
\eqalignno{
\left({\varrho\over1+\varrho}\right)^\ell
\varrho^\alpha p''(\varrho)&=
\xi^2{d\over d\xi}\left(\xi^2 q'(\xi)
+\xi\left({\ell\xi\over \xi-1}+\alpha\right)q(\xi)\right)\cr
&+\xi\left({\ell\xi\over\xi-1}+\alpha\right)\left(\xi^2q'(\xi)
+\xi \left({\ell\xi \over \xi-1}+\alpha \right)q(\xi)\right)\cr
&=\xi^4q''(\xi)+2\xi^3\left({\ell\xi\over\xi-1}+\alpha
+1\right)q'(\xi)\cr
&+\xi^2\left({\ell(\ell-1)\over(\xi-1)^2}+
{2\ell(\ell+\alpha\xi)\over\xi-1}+\alpha
+\alpha^2+\ell+\ell^2\right)q(\xi).&(33.12)
}
$$
We have, after some arrangements,
$$
\eqalignno{
&\left({\varrho\over1+\varrho}\right)^\ell
\varrho^\alpha\left((\varrho^2+\varrho)
\left({d\over d\varrho}\right)^2
+(2\varrho+1){d\over d\varrho}
+{\ell^2\over\varrho+1}\right)p(\varrho)\cr
=&\,\xi\Big\{\xi(1-\xi)q''(\xi)
+(2\alpha-(2\alpha+2\ell+1)\xi)
q'(\xi)-(\alpha+\ell)^2q(\xi)\Big\}
+\alpha(\alpha-1)q(\xi).\qquad&(33.13)
}
$$
Assertions $(33.7)$ and $(33.13)$ imply that 
$(33.1)$ necessitates
$$
\left[\xi(1-\xi)\left({d\over d\xi}\right)^2
+(2\alpha-(2\alpha+2\ell+1)\xi){d\over d\xi}
-(\alpha+\ell)^2\right]q(\xi)=0.\eqno(33.14)
$$
\medskip
\noindent
{\csc Notes:} The factor $H_\ell$ is due to 
Roelcke [31, Teil II, (7.11)]; however, the construction
$(33.4)$ via the Cartan decomposition does not
seem to have been reported explicitly before. This extends
to $\r{PSL}(2,\B{C})$, though in a matrix form.
\medskip
\noindent
{\bf 34. Green's function defined.}
In this way we are led to our definition of the
Green function for  
$\Omega_\ell+\alpha(\alpha-1)$ on $\B{H}^2$:
$$
\displaystyle{g_{\alpha,\ell}(z,w)
=p_{\alpha,\ell}\big(\varrho(z,w)\big)H_\ell(z,w),}
\atop\displaystyle
{p_{\alpha,\ell}(\varrho)={\Gamma(\alpha+\ell)
\Gamma(\alpha-\ell)\over
4\pi\Gamma(2\alpha)}\left(1+{1\over\varrho}
\right)^\ell\varrho^{-\alpha}\,
{}_2F_1\left(\alpha+\ell,\alpha+\ell;2\alpha;
-{1\over\varrho}\right),}\eqno(34.1)
$$
under the assumption $(33.2)$. We note that 
$$
{p_{\alpha,\ell}(\varrho)=p_{\alpha,-\ell}(\varrho),}\atop
{\overline{g_{\alpha,\ell}(z,w)}=g_{\bar\alpha,-\ell}(z,w)=
g_{\bar\alpha,\ell}(w,z),}
\eqno(34.2)
$$
since 
$$
{}_2F_1(a,b\,;c\,;\xi)
=(1-\xi)^{c-a-b}{}_2F_1\big(c-a,c-b\,;c\,;\xi\big),\quad
|\arg(1-\xi)|<\pi.\eqno(34.3)
$$
\par
In passing, we observe that for any $\r{h}\in\r{G}$ we have
$$
g_{\alpha,\ell}(z,\r{h}(w))\left({\jmath(\r{h},w)
\over\jmath(\r{h},\overline{w})}\right)^\ell
=g_{\alpha,-\ell}
(w,\r{h}^{-1}(z))\left({\jmath(\r{h}^{-1},z)
\over\jmath(\r{h}^{-1},\overline{z})}\right)^{-\ell},
\eqno(34.4)
$$
whence 
$$
\big(\Omega_{-\ell}+\alpha(\alpha-1)\big)_w\left(
g_{\alpha,\ell}(z,\r{h}(w))\left({\jmath(\r{h},w)
\over\jmath(\r{h},\overline{w})}\right)^\ell\right)=0,\quad
z\ne\r{h}(w).
\eqno(34.5)
$$
The relation $(34.4)$ is the
same as the transformation property of $H_\ell$; or more
basically, $(3.4)$ and $(33.4)$ give the equation. As to
$(34.5)$, it is as well a consequence of $(32.5)$.
\par
We shall show that $g_{\alpha,\ell}$ indeed
serves our purpose. Thus we invoke 
the Mellin--Barnes integral 
representation for the Gaussian
hypergeometric function and have
$$
p_{\alpha,\ell}(\varrho)=
{\Gamma(\alpha-\ell)\over
8\pi^2i\Gamma(\alpha+\ell)}(\varrho+1)^\ell
\int_{(-{1\over2})}
{\Gamma^2(s+\alpha+\ell)
\Gamma(-s)\over\Gamma(s+2\alpha)
}\varrho^{-s-\alpha-\ell}ds,
\eqno(34.6)
$$
where the contour is $\Re s=-{1\over2}$; and
the integrand is of exponential decay. We have, on $(33.2)$,
$$
\eqalignno{
p_{\alpha,\ell}(\varrho)&=-{\log\varrho\over4\pi}
-{1\over4\pi}\left(c_E+{\Gamma'\over\Gamma}(\alpha+\ell)+
{\Gamma'\over\Gamma}(\alpha-\ell)\right)
+O\big(\varrho|\log\varrho|\big),\quad
\varrho\to+0,\qquad&(34.7)\cr
p_{\alpha,\ell}'(\varrho)&=-{1\over 4\pi\varrho}
-{1\over4\pi}(\alpha+\ell)(\alpha-\ell-1)\log\varrho+O(1),
\quad\varrho\to+0,&(34.8)\cr
p_{\alpha,\ell}(\varrho)&
={\Gamma(\alpha+\ell)\Gamma(\alpha-\ell)
\over4\pi\Gamma(2\alpha)}
 \varrho^{-\alpha}\big(1+O(\varrho^{-1})\big),\quad
\varrho\to+\infty,&(34.9)\cr
p_{\alpha,\ell}'(\varrho)&
=-\alpha{\Gamma(\alpha+\ell)\Gamma(\alpha-\ell)
\over4\pi\Gamma(2\alpha)}
 \varrho^{-\alpha-1}\big(1+O(\varrho^{-1})\big),\quad
\varrho\to+\infty,&(34.10)
}
$$ 
where $c_E$ is the Euler constant. To show $(34.7)$, we shift the
contour in $(34.6)$ to $\Re s=
-\alpha-\ell-{3\over2}$; as to $(34.8)$, we first 
differentiate inside the integral and shift the contour 
in the same way; further,
the definition $(34.1)$ readily implies $(34.9)$--$(34.10)$.
\par
We then introduce the integral operator
$$
\c{G}_{\alpha,\ell}h(z)
=\int_{\B{H}^2}g_{\alpha,\ell}(z,w)f(w)d\mu(w).\eqno(34.11)
$$
We claim that $\c{G}_{\alpha,\ell}$ is the left inverse of
$\Omega_\ell+\alpha(\alpha-1)$: Provided $(33.2)$,
$$
\c{G}_{\alpha,\ell}\big(\Omega_\ell+\alpha(\alpha-1)\big)f=f,
\eqno(34.12)
$$
for any smooth function $f$ on $\B{H}^2$
with which the left side converges absolutely.
In fact, we can
proceed exactly the same way as in the Euclidean plane
situation, with the
observation that the left side of $(34.12)$ equals
$$
\eqalignno{
\lim_{\rho\to0}\int_{\B{H}^2\backslash\e{D}}&
\Big\{g_{\alpha,\ell}(z,w)(\Omega_0 f)(w)
-((\Omega_0)_wg_{\alpha,\ell}(z,w))f(w)\Big\}d\mu(w)\cr
&+2\ell i\lim_{\rho\to0}\int_{\B{H}^2\backslash\e{D}}
v{\partial\over\partial u}\big\{g_{\alpha,\ell}(z,w)f(w)\big\}
d\mu(w),&(34.13)
}
$$
where $\e{D}$ is an Euclidean disk of radius $\rho$ with the centre at
$z$, and we have used $(12.6)$ and $(34.5)$. The
first limit equals $f(z)$ by Green's formula $(1.9)$ and $(34.7)$;
and the second limit vanishes. In particular, we have
$$
\c{G}_{\alpha,\ell}(\Im z)^{1/2+\nu}=
{(\Im z)^{1/2+\nu}\over
\big(\alpha-{1\over2}\big)^2-\nu^2},
\quad |\Re\nu|<\Re\alpha-{1\over2}.\eqno(34.14)
$$
The range of $\nu$ follows from $(34.9)$. 
\medskip
\noindent
{\csc Notes:} Our $g_{\alpha,\ell}$ is to be compared with 
Hejhal's $k_s(z;w)$ on [11, Vol.\ 2, p.\ 350]
which is, in our notation,
$$
\eqalignno{
-&\,{\Gamma(\alpha+\ell)
\Gamma(\alpha-\ell)\over
4\pi\Gamma(2\alpha)}
\left(1-\left|{z-w\over z-\overline{w}}\right|^2
\right)^\alpha\cr
&\quad\times {}_2F_1\left(\alpha+\ell,\alpha-\ell\,;2\alpha\,;
1-\left|{z-w\over z-\overline{w}}\right|^2\right)
H_\ell(z,w),&(34.15)
}
$$
and equal to $-g_{\alpha,\ell}(z,w)$, as
$$
{}_2F_1(a,b;c\,;\xi)=(1-\xi)^{-a}
{}_2F_1\big(a,c-b\,; c\,;\xi/(\xi-1)\big),\quad |\arg(1-\xi)|<\pi.
\eqno(34.16)
$$
The reason for our use 
of the definition $(34.1)$ is in that
we are able to exploit $(34.6)$, which is of rapid convergence
and thus convenient to manipulate.
For $(34.3)$ and $(34.16)$ see Lebedev [18, Section 9.5].
A further comparison should be made with Roelcke [31, Teil II, \S7].
See also Bruggeman [4, p.\ 70]. Here is an extra observation:
The Green function $g_{\alpha,0}$ appears in a peculiar way as
the kernel function that connects the
fourth power moment of the Riemann zeta-function to
the spectral decomposition of $L^2(\varGamma\backslash\r{G})$;
moreover, $g_{\alpha,0}$ can be expressed as a convolution of
Bessel functions of representation. For these 
see [5] and [23, Part XII].
There is an extension to the complex situation; 
see [24]. Commenting further, the product of two automorphic 
$L$-functions is naturally related to the product of two instances
of $(29.4)$, which can be viewed as a Mellin transform of
the convolution of two instances of $(29.1)$, and the
latter is essentially the Gaussian hypergeometric function ${}_2F_1$. 
The appearance of ${}_2F_1$ 
in the spectral decomposition [22, Theorem 4.1] 
of the fourth moment
of the Riemann zeta-function is somewhat mysterious but might
be understood as an indication that
the unitary map $r_\r{w}^{V_1}r_\r{w}^{V_2}$
of the space $L^2(\B{R}^\times, d^\times\!/\pi)$ 
into itself is behind the scene.
\medskip
\noindent
{\bf 35. Automorphic Green's function.}
If $f$ is in $L^2_\ell(\varGamma\backslash
\B{H}^2)$ we have formally
$$
\c{G}_{\alpha,\ell}f(z)=\int_{\varGamma\backslash\B{H}^2}
G_{\alpha,\ell}(z,w)f(w)d\mu(w),\eqno(35.1)
$$
where
$$
G_{\alpha,\ell}(z,w)=
\sum_{\gamma\in\varGamma}g_{\alpha,\ell}(z,\gamma(w))
\left({\jmath(\gamma,w)\over\jmath(\gamma,
\overline{w})}\right)^\ell.\eqno(35.2)
$$
By $(34.9)$ we have,  with $z$ 
being bounded and $w\in\e{F}$ tending to the cusp,
$$
\eqalignno{
G_{\alpha,\ell}(z,w)&\ll\sum_{\gamma
\in\varGamma_\infty\backslash\varGamma}
\sum_{n=-\infty}^\infty p_{\alpha,\ell}
\big(\r{s}^n\gamma(z),w\big)\cr
&\ll\sum_{\gamma\in\varGamma_\infty\backslash\varGamma}
\sum_{n=-\infty}^\infty
\bigg({n^2+(\Im w)^2\over 
\Im w\,\Im\gamma(z)}\bigg)^{-\Re\alpha}\cr
&\ll(\Im w)^{1-\Re\alpha}
\sum_{\gamma\in\varGamma_\infty\backslash\varGamma}
\big(\Im\gamma(z)\big)^{\Re\alpha}\ll (\Im w)^{1-\Re\alpha},
&(35.3)
}
$$
where $\r{s}=\left[{1\atop}{1\atop1}\right]$.
Since the last sum has to be convergent,
we need to have $\Re\alpha>1$; see $(33.2)$. 
Thus, for the basis vectors 
$$
\psi_V^{(\ell)}(z)=\lambda_V^{(\ell)}(\r{n}[x]\r{a}[y]),
\quad z=x+iy,\eqno(35.4)
$$ 
of the cuspidal subspace of 
$L^2_\ell(\varGamma\backslash\B{H}^2)$ which are
induced by $(22.1)$--$(22.4)$, we have, on
noting that $\psi_V^{(\ell)}$ is bounded throughout $\B{H}^2$,
$$
\c{G}_{\alpha,\ell}\psi_V^{(\ell)}={\psi_V^{(\ell)}
\over\big(\alpha-{1\over2}\big)^2-\nu_V^2},\eqno(35.5)
$$
under the convention $(35.1)$.
\par
We extend this to the Eisenstein series of 
weight $2\ell$ on $\B{H}^2$: Namely, for the function
$$
e_\ell(z,\nu)=E_\ell(\r{n}[x]\r{a}[y],\nu),\eqno(35.6)
$$
we shall show that
$$
\c{G}_{\alpha,\ell}e_\ell(z,\nu)={e_\ell(z,\nu)\over
\big(\alpha-{1\over2}\big)^2-\nu^2},
\quad e_\ell(z,\nu)\ne\infty,\, 
|\Re\nu|<\Re\alpha-{1\over2}.\eqno(35.7)
$$
To this end, we assume
temporarily that ${1\over2}<\Re\nu<\Re\alpha-{3\over2}$ holds;
this means that we work in a sub-domain of $(33.2)$ 
for the moment. Then $(34.14)$ implies that 
$$
\eqalignno{
{e_\ell(z,\nu)\over \big(\alpha-{1\over2}\big)^2-\nu^2}&=
\sum_{\gamma\in\varGamma_\infty\backslash\varGamma}
\left({\jmath(\gamma,z)\over
\jmath(\gamma,\bar{z})}\right)^{-\ell}
\int_{\B{H}^2} g_{\alpha,\ell}(\gamma(z),w)
(\Im w)^{1/2+\nu}d\mu(w)\cr
&=\sum_{\gamma\in\varGamma_\infty\backslash\varGamma}
\int_{\B{H}^2} g_{\alpha,\ell}(z,w)
(\Im \gamma(w))^{1/2+\nu}\left({\jmath(\gamma,w)\over
\jmath(\gamma,\bar{w})}\right)^\ell d\mu(w)\cr
&=\int_{\B{H}^2} g_{\alpha,\ell}(z,w)
e_\ell(w,\nu)d\mu(w).&(35.8)
}
$$
The last line is due to absolute convergence coupled with
the observation that by $(18.2)$, for $\ell=0$,
$$
\eqalignno{
\sum_{\gamma\in\varGamma_\infty\backslash\varGamma}
(\Im \gamma(w))^{1/2+\Re\nu}
&\le e_0(iv,\Re\nu)=e_0(i/v,\Re\nu)\cr
&\ll v^{1/2+\Re\nu}+v^{-1/2-\Re\nu}.&(35.9)
}
$$
Thus we have
$$
{e_\ell(z,\nu)\over \big(\alpha-{1\over2}\big)^2-\nu^2}
=\int_{\varGamma\backslash\B{H}^2}
G_{\alpha,\ell}(z,w)e_\ell(w,\nu)d\mu(w).\eqno(35.10)
$$
This converges uniformly for $\alpha,\nu$ in the domain
indicated in $(35.7)$, 
and by analytic continuation 
we end the proof of $(35.7)$.
\par
Further, we note that
$$
\hbox{$G_{\alpha,\ell}(z,w)$ is 
of weight $2\ell$ with respect to $z$}.\eqno(35.11)
$$
In fact, we have, for any $\tau\in\varGamma$ and
in the region of absolute convergence,
$$
G_{\alpha,\ell}(\tau(z),w)=\sum_{\gamma\in
\varGamma}p_{\alpha,\ell}
\big(\varrho(\tau(z),\tau\gamma(w))\big)
H_\ell(\tau(z),\tau\gamma(w))
\left({\jmath(\tau\gamma,w)
\over\jmath(\tau\gamma,\overline{w})}\right)^\ell.\eqno(35.12)
$$ 
To the first factor in the summand we apply the second line in
$(2.1)$, and to the remaining factors
the third and the fourth identities in $(1.2)$.
\medskip
\noindent
{\csc Notes:} The discussion of the previous and the present 
sections corresponds to [22, Section 1.3]. 
Some simplifications have been made. 
\medskip
\noindent
{\bf 36. Iterated kernel.}
We now turn to the function
$$
G^{(\ell)}(z,w)=G_{\alpha,\ell}(z,w)-G_{\beta,\ell}(z,w),\quad
\Re\beta>\Re\alpha,\eqno(36.1)
$$
in which $(33.2)$ is still imposed. 
This is related to the iteration of
the integral transform $(34.11)$, for
we have the Hilbert identity
$$
\c{G}_{\alpha,\ell}
\c{G}_{\beta,\ell}={\c{G}_{\alpha,\ell}-\c{G}_{\beta,\ell}\over
(\beta-\alpha)(\beta+\alpha-1)},\eqno(36.2)
$$
as is well indicated by $(34.12)$; however, this fact is not
needed in our discussion below.
\par
The assertion $(34.7)$ implies that $G^{(\ell)}(z,w)$ is
continuous on the diagonal $z=w$; 
and by $(35.3)$ we see that $G^{(\ell)}(z,w)=
\overline{G^{(-\ell)}(w,z)}$ is bounded when
$z\in\e{F}$ tends to the cusp while $w$ is bounded.
Hence, in view of $(35.11)$, we may apply 
$(17.7)$ and $(22.10)$ to $G^{(\ell)}(z,w)$ as a function
of $z$. 
We have that for each fixed $w\in\B{H}^2$
$$
G^{(\ell)}(z,w)={}^{0}G(z,w)+{}^c G(z,w),\eqno(36.3)
$$
where
$$
\eqalignno{
{}^0G(z,w)&={\sum_V}^{(\ell)} k(\nu_V)\psi_V^{(\ell)}(z)\,
\overline{\psi_V^{(\ell)}(w)},&(36.4)\cr
{}^c G(z,w)&={1\over4\pi i}\int_{(0)}
k(\nu)e_\ell(z,\nu)\overline{\,e_\ell(w,\nu)}d\nu
+{3\over\pi}\delta_0k\big(\txt{1\over2}\big),&(36.5)
}
$$
with convergence in the mean in
$L^2_\ell(\varGamma\backslash\B{H}^2)_z$. Here
$\delta_0$ is equal to $1$ if $\ell=0$ and to $0$ otherwise;
and
$$
k(\nu)={1\over\big(\alpha-{1\over2}\big)^2-\nu^2}
-{1\over\big(\beta-{1\over2}\big)^2-\nu^2}.\eqno(36.6)
$$
In fact $(36.4)$ follows from $(35.5)$, and the integrated part
of $(36.5)$ from $(35.7)$. As to the constant term of
$(36.5)$, we note that the assertion
$\big\langle{G^{(0)}(z,\cdot),1}\big\rangle
=k\big({1\over2}\big)$ follows from a combination of
$(35.7)$, $\ell=0$, and
$e_0\big(z,-{1\over2}\big)=1$; the latter comes from the fact
that in $(18.2)$, $\ell=0$,
we have $\varphi_\varGamma\big(\!-{1\over2}\big)=0$ as well as
$\e{A}^\delta\phi_0\big(\r{g},-{1\over2}\big)=0$ because of
$(16.1)$.
\par
We shall demonstrate in the next section
that the series $(36.4)$ converges
absolutely and uniformly over the product domain $\e{F}\times
\e{F}$. To achieve this, we shall appeal to the
well-known theorem of Mercer
on the spectral expansion of positive symmetric kernels. 
As a prerequisite, we need to
have that 
$$
\hbox{$
{}^0G(z,w)=\big\{G^{(\ell)}-{}^cG\big\}(z,w)$
is continuous and bounded on $\e{F}
\times\e{F}$. }\eqno(36.7)
$$
\medskip
\noindent
{\csc Notes:} For Mercer's theorem see 
Riesz and Sz-Nagy [30, Section 98]. Also we note 
that the functions
$$
\eqalignno{
p_{\alpha,\ell}^{(n)}(\varrho)&=
{\Gamma(\alpha-\ell)\over
8\pi^2i\Gamma(n+1)\Gamma(n+\alpha+\ell)}\cr
\times(\varrho+1)^\ell&\int_{(-{1\over2})}
{\Gamma(s+n+\alpha+\ell)\Gamma(s+\alpha+\ell)
\Gamma(-s)\over\Gamma(s+n+2\alpha)
}\varrho^{-s-\alpha-\ell}ds,
&(36.8)
}
$$
with integers $n\ge0$, satisfy
$$
p_{\alpha,\ell}^{(n)}-p_{\alpha+1,\ell}^{(n)}=(n+1)(n+2\alpha)
p_{\alpha,\ell}^{(n+1)},
\quad p_{\alpha,\ell}^{(0)}=p_{\alpha,\ell},\eqno(36.9)
$$
which is an extension of [22, $(1.3.12)$]. 
The fast converging
expression $(36.8)$ serves the same purpose as $(34.6)$ does,
when discussing the spectral
resolution of $\Omega_\ell$ via
iterations of $\c{G}_{\alpha,\ell}$.
See the notes to the next section.
\medskip
\noindent
{\bf 37. Control of divergence.} We shall prove that
the functions $G^{(\ell)}(z,w)$ 
and ${}^cG(z,w)$ diverge in the same mode as $z,w\in\e{F}$ 
tend to the cusp, and confirm the assertion $(36.7)$. 
Our argument is analogous to that of
[22, Lemmas 1.7 and 1.8].
\par
We shall first prove an approximation to $G^{(\ell)}$:
Uniformly for $z=x+iy,w=u+iv\in\e{F}$,
$$ 
\left|G^{(\ell)}(z,w)-{1\over2\pi i}\int_{(0)}
k(\nu)y^{{1/2}+\nu}
v^{{1/2}-\nu}d\nu\right|\ll1.\eqno(37.1)
$$
To this end, we write, with $z,w\in\B{H}^2$,
$$
\eqalign{
G^{(\ell)}(z,w)&=
(S_\alpha-S_\beta)(z,w)+(T_\alpha-T_\beta)(z,w),\cr
S_\alpha(z,w)&=\sum_{n=-\infty}^\infty g_{\alpha,\ell}(z,n+w),\cr
T_\alpha(z,w)&=
\sum_{\scr{\gamma\in\varGamma_\infty\backslash\varGamma}
\atop\scr{\gamma\not\in\varGamma_\infty}}
S_\alpha\big(z,\gamma(w)\big)
\left({\jmath(\gamma,w)\over
\jmath(\gamma,\overline{w})}\right)^\ell.
}\eqno(37.2)
$$ 
By the Euler--Maclaurin sum formula, we have,
with $\xi(t)=t-[t]-{1\over2}$,
$$
\eqalignno{
S_\alpha(z,w)&=\int_{-\infty}^\infty g_{\alpha,\ell}(z,t+w)dt
+\int_{-\infty}^\infty \xi(t){\partial\over
\partial t}g_{\alpha,\ell}(z,t+w)dt\cr
&=\big\{S_\alpha^{(0)}+S_\alpha^{(1)}\big\}(z,w),&(37.3)
}
$$
say. We have the decomposition
$$
\eqalign{
G^{(\ell)}(z,w)&=\Big\{(S^{(0)}_\alpha-S^{(0)}_\beta)
+(T^{(0)}_\alpha-T^{(0)}_\beta)
+(T^{(1)}_\alpha-T_\beta^{(1)})\Big\}(z,w),\cr
T^{(0)}_\alpha(z,w)
&=\sum_{\scr{\gamma\in\varGamma_\infty\backslash\varGamma}
\atop\scr{\gamma\not\in\varGamma_\infty}}
S^{(0)}_\alpha\big(z,\gamma(w)\big)
\left({\jmath(\gamma,w)\over
\jmath(\gamma,\overline{w})}\right)^\ell,\cr
T^{(1)}_\alpha(z,w)&
=\sum_{\gamma\in\varGamma_\infty\backslash\varGamma}
S^{(1)}_\alpha\big(z,\gamma(w)\big)
\left({\jmath(\gamma,w)\over
\jmath(\gamma,\overline{w})}\right)^\ell.
}\eqno(37.4)
$$
We assert that for any $z,w\in\B{H}^2$ the difference
$\big(S^{(0)}_\alpha-S^{(0)}_\beta\big)(z,w)$ 
equals the integrated term
in $(37.1)$.  In fact, the formula $(34.14)$ implies that
$$
\int_0^\infty \left(\int_{-\infty}^\infty 
g_{\alpha,\ell}(z,u+iv)du\right)
v^{\nu-3/2}dv={y^{1/2+\nu}\over
\big(\alpha-{1\over2}\big)^2-\nu^2}.\eqno(37.5)
$$
Since the inner integral is a continuous function of $v$, we have,
via the Mellin inversion procedure,
$$
\int_{-\infty}^\infty 
g_{\alpha,\ell}(z,u+iv)du=
{1\over2\pi i}\int_{(0)}{
y^{1/2+\nu}v^{1/2-\nu}\over
\big(\alpha-{1\over2}\big)^2-\nu^2}d\nu,\eqno(37.6)
$$
which proves the claim. To deal with $T^{(0)}_\alpha$,
we shift the last contour to $\Re\nu=-\Re\alpha$:
$$
S_\alpha^{(0)}(z,w)={y^{1-\alpha} v^\alpha\over 2\alpha-1}
+{1\over2\pi i}\int_{(-\Re\alpha)}
{y^{1/2+\nu}v^{1/2-\nu}\over
\big(\alpha-{1\over2}\big)^2-\nu^2}d\nu.\eqno(37.7)
$$
Thus, we have that for any $z,w\in\B{H}^2$
$$
\eqalignno{
&T^{(0)}_\alpha(z,w)={y^{1-\alpha}\over 2\alpha-1}
\big(e_\ell(w,\alpha-\txt{1\over2})-v^\alpha\big)\cr
&+{1\over2\pi i}\int_{(-\Re\alpha)}{
y^{1/2+\nu}\over
\big(\alpha-{1\over2}\big)^2-\nu^2}\big(e_\ell(w,-\nu)
-v^{1/2-\nu}\big)d\nu.&(37.8)
}
$$
On the other hand, we have, for
$z,w\in\B{H}^2$,
$$
S^{(1)}_\alpha(z,w)=\int_{-\infty}^\infty\xi(t)\left(
p'(\varrho)\varrho_uH_\ell-\ell
p(\varrho)H_{\ell-1}{\overline{w}-z+\overline{z}-w
\over (\overline{w}+t-z)^2}\right)dt,\eqno(37.9)
$$
where $p=p_{\alpha,\ell}$, $\varrho=\varrho(z,w+t)$,
$H_\ell=H_\ell(z,w+t)$; thus, 
$$
(S^{(1)}_\alpha-S^{(1)}_\beta)(z,w)
\ll \int_{-\infty}^\infty\left(|t||p'|+
|p|{(y+v)/\sqrt{yv}\over t^2+(y+v)^2/yv}\right)dt,\eqno(37.10)
$$
where $
p=(p_{\alpha,\ell}-p_{\beta,\ell})(t^2+\varrho(iy,iv))$,
$p'=(p_{\alpha,\ell}'-p_{\beta,\ell}')(t^2+\varrho(iy,iv))$.
We have, by $(34.7)$--$(34.10)$,
$$
p(\tau)\ll\min\big\{1, \tau^{-\Re\alpha}\big\},\quad p'(\tau)\ll
\min\big\{|\log\tau|,\tau^{-\Re\alpha-1}\big\}.\eqno(37.11)
$$
Considering the cases $y/v\le{1/2}$, ${1\over2}\le y/v\le2$, and
$y/v\ge2$ separately, we find that
$$
\big(S^{(1)}_\alpha-S^{(1)}_\beta\big)(z,w)\ll 
\left({yv\over (y+v)^2}\right)^{\Re\alpha}.\eqno(37.12)
$$
Hence we have that for any $z,w\in\B{H}^2$
$$
T_\alpha^{(1)}(z,w)\ll \left({yv\over (y+v)^2}\right)^{\Re\alpha}
+y^{-\Re\alpha}\left(e_0\big(w,\Re\alpha-\txt{1\over2}\big)
-v^{\Re\alpha}\right).\eqno(37.13)
$$
Inserting the bound $(15.5)$ and 
the expansion $(18.2)$ into $(37.8)$ and $(37.13)$,
with $z,w\in\e{F}$,
we end the proof of $(37.1)$.
\par
We turn to ${}^cG(z,w)$. We denote the sum over 
$n$ in $(18.2)$, for $\r{g}=\r{n}[x]\r{a}[y]$, 
by $\tilde{e}_\ell(z,\nu)$. Invoking the functional equation $(18.4)$
we get the decomposition
$$
\eqalignno{
{}^cG(z,w)&={1\over{2\pi i}}\int_{(0)}
k(\nu)y^{{1/2}-\nu}v^{{1/2}+\nu}d\nu\cr
&+\big\{{}^cG^{(0)}+{}^cG^{(1)}
+{}^cG^{(2)}+{}^cG^{(3)}\big\}(z,w), &(37.14)
}
$$
where
$$
\eqalignno{
&{}^cG^{(0)}(z,w)={(-1)^\ell\over{2\pi}i}\int_{(0)}
k(\nu)(yv)^{{1/2}-\nu}
{\Gamma^2\big(\nu+{1\over2}\big)
\varphi_\varGamma(\nu)
\over\Gamma\big(\nu+|\ell|+{1\over2}\big)
\Gamma\big(\nu-|\ell|+{1\over2}\big)}
d\nu+{3\over\pi}\delta_0
k\big(\txt{1\over2}\big),\qquad\cr
&{}^cG^{(1)}(z,w)={1\over{2\pi}i}\int_{(0)}
k(\nu)y^{{1/2}-\nu}\tilde{e}_\ell(w,\nu)d\nu,\quad
{}^cG^{(2)}(z,w)=\overline{{}^cG^{(1)}(w,z)},&(37.15)\cr
&{}^cG^{(3)}(z,w)={1\over{4\pi i}}\int_{(0)}k(\nu)
\overline{\tilde{e}_\ell(z,\nu)}\,\tilde{e}_\ell(w,\nu)d\nu.
}
$$
In what follows we shall use well-known
estimations of $1/\zeta(s)$ for $\Re s\ge1$ and
$\zeta(s)$ for $\Re s\le1$, without mentioning. 
Also we may assume that 
$yv\ge1$ without loss of generality.
\par
The treatment of ${}^cG^{(0)}$ is easy. It suffices to move the
contour to $(\alpha)$. When $\ell=0$, we pass the simple
pole at $\nu={1\over2}$, whose contribution is cancelled
by the constant term on the right side. We have
$$ 
{}^cG^{(0)}(z,w)\ll (yv)^{1-\Re\alpha}.\eqno(37.16)
$$ 
We shall consider ${}^cG^{(2)}(z,w)$. 
Invoking $(15.5)$, we have, uniformly for $y\ge1$ and
$\Re\nu\ge0$,
$$
\eqalignno{
\tilde{e}_\ell(z,\nu)
&\ll{(1+|\nu|)y^{-1/2-\Re\nu}
\over{|\zeta(1+2\nu)|}}
\exp\left(\!-{\pi{y}\over1+|\nu|}\right)\sum_{n=1}^\infty
{d(n)\over n}\exp\left(\!-{\pi{yn}\over 1+|\nu|}\right)\cr
&\ll{(1+|\nu|)y^{-1/2-\Re\nu}
\over{|\zeta(1+2\nu)|}}
\left(\log\big(2+|\nu|/ y\big)\!\right)^2
\exp\left(\!-{\pi{y}\over1+|\nu|}\right),&(37.17)
}
$$
where $d(n)$ is the number of divisors of $n$. 
Hence, shifting the paths 
in ${}^cG^{(1)}(z,w)$ and ${}^cG^{(2)}(z,w)$ to
$+\infty$, we get
$$
\eqalign{
{}^cG^{(1)}(z,w)&\ll{y}^{1-\Re\alpha}
\exp\big(-v/(|\alpha|+|\beta|)\big),\cr
{}^cG^{(2)}(z,w)&\ll{v}^{1-\Re\alpha}
\exp\big(-y/(|\alpha|+|\beta|)\big).
}\eqno(37.18)
$$
As to ${}^cG^{(3)}(z,w)$, we apply $(37.17)$ with $\Re\nu=0$,
$$
\eqalignno{
{}^cG^{(3)}(z,w)&\ll(yv)^{-{1/2}}
\int_{(0)}{{\exp\big(\!-\pi(y+v)/(1+|\nu|)\big)}
\over{|\zeta(1+2\nu)|^2
(1+|\nu|)^2}}\big(\log(2+|\nu|)\big)^4|d\nu|\cr
&\ll (yv)^{-1}(\log2yv)^6.&(37.19)
}
$$ 
Collecting $(37.14)$, $(37.18)$, $(37.19)$, we get,
uniformly for $z,w\in\e{F}$,
$$
\left|{}^cG(z,w)-{1\over{2\pi i}}\int_{(0)}
k(\nu)y^{{1/2}-\nu}v^{{1/2}+\nu}d\nu\right|\ll 1.\eqno(37.20)
$$
Combined with $(37.1)$, this yields the assertion $(36.7)$.
\par
We now have that
$$
\int_{\varGamma\backslash\B{H}^2}\left(
\int_{\varGamma\backslash\B{H}^2}
\big|{}^0 G(z,w)\big|^2d\mu(z)\right)d\mu(w)
={\sum_V}^{(\ell)} |k(\nu_V)|^2\eqno(37.21)
$$
is convergent. We see via Fubuni's
theorem that $(36.4)$ converges in the mean
in the space $L^2\big((\varGamma\backslash\B{H}^2)
\times(\varGamma\backslash\B{H}^2)\big)$. Hence, by
Mercer's theorem we conclude that
$(36.4)$ converges uniformly over
$\e{F}\times\e{F}$. The same holds with $(36.5)$,
as has been proved in the above.
\medskip
\noindent
{\csc Notes:} A way to approach relatively 
directly to the spectral resolution 
of $\Omega_\ell$ is to elaborate a little 
the present section.
A key is the integral expression $(36.8)$. 
With it, one may proceed
in much the same way as in [22, Sections 1.3--1.4]. There is an
issue which is peculiar to the situation $\ell\ne0$ and 
related to holomorphic cusp forms, i.e., those 
vectors in $L^2_\ell(\varGamma\backslash\r{G})$ which
correspond to eigenvalues $l-{1\over2}$, $l\in\B{N}$, 
$l\le\ell$. It can,
however, be handled by following Section 21.
This means, in particular, that the argument eventually
involves the Maass operators disguised as $(32.1)$.
Another extension is to consider Hecke congruence subgroups
in place of the full modular group.
Then the above corresponds to the control of the
scattering of the relevant automorphic Green function
at incongruent cusps; and our method works without
essential changes. See [26, Part II] 
and [27, Vol.\ 2, Appendix 1].
\medskip
\noindent
{\bf 38. Trace of the Casimir operator.} 
Having confirmed $(36.7)$,
we are now able to consider
$$
\int_{\e{F}}\big\{G^{(\ell)}-{}^cG\big\}(z,z)
d\mu(z)={\sum_V}^{(\ell)}k(\nu_V).\eqno(38.1)
$$
The right side is a trace of $\Omega_\ell$ the restriction
of the Casimir operator to $L^2_\ell(\varGamma\backslash
\r{G})$, while the left side is defined and thus can be computed,
in terms of elements of $\varGamma$.
This leads us to an instance of Selberg's trace formulas. However,
we shall develop only a trace formula for the operator $\Omega_0$,
the reason for which is in the identity
$$
{\sum_V}^{(\ell)}k(\nu_V)=
{\sum_V}^{(0)}k(\nu_V)+\sum_{l=1}^{\ell}
(\ell-l+1)\vartheta_\varGamma(l)
k\big(l-\txt{1\over2}\big),\eqno(38.2)
$$
where $\vartheta_\varGamma(l)$ is defined in Section 21. 
This follows from $(22.8)$ or $(32.12)$. 
In the first sum on the right side,
$V$ runs over all irreducible representations in the
unitary principal series.
\par
Hence, what counts is to compute the part $\sum^{(0)}_V$ 
which is a trace of $\Omega_0$. With this,  we shall assume
hereafter that
$$
\ell=0,\quad \Re\beta>\Re\alpha>1.\eqno(38.3)
$$
We also simplify the notation:
$$
p_\alpha=p_{\alpha,0}, \quad g_\alpha=g_{\alpha,0}
\eqno(38.4)
$$
\par
To compute the geometric side of $(38.1)$, for $\ell=0$,
we introduce the classification of the elements $\gamma
\in\varGamma$:
$$
\varGamma={\cal C}^{(0)}\sqcup{\cal C}^{(1)}
\sqcup{\cal C}^{(2)}\sqcup{\cal C}^{(e)}.
\eqno(38.5)
$$
Here ${\cal C}^{(0)}=\{1\}$, and 
$$
\eqalign{
&{\cal C}^{(1)}=\big\{\hbox{$\gamma$ has a single 
fixed point on $\B{R}\cup\{\infty\}$}\big\}: 
\hbox{parabolic,}\cr
&{\cal C}^{(2)}=\big\{\hbox{$\gamma$ has two different
fixed points
on $\B{R}\cup\{\infty\}$}\big\}: \hbox{hyperbolic},\cr
&{\cal C}^{(e)}=\big\{\hbox{$\gamma$ has a single
fixed point inside $\B{H}^2$}\big\}: \hbox{elliptic.}
}\eqno(38.6)
$$
We divide the sum $(35.2)$, for $\ell=0$, 
according to $(38.5)$, and
denote the corresponding parts of $G^{(0)}$ 
by $R^{(0)}$, $R^{(1)}$,
$R^{(2)}$, and $R^{(e)}$, respectively. The subgroup
$\varGamma_\infty$ is contained in ${\cal C}^{(0)}
\sqcup{\cal C}^{(1)}$; thus $R^{(2)}$ and $R^{(e)}$ are
bounded on $\e{F}\times\e{F}$, which is obvious from the
proof of $(37.1)$, or more precisely from the
boundedness of $T_\alpha$ defined in $(37.2)$. We have
$$
\eqalignno{
\mathop{{\sum}^\r{u}}_{V\hfil}k(\nu_V)
&=\int_{\eusm F}\big\{R^{(0)}+R^{(2)}
+R^{(e)}\big\}(z,z)d\mu(z)\cr
&+\lim_{Y\to\infty}\int_{{\eusm F}_Y}
\big\{R^{(1)}-{}^cG\big\}(z,z)d\mu(z),&(38.7)
}
$$
with $\e{F}_Y=\e{F}\cap\big\{\Im z\le Y\big\}$. Here
${\sum}^\r{u}$ indicates that the sum is restricted to the
unitary principal series, that is, the same as ${\sum}^{(0)}$.
\par
It is immediate to see that $(34.7)$ gives
$$
\int_{\eusm F}R^{(0)}(z,z)d\mu(z)
=-{1\over6}\left({\Gamma'\over\Gamma}
(\alpha)-{\Gamma'\over\Gamma}(\beta)\right).
\eqno(38.8)
$$
\medskip
\noindent
{\csc Notes:} The formula $(38.1)$ with $(38.2)$ 
could be used to evaluate $\vartheta_\varGamma(l)$ 
for any $l\in\B{N}$, supplying an alternative proof
of the dimension formula $(21.8)$. This is a typical instance of
Selberg's general observation  
made in his seminal work [32].
Hejhal executed this programme in
[11, Vol.\ 2, Chapters 9--10] with a fairly general choice of the
underlying discrete subgroups;  
his discussion is built on Roelcke's investigation [31], though. 
It should be worth mentioning that in the actual computation of
$\vartheta_\varGamma(l)$ it suffices to consider 
the situation with $\ell=l$, $\alpha$
in an immediate neighbourhood of $l$ and
$\beta$ tending to $+\infty$, which simplifies the
discussion substantially. 
\medskip
\noindent
{\bf 39. Parabolic terms.} We shall deal with $R^{(1)}$. To this
end, we note that
$$
{\cal C}^{(1)}=
\bigsqcup_{\gamma\in\varGamma_\infty\backslash\varGamma}
\left\{\gamma^{-1}
\r{s}^n\gamma:\,{\Bbb Z}\ni n\not=0\right\},\eqno(39.1)
$$
where $\r{s}$ is as in $(35.3)$. In fact, the fixed point of any
$\omega\in{\cal C}^{(1)}$ is a rational number or the cusp. Thus
there exists a $\gamma\in\varGamma$ such that the fixed
point of $\gamma\omega\gamma^{-1}$ is the cusp, and
$\gamma\omega\gamma^{-1}=\r{s}^n$ with a certain 
$n\in\B{Z}$. Also, if $\gamma^{-1}\r{s}^m\gamma
=\gamma'^{-1}\r{s}^n\gamma'$ with $m,n\in\B{Z}$, then 
 we get $m=n$ and $\gamma'\gamma^{-1}
\in\varGamma_\infty$, which confirms $(39.1)$. Thus we have
$$
\eqalignno{
R^{(1)}(z,z)&=2\sum_{\gamma\in\varGamma_\infty
\backslash\varGamma}\;
\sum_{n=1}^\infty (p_\alpha-p_\beta)
\big((n/2\Im\gamma(z))^2\big)\cr
&=\big\{R^{(1)}_\alpha-R^{(1)}_\beta\big\}(z,z),&(39.2)
}
$$
say. By $(34.6)$, for $\ell=0$, we have
$$
R_\alpha^{(1)}(z,z)
={2^{2\alpha}\over4\pi^2i}
\int_{\left(-{1\over2}\right)}
{2^{2\xi}\Gamma^2(\xi+\alpha)
\Gamma(-\xi)\over\Gamma(\xi+2\alpha)}
\zeta\big(2(\xi+\alpha)\big)e_0\big(z,2(\xi+\alpha)
-\txt{1\over2}\big)d\xi.\eqno(39.3)
$$
We integrate this over $\e{F}_Y$.  
We observe that
$$
\eqalignno{
\int_{\e{F}_Y}e_0\big(z,2(\xi+\alpha)
-\txt{1\over2}\big)d\mu(z)&=
\int_{\e{F}_Y}{(-\Omega_0e_0)\big(z,2(\xi+\alpha)
-\txt{1\over2}\big)\over
2(\xi+\alpha)\big(2(\xi+\alpha)-1\big)}d\mu(z)\cr
&=\int_{-{1\over2}}^{1\over2}
{(\partial/\partial y)_{y=Y}e_0\big(z,2(\xi+\alpha)
-\txt{1\over2}\big)\over
2(\xi+\alpha)\big(2(\xi+\alpha)-1\big)}dx, &(39.4)
}
$$
where we have used Green's formula $(1.9)$, i.e.,
$(32.7)$ for $\ell=0$, as well as the cancellation 
noted after $(32.11)$. We find, via $(18.2)$ for $\ell=0$,
that
$$
\eqalignno{
\int_{{\eusm F}_Y}R^{(1)}_\alpha(z,z)d\mu(z)
=&\,{2^{2\alpha}\over4\pi^2i}\int_{\left(-{1\over2}\right)}
{2^{2\xi}\Gamma^2(\xi+\alpha)\Gamma(-\xi)
\zeta\big(2(\xi+\alpha)\big)
\over\Gamma(\xi+2\alpha)}\cr
\times&\bigg({Y^{2(\xi+\alpha)-1}\over2(\xi+\alpha)-1}
-{Y^{-2(\xi+\alpha)}\over2(\xi+\alpha)}
\varphi_\varGamma\big(2(\xi+\alpha)-\txt{1\over2}\big)
\bigg)d\xi.&(39.5)
}
$$
We shift the contour to $\Re\xi={1\over4}-\Re\alpha$;
the relevant singularity is the pole of order 2
at $\xi={1\over2}-\alpha$. We get
$$
\int_{{\eusm F}_Y}R^{(1)}_\alpha(z,z)d\mu(z)
={1\over2\alpha-1}\Big
\{\log \left(\txt{1\over2}Y\right)
-{\Gamma'\over\Gamma}\big(\alpha+\txt{1\over2}\big)\Big\}
+{1\over(2\alpha-1)^2}+O\big(Y^{-{1/2}}\big).
\eqno(39.6)
$$
\par
On the other hand, we have
$$
\int_{{\eusm F}_Y}{}^cG(z,z)d\mu(z)
={1\over4\pi i}\int_{(0)}k(\nu)\int_{{\eusm F}_Y}
|e_0(z,\nu)|^2d\mu(z)d\nu+h\big(\txt{1\over2}\big)
+O\big(Y^{-1}\big).\eqno(39.7)
$$
We then need the following formula: 
For $\nu_1,\nu_2\in i\B{R}$, $\nu_1\ne\pm\nu_2$,
$$
\eqalignno{
\int_{\e{F}_Y}&e_0(z,\nu_1)e_0(z,-\nu_2)d\mu(z)
={{Y^{\nu_1-\nu_2}-Y^{\nu_2-\nu_1}}\over{\nu_1-\nu_2}}
+{{Y^{\nu_2-\nu_1}}\over{\nu_1-\nu_2}}
\big(1-\varphi_\varGamma(\nu_1)
\varphi_\varGamma(-\nu_2)\big)\cr
&-{{Y^{-\nu_1-\nu_2}}\over{\nu_1+\nu_2}}
\varphi_\varGamma(\nu_1)
+{{Y^{\nu_1+\nu_2}}\over{\nu_1+\nu_2}}
\varphi_\varGamma(-\nu_2)
-\int_{{\eusm{F}}\backslash{\eusm{F}}_Y}
\tilde{e}_0(z,\nu_1)\tilde{e}_0(z,-\nu_2)d\mu(z),
&(39.8)
}
$$ 
where $\tilde{e}_0$ is as in $(37.15)$.
To prove this, we apply $(1.9)$ to see that
the left side equals
$$
\big(\nu_2^2-\nu_1^2\big)^{-1}\int_{-{1\over2}}^{1\over2}
\left[e_0(z,\nu_1){\partial e_0\over{\partial{y}}}(z,-\nu_2)
-e_0(z,-\nu_2){\partial e_0\over{\partial{y}}}
(z,\nu_1)\right]_{y=Y}dx,\eqno(39.9)
$$ 
which is
$$
\eqalignno{
\big(\nu_2^2-\nu_1^2\big)^{-1}
\Bigg\{&-(\nu_1+\nu_2)Y^{\nu_1-\nu_2}
+(\nu_1-\nu_2)\varphi_\varGamma(\nu_1)
Y^{-\nu_1-\nu_2}\cr
&-(\nu_1-\nu_2)
\varphi_\varGamma(-\nu_2)Y^{\nu_1+\nu_2}
+(\nu_1+\nu_2)\varphi_\varGamma(\nu_1)
\varphi_\varGamma(-\nu_2)Y^{-\nu_1+\nu_2}\cr
&+\int_{-{1\over2}}^{1\over2}
\left[\tilde{e}_0(z,\nu_1){\partial\tilde{e}_0
\over{\partial{y}}}(z,-\nu_2)
-\tilde{e}_0(z,-\nu_2){\partial\tilde{e}_0\over{\partial{y}}}
(z,\nu_1)\right]_{y=Y}dx\Bigg\}.\qquad
&(39.10)
}
$$ 
On noting that $\Omega_0\tilde{e}_0(z,\nu)
=\big({1\over4}-\nu^2\big)\tilde{e}_0(z,\nu)$,
another application of $(1.9)$ leads us to $(39.8)$. Thus the
first term on the right of $(39.7)$ is
$$
\eqalignno{
{1\over4\pi i}\int_{(0)}&k(\nu)
\bigg(2\log Y-{\varphi'_\varGamma
\over\varphi_\varGamma}(\nu)
-\int_{{\eusm{F}}\backslash{\eusm{F}}_Y}
|\tilde{e}_0(z,\nu)|^2d\mu(z)\bigg)d\nu\cr 
&+{1\over4\pi i}\int_{(0)} \left(Y^{2\nu}
\varphi_\varGamma(-\nu)-1\right)k(\nu){d\nu\over\nu},
&(39.11)
}
$$
in which the second line is a result of removing the superfluous
singularity at $\nu_1+\nu_2=0$ in the expression $(39.8)$.
The last integral is seen to be ${1\over4}k(0)+O(Y^{-1/2})$ by
shifting the contour to $\Re\nu=-{1\over4}$ and then back to
the original. Thus we have
$$
\eqalignno{
\int_{{\eusm F}_Y}&{}^cG(z,z)d\mu(z)
=\left({1\over2\alpha-1}-{1\over2\beta-1}\right)\log Y\cr
&-{1\over4\pi i}\int_{-\infty}^\infty k(\nu)
{\varphi'_\varGamma
\over\varphi_\varGamma}(\nu)d\nu
+k\big(\txt{1\over2}\big)+{1\over4}k(0)
+O\big(Y^{-{1/2}}\big).&(39.12)
}
$$
Combined with $(39.6)$, this gives that 
$$
\eqalignno{
&\lim_{Y\to\infty}\int_{{\eusm F}_Y}[R^{(1)}-{}^cG](z,z)
d\mu(z)={1\over4\pi i}\int_{-\infty}^\infty k(\nu)
{\varphi'_\varGamma
\over\varphi_\varGamma}(\nu)d\nu-k\big(\txt{1\over2} \big)\cr
&\quad-{1\over2\alpha-1}\left({\Gamma'\over\Gamma}
(\alpha+\txt{1\over2})+\log2\right)
+{1\over2\beta-1}\left({\Gamma'\over\Gamma}
(\beta+\txt{1\over2})+\log2\right).&(39.13)
}
$$
The first term on the right side equals
$$
\eqalignno{
&-{1\over2\pi i}\int_{-\infty}^\infty k(\nu)
\left(2{\zeta'\over\zeta}(1+2\nu)+{1\over \nu}+
{\Gamma'\over\Gamma}
\big(\txt{1\over2}+\nu\big)-\log\pi\right)d\nu\cr
=&-{1\over2\alpha-1}
\left(2{\zeta'\over\zeta}(2\alpha)+{2\over2\alpha-1}+
{\Gamma'\over\Gamma}(\alpha)-\log\pi\right)\cr
&+{1\over2\beta-1}
\left(2{\zeta'\over\zeta}(2\beta)+{2\over2\beta-1}+
{\Gamma'\over\Gamma}(\beta)-\log\pi\right).&(39.14)
}
$$
The insertion of the term $1/\nu$ on the left side is to
 remove the superfluous singularity at $\nu=0$; and
the right side is the result of shifting the contour to $+\infty$.
\medskip
\noindent
{\csc Notes:} In literature the formula 
$(39.8)$ is termed the Maass--Selberg identity. 
\medskip
\noindent
{\bf 40. Hyperbolic terms.} In order to deal with $R^{(2)}$,
we need a preparation: We shall introduce a decomposition
of ${\cal{C}}^{(2)}$ which corresponds to $(39.1)$. Thus, let
$\omega$ be an arbitrary element in ${\cal C}^{(2)}$, and let us
suppose that $\r{g}_0\in\r{G}$ maps the cusp and the origin to the
two fixed points of $\omega$. Then there exists a positive
constant $c(\omega)\ne1$ such that
$\r{g}_0^{-1}\omega\r{g}_0(z)=c(\omega)z$ for
$z\in\B{H}^2$. In particular, those hyperbolic elements
$\varGamma$ which share the same pair of fixed points
make up an infinite cyclic group. 
Hence, we may classify the elements
of ${\cal C}^{(2)}$ in terms of their fixed points:
$$
{\cal C}^{(2)}=\bigsqcup_\omega 
\Big\{[\omega]\backslash\{1\}\Big\},\quad [\omega]=
\{\omega^n:\, n\in\B{Z}\},\eqno(40.1)
$$
where $\{\omega\}$ is a representative set. Further, we classify
subgroups $[\omega]$ according to the 
$\varGamma$-conjugacy:
$$
{\cal C}^{(2)}=\bigsqcup_{\omega_0}
\bigsqcup_{\gamma\in B(\omega_0)\backslash\varGamma}
\gamma^{-1}\Big\{[\omega_0]\backslash\{1\}\Big\}\gamma,
\eqno(40.2)
$$
where $B(\omega_0)$ is the normaliser of 
$[\omega_0]$ in $\varGamma$. Here it should be noted that
the centraliser of any hyperbolic $\omega\ne1$ is
$[\omega]$. In fact, if $\gamma\omega
=\omega\gamma$, then $\gamma,\, \omega$ share the
same fixed points, and $\gamma$ is a power of $\omega$; 
to confirm this, one may use the relation $\r{g}_0^{-1}
\gamma\r{g}_0(c(\omega)z)=c(\omega)
\r{g}_0^{-1}\gamma\r{g}_0(z)$ with $\r{g}_0$ as above.
On the other hand, if $\xi\in B(\omega_0)$, then
$\xi^{-1}\omega_0\xi$ is a generator $[\omega_0]$. Thus,
there are two cases in general: Either there exists a 
$\delta\in B(\omega_0)$ such that 
$\delta^{-1}\omega_0\delta=\omega_0^{-1}$ or there does not.
In the former case, $\delta$ is an involution exchanging 
the fixed points of $\omega_0$. Hence
$[B(\omega_0):[\omega_0]]=2$, and
$$
\bigsqcup_{\gamma\in B(\omega_0)
\backslash\varGamma}\gamma^{-1}
\Big\{[\omega_0]\backslash\{1\}\Big\}\gamma
=\bigsqcup_{\gamma\in [\omega_0]\backslash\varGamma}
\big\{\gamma^{-1}\omega_0^n\gamma: n\in
\B{N}\big\}.\eqno(40.3)
$$
In the latter case, we have $B(\omega_0)=[\omega_0]$, and
$$
\bigsqcup_{\gamma\in B(\omega_0)
\backslash\varGamma}\gamma^{-1}
\Big\{[\omega_0]\backslash\{1\}\Big\}\gamma
=\bigsqcup_{\pm}\bigsqcup_{\gamma\in [\omega_0]
\backslash\varGamma}\big\{\gamma^{-1}\omega_0^{\pm n}
\gamma:{\Bbb Z}\ni n>0\big\}.\eqno(40.4)
$$ 
For example, $\omega=\left[{2\atop1}{1\atop1}\right]$ 
belongs to the former case, since
$\delta=\left[{\atop1}{-1\atop}\right]$ gives $\delta
\omega\delta^{-1}=\omega^{-1}$. On the other hand,
$\omega=\left[{2\atop1}{3\atop2}\right]$ belongs to the
latter case. Hence, we classify the generators
$\big\{\omega_0^{\pm1}\big\}$ according to the 
$\varGamma$-conjugacy, and designate
$\{\varpi\}$ as a representative set:
$$
{\cal C}^{(2)}=\bigsqcup_{\varpi}
\bigsqcup_{\gamma\in [\varpi]\backslash\varGamma}
\big\{\gamma^{-1}\varpi^n\gamma:
{\Bbb Z}\ni n>0\big\}.\eqno(40.5)
$$
Each class containing $\varpi$ is termed a prime class.
\par
Then, we introduce the norm of a hyperbolic
element $\omega\ne1$: With $c(\omega)$ as above,
$$
N(\omega)=\max\big(c(\omega),c(\omega^{-1})\big)>1,
\eqno(40.6)
$$
which is a function of the conjugacy class to which $\omega$
belongs. Using the hyperbolic distance $(2.1)$, we have
$$
N(\omega)=\exp\Big(\inf_{z\in{\eusm H}}
d(z,\omega(z))\Big),\eqno(40.7)
$$
since
$$
\eqalignno{
d\big(z,\omega(z)\big)=\,&d\big(\r{g}_0^{-1}(z),
c(\omega)\r{g}_0^{-1}(z)\big)\cr
=\,&2\,{\rm arcsinh}\left({\big|\r{g}_0^{-1}(z)-c(\omega)
\r{g}_0^{-1}(z)\big|\over2\sqrt{c(\omega)}\,
\Im\r{g}_0^{-1}(\omega)}\right)\cr
\ge\,& 2\,{\rm arcsinh}\Big(
\txt{1\over2}\big|c(\omega)^{1/2}
-c(\omega)^{-{1/2}}\big|\big)=\log N(\omega),&(40.8)
}
$$
in which the minimum is attained with $z$ such that
$\r{g}_0^{-1}(z)\in i\B{R}$. We shall use the term
$$
\hbox{$N(\varpi)$: a pseudo-prime.}\eqno(40.9)
$$
The reason why we liken $N(\varpi)$ to a prime number will be
revealed in Section 42. The distance $d(z,\omega(z))$
is the hyperbolic length of the closed geodesic 
on the Riemann surface $\varGamma\backslash\B{H}^2$
which starts at $z$ and returns to $z$. Hence $\log N(\varpi)$ is
the length of the shortest closed geodesic among those
associated with the prime class represented by $\varpi$.
\par
With this, we let $R^{(2)}_\alpha$ be the part of $R^{(2)}$
corresponding to $p_\alpha$, similarly to $(39.2)$. Then we have,
by $(40.5)$,
$$
\eqalign{
\int_{\eusm F}R_\alpha^{(2)}(z,z)d\mu(z)&=
\int_\e{F}\sum_\varpi\sum_{\gamma\in[\varpi]\backslash
\varGamma}\,\sum_{n=1}^\infty p_\alpha
\big(\varrho(z,\gamma^{-1}\varpi^n\gamma(z))\big)d\mu(z)\cr
&=\sum_\varpi\sum_{n=1}^\infty 
\int_{[[\varpi]]}p_\alpha\big(\varrho(z,\varpi^nz)\big)d\mu(z),
}\eqno(40.10)
$$
where $[[\varpi]]=\bigsqcup_{\gamma\in[\varpi]\backslash
\varGamma}\gamma\e{F}=[\varpi]\backslash\B{H}^2$.
We choose an $\eta\in\r{G}$ such that
$\eta\varpi\eta^{-1}(z)=N(\varpi)z$. Then
$\eta[[\varpi]]=[[\eta\varpi\eta^{-1}]]$, and we have
$$
\eqalignno{
\int_{[[\varpi]]}p_\alpha
\big(\varrho(z,\varpi^n(z))\big)d\mu(z)
&=\int_{[[\eta\varpi\eta^{-1}]]}
p_\alpha\big(\varrho(z,\eta\varpi^n\eta^{-1}(z))\big)
d\mu(z)\cr&=\mathop{\int\!\int}_{\scr{1\le|z|\le N(\varpi)}\atop
\scr{z\in\B{H}^2}}p_\alpha
\big(\varrho(z, N(\varpi)^nz)\big){dxdy\over y^2}.&(40.11)
}
$$
For any $a,b>1$, we have, by $(34.6)$ for $\ell=0$,
$$
\eqalignno{
&\mathop{\int\!\int}_{\scr{1\le|z|\le a}\atop\scr{0<\Im z}}
p_\alpha\big(\varrho(z,bz)\big){dxdy\over y^2}
=2\log a\,\int_0^{{1\over2}\pi} p_\alpha\left(
{1\over4\sin^2\theta}\big(b^{1/2}-b^{-{1/2}}\big)^2
\right){d\theta\over\sin^2\theta}\cr
&={\log a\over4\pi^2i}\int_{\left(-{1\over2}\right)}
{\Gamma^2(\xi+\alpha)\Gamma(-\xi)\over\Gamma(\xi+2\alpha)}
\big(\txt{1\over2}(b^{1/2}-b^{-{1/2}})\big)^{-2(\xi+\alpha)}
\int_0^{{1\over2}\pi}(\sin\theta)^{2(\xi+\alpha-1)}d\theta d\xi\cr 
&={\log a\over8\pi^{3/2}i}\int_{\left(-{1\over2}\right)}
{\Gamma(\xi+\alpha)\Gamma\big(\xi+\alpha-{1\over2}\big)
\Gamma(-\xi)\over\Gamma(\xi+2\alpha)}
\big(\txt{1\over2}\big(b^{1/2}-b^{-{1/2}}\big)
\big)^{-2(\xi+\alpha)}d\xi\cr
&={\log a\over 4\pi^{1/2}}{\Gamma\big(\alpha-{1\over2}\big)
\Gamma(\alpha)\over\Gamma(2\alpha)}
\big(\txt{1\over2}\big(b^{1/2}-b^{-{1/2}}\big)\big)^{-2\alpha}
{}_2F_1\left(\alpha-\txt{1\over2},\alpha;2\alpha;
-\big(\txt{1\over2}\big(b^{1/2}-b^{-1/2}\big)\big)^2\right)
\qquad\cr
&={b^{{1/2}-\alpha}\log a\over(2\alpha-1)
(b^{1/2}-b^{-{1/2}})}.&(40.12)
}
$$
Hence we have that
$$
\eqalignno{
\int_{\eusm F}R^{(2)}(z,z)d\mu(z)
&={1\over2\alpha-1}\sum_{\varpi}\sum_{n=0}^\infty
{\log N(\varpi)\over N(\varpi)^{\alpha+n}-1}\cr
&-{1\over2\beta-1}\sum_{\varpi}\sum_{n=0}^\infty
{\log N(\varpi)\over N(\varpi)^{\beta+n}-1}.&(40.13)
}
$$
We stress that
the double sums are absolutely and uniformly convergent
for $\Re\alpha>1$ and $\Re\beta>1$, respectively, and
represent regular functions. In fact, by the same token
as the remark made immediately after $(38.6)$ and
by the observation that
$(34.15)$ for $\ell=0$ implies that
$p_\alpha(\varrho)>0$ for $\alpha>0$,
the series defining $R^{(2)}_\alpha(z,z)$
is absolutely and uniformly convergent and bounded
for $\alpha>1$ and $z\in\e{F}$.
\medskip
\noindent
{\csc Notes:} The assertion $(40.12)$ depends on 
quadratic transformations of the
Gaussian hypergeometric function, an account of
which is given in [18, Section 9.6]. The particular
formula we need here is 
$$
{}_2F_1\big(\alpha-\txt{1\over2},\alpha;2\alpha;z\big)=
\left(\txt{1\over2}(1+\sqrt{1-z})\right)^{1-2\alpha},
\quad|\arg(1-z)|<\pi.\eqno(40.14)
$$
See $(9.8.3)$ loc.\ cit.
\medskip
\noindent
{\bf 41. Elliptic terms.} The fixed points of elliptic elements
in $\varGamma$ are $\varGamma$-images of the points
$\exp\big({1\over2}\pi i\big)$ and 
$\exp\big({2\over3}\pi i\big)$. The subgroups consisting 
of elements which fix these points are the cyclic groups
 $[\omega_2]$ and $[\omega_3]$ with
$\omega_2=\left[{\atop1}
{-1\atop}\right]$ and $\omega_3=
\left[{\atop1}{-1\atop\phantom{-}1}\right]$, respectively.
Also, their normalisers are themselves. Thus we have
$$
\eqalignno{
\int_{\eusm F}R^{(e)}(z,z)d\mu(z)
=&\int_{[[\omega_2]]}(p_\alpha-p_\beta)
\big(\varrho(z,\omega_2(z))\big)d\mu(z)\cr
+&\sum_{j=1,2}\int_{[[\omega_3]]}
(p_\alpha-p_\beta)\big(\varrho(z,\omega^j_3(z))
\big)d\mu(z).&(41.1)
}
$$
To compute the right side, we shall consider, more
generally, the following expression
$$
T_\alpha(\omega^j)=\int_{[[\omega]]}
p_\alpha\big(\varrho(z,\omega^j(z))\big)
d\mu(z),\qquad1\le j\le m-1,
\eqno(41.2)
$$
where $\omega$ is elliptic and of order $m\ge2$. We have
$$
\int_{\eusm F}R^{(e)}(z,z)d\mu(z)
=T_\alpha(\omega_2)-T_\beta(\omega_2)
+\sum_{j=1,2}\big(T_\alpha(\omega^j_3)
-T_\beta(\omega^j_3)\big).\eqno(41.3)
$$
In $(41.2)$ we replace the domain
$[[\omega]]$ by $\omega^\nu[[\omega]]$,
$\nu\in\B{Z}$, and get
$$
T_\alpha(\omega^j)
={1\over m}\int_{\B{H}^2}p_\alpha
\big(\varrho(z,\omega^j(z))\big)d\mu(z).\eqno(41.4)
$$
Exploiting the $\varGamma$-conjugation, we may assume that
the fixed point of $\omega$ is $\exp\big({1\over2}\pi i\big)$,
that is, $\omega^j=\left[{\eta\atop\lambda}
{-\lambda\atop\phantom{-}\eta}\right]$ with
$\eta=\cos((j/m)\pi)$, $\lambda=\sin((j/m)\pi))$.
Hence
$$
T_\alpha(\omega^j)={1\over m}\int_{-\infty}^\infty\int_0^\infty
p_\alpha\left(\txt{1\over4}\lambda^2y^{-2}|z^2
+1|^2\right){dydx\over y^2}.\eqno(41.5)
$$
We perform the change of variable $y\mapsto(x^2+1)^{1/2}u$,
and divide the inner integral at $u=1$, and then apply the
change of variable $u\mapsto1/u$  in the integral 
over the unit interval. Further, we apply $v=u-1/u$ and
$v\mapsto 2|x|(x^2+1)^{-{1/2}}w^{1/2}$.  The
inner integral is transformed into
$$
{|x|\over x^2+1}\int_0^\infty 
p_\alpha\big((\lambda x)^2(w+1)\big)
w^{-{1/2}}dw.\eqno(41.6)
$$
This is computed in much the same as $(40.12)$, and we have
$$
\eqalignno{
T_\alpha(\omega^j)&={1\over (2\alpha-1)\lambda m}\int_0^\infty
{\big(\lambda x+\sqrt{(\lambda x)^2+1}\,\big)^{1-2\alpha}
\over x^2+1}dx\cr
&={2\over(2\alpha-1)m}\int_0^\infty
{\big(\max(\xi,1/\xi)\big)^{1-2\alpha}\over(\xi-1/\xi)^2
+4\lambda^2}d\xi.&(41.7)
}
$$
The second line is due to the transformation
$2\lambda x\mapsto \xi-1/\xi$. Invoking the representation
$$
{1\over2\alpha-1}\big(\max(\xi,1/\xi)\big)^{1-2\alpha}
={1\over2\pi}\int_{-\infty}^\infty
{\xi^{2it}\over{\big(\alpha-{1\over2}\big)^2+t^2}}dt,
\eqno(41.8)
$$
we have that
$$
\eqalignno{
T_\alpha(\omega^j)
&={1\over\pi m}\int_{-\infty}^\infty 
{1\over\big(\alpha-{1\over2}\big)^2+t^2}
\int_0^\infty{\xi^{2it}\over
\big(\xi-1/\xi\big)^2+4\lambda^2}d\xi dt\cr
&={1\over\pi m}\int_{-\infty}^\infty{1\over
\big((\alpha-{1\over2})^2+t^2\big)(1+e^{-2\pi t})}
\int_{-\infty}^\infty{\xi^{2it}\over\big(\xi-1/\xi\big)^2
+4\lambda^2}d\xi dt.
&(41.9)
}
$$
In the last integral we set $\arg\xi=0$ for $\xi>0$ and
$\arg\xi=\pi$ for $\xi<0$. Shifting the path to $+i\infty$,
$$
T_\alpha(\omega^j)={1\over2m\lambda}\int_{-\infty}^\infty
{e^{-2(j/m)\pi t}\over\big((\alpha-{1\over2})^2+t^2\big)
(1+e^{-2\pi t})}dt.\eqno(41.10)
$$
Again shifting the path to $+i\infty$,
$$
\eqalignno{
T_\alpha(\omega^j)&={\pi\over(2\alpha-1)m\lambda}
{e^{(1-2\alpha)(j/m)\pi i}\over1-e^{-2\alpha\pi i}}\cr
&+{i\over2(2\alpha-1)m\lambda}\sum_{l=0}^\infty
e^{-(2l+1)(j/m)\pi i}
\left({1\over\alpha+l}-{1\over1-\alpha+l}\right),&(41.11)
}
$$
which gives a meromorphic continuation with respect to $\alpha$.
The sum is
$$
\eqalignno{
&{1\over m}\sum_{l=0}^{m-1}e^{-(2l+1)(j/m)\pi i}
\sum_{f=0}^\infty\left({1\over(\alpha+l)/m+f}
-{1\over(1-\alpha+l)/m+f}\right)\cr
=&{1\over m}\sum_{l=0}^{m-1}e^{-(2l+1)(j/m)\pi i}
\left({\Gamma'\over\Gamma}\big((1-\alpha+l)/m\big)
-{\Gamma'\over\Gamma}\big((\alpha+l)/m\big)\right)\cr
=&{1\over m}\sum_{l=0}^{m-1}\left(
e^{(2l+1)(j/m)\pi i}{\Gamma'\over\Gamma}
\big(1-(\alpha+l)/m\big)-e^{-(2l+1)(j/m)\pi i}
{\Gamma'\over\Gamma}\big((\alpha+l)/m\big)\right)\qquad\cr
=&{2i\over m}\sum_{l=0}^{m-1}\left(\sin\big((2l+1)(j/m)\pi\big)
{\Gamma'\over\Gamma}\big((\alpha+l)/m\big)+{\pi\over 2i}
{e^{(2l+1)(j/m)\pi i}\over\tan\big(\pi(\alpha+l)/m\big)}
\right),&(41.12)
}
$$
where the second line depends on Weierstrass' product representation
of the $\Gamma$-function and the fourth line on the
functional equation $\Gamma(s)\Gamma(1-s)=\pi/\sin(\pi s)$.
As to the first term on the right side of $(41.11)$, we have
$$
{e^{(1-2\alpha)(j/m)\pi i}\over1-e^{-2\alpha\pi i}}
={1\over 2im}\sum_{l=0}^{m-1}{e^{(2l+1)(j/m)\pi i}
\over\tan\big(\pi(\alpha+l)/m\big)}.\eqno(41.13)
$$
In fact, the difference of the two sides is entire, since they have
the same set of poles and respective principal parts are identical.
Noting the periodicity and
taking $\alpha$ to $i\infty$, we get the assertion.
Hence we have that
$$
T_\alpha(\omega^j)={1\over(1-2\alpha)m^2}
\sum_{l=0}^{m-1}{\sin\big((2l+1)(j/m)\pi\big)
\over\sin\big((j/m)\pi\big)}
{\Gamma'\over\Gamma}\big((\alpha+l)/m\big),\eqno(41.14)
$$
and
$$
\sum_{j=1}^{m-1}T_\alpha(\omega^j)
={1\over (1-2\alpha)m^2}\sum_{l=0}^{m-1}(m-2l-1)
{\Gamma'\over\Gamma}\big((\alpha+l)/m\big),\eqno(41.15)
$$
since
$$
\sum_{j=1}^{m-1}{\sin\big((2l+1)j\pi/m\big)
\over\sin\big((j/m)\pi\big)}
=\sum_{j=1}^{m-1}\sum_{h=-l}^le\big((j/m)h\big)=m-2l-1.
\eqno(41.16)
$$
\medskip
\noindent
{\bf 42. Selberg's zeta-function and trace formula.} We sum up the
discussion developed in the last four sections. 
Thus we classify hyperbolic elements of 
$\varGamma$ as in $(40.5)$, which defines the prime classes 
$\{\varpi\}$. The corresponding pseudo-primes
$N(\varpi)$ are defined by $(40.9)$. 
\par
With this, we introduce the Selberg zeta-function associated with
$\varGamma$ by
$$
\zeta_\varGamma(s)=\prod_\varpi\prod_{n=0}^\infty
\bigg(1-{1\over N(\varpi)^{s+n}}\bigg).\eqno(42.1)
$$
We have shown already that
$$
{\zeta_\varGamma'\over\zeta_\varGamma}(s)
=\sum_\varpi\sum_{n=0}^\infty
{\log N(\varpi)\over N(\varpi)^{s+n}-1},
\quad \Re s>1,\eqno(42.2)
$$
is absolutely convergent. We introduce also
$$
\eqalignno{
\Psi_\varGamma(s)&
={\zeta_\varGamma'\over\zeta_\varGamma}(s)
-2{\zeta'\over\zeta}(2s)
-{1\over s-1}-{1\over s}-{2\over 2s-1}\cr
&+\log(2\pi)-2{\Gamma'\over\Gamma}(2s)+W(s),
&(42.3)
}
$$
with
$$
\eqalignno{
W(s)=&-{1\over6}(2s-1)
{\Gamma'\over\Gamma}(s)
+{1\over4}{\Gamma'\over\Gamma}
\big(\txt{1\over2}(s+1)\big)
-{1\over4}{\Gamma'\over\Gamma}
\big(\txt{1\over2}s\big)\cr
&+{2\over9}{\Gamma'\over\Gamma}
\big(\txt{1\over3}(s+2)\big)
-{2\over9}{\Gamma'\over\Gamma}\big(\txt{1\over3}s\big).
&(42.4)
}
$$

Collecting $(38.7)$, $(38.8)$, $(39.13)$, $(39.14)$,
$(40.13)$, $(41.3)$, $(41.15)$ as well as invoking the
duplication formula for the $\Gamma$-function, we obtain
a version of
the Selberg trace formula for the full modular group $\varGamma$:
The function $\Psi_\varGamma$ exists 
as a meromorphic function over the entire complex plane, and
we have, for arbitrary $\alpha,\beta\in\B{C}$,
$$
{\Psi_\varGamma(\alpha)\over2\alpha-1}
-{\Psi_\varGamma(\beta)\over2\beta-1}=
\mathop{{\sum}^\r{u}}_{V\hfil}
\left\{{1\over\big(\alpha-{1\over2}\big)^2
-\nu_V^2}-{1\over\big(\beta-{1\over2}\big)^2
-\nu_V^2}\right\}.\eqno(42.5)
$$
In fact, we have already established this for
$1<\Re\alpha<\Re\beta$. Since this implies in particular that
we have $\sum_V^{\r{u}}|\nu_V|^{-4}<+\infty$, the right
side is a meromorphic function over $\B{C}^2$, 
by analytic continuation. One may put $\alpha=s$ and 
$\beta=2$, and find that $\Psi_\varGamma(s)$ is regular 
for all $s$ save for 
the simple poles $\big\{{1\over2}\pm \nu_V: 
\hbox{$V$ in the unitary
principal series}\big\}$.
\smallskip
The definition $(42.1)$ conjures up the 
Euler product representation for the Riemann zeta-function.
Thus we shall look for analogies between the
two zeta-functions. First we observe that $W(s)$ has poles
at non-positive integers; the residues at $s=0$ is equal to $1$, and
all other poles have negative integers as their residues, which can be
seen by classifying poles of the expression $(42.4)$ 
according to $\bmod\, 6$. Hence the poles of
$(\zeta'_\varGamma/\zeta_\varGamma)(s)$ are all of order $1$,
and their residues are equal to integers. This means that
$\zeta_\varGamma(s)$ exists as a meromorphic function over
$\B{C}$. To be more precise, the function
$$
\xi_\varGamma(s)={(2\pi)^s\zeta_\varGamma(s)\over\zeta(2s)
\Gamma(2s)s(s-1)(2s-1)}
\exp\bigg(\int_{0}^{s-{1\over2}} W\big(\eta+\txt{1\over2}\big)
d\eta\bigg)\eqno(42.6)
$$
is entire, and its zeros are ${1\over2}\pm\nu_V$,
with $V$ as above. In fact, we have
$$
{\xi_\varGamma'\over\xi_\varGamma}(s)
=\Psi_\varGamma(s).\eqno(42.7)
$$
Taking account of the nature of the poles of $W(s)$, the function
$$
\zeta_\varGamma(s)\over\zeta(2s)\Gamma(2s)(2s-1)
\eqno(42.8)
$$
is entire. Its real zeros are precisely at
$1$, which is simple, and at negative integers; and the set of
its complex zeros coincides with that of $\xi_\varGamma(s)$.
We have, in a neighbourhood of $s=1$,
$$
{\zeta_\varGamma'\over\zeta_\varGamma}(s)={1\over s-1}
+O(1).\eqno(42.9)
$$
Further, the functional equation
$$
\zeta_\varGamma(s)
=\chi_p(s)\chi_e(s)\chi_1(s)\zeta_\varGamma(1-s)
\eqno(42.10)
$$
holds, where
$$
\eqalign{
\chi_p(s)&=-{(2\pi)^{1-2s}\zeta(2s)\Gamma(2s)\over
\zeta(2(1-s))\Gamma(2(1-s))},\cr
\chi_e(s)&=\bigg({\sin{1\over2}\pi(1-s)
\over\sin{1\over2}\pi s}\bigg)^{1/2}
\bigg({\sin{1\over3}\pi(1-s)\over
\sin{1\over3}\pi s}\bigg)^{2/3},\cr
\chi_1(s)&=\exp\bigg(\txt{1\over3}\pi
\int_0^{s-{1\over2}}\eta\tan(\pi\eta)d\eta\bigg).
}\eqno(42.11)
$$
The factors $\chi_p$, $\chi_e$, $\chi_1$ come from
parabolic, elliptic elements, and the unit, respectively. The function
$\chi_e(s)\chi_1(s)$ is meromorphic taking the value $1$ at
$s={1\over2}$. In fact, the right side of the
trace formula $(42.5)$ is
invariant against the transform $\alpha\mapsto 1-\alpha$, and
we have
$$
\Psi_\varGamma(s)=-\Psi_\varGamma(1-s).\eqno(42.12)
$$
Integrating this, we have
$$
\xi_\varGamma(s)=\xi_\varGamma(1-s).\eqno(42.13)
$$
Also, we have
$$
\eqalignno{
W&\big(\txt{1\over2}+\eta\big)
+W\big(\txt{1\over2}-\eta\big)=
{\eta\over3}{d\over d\eta}
\log\,\sin\big(\pi(\txt{1\over2}+\eta)\big)\cr
&-{1\over2}{d\over d\eta}\log{\sin\big(\txt{1\over2}
\pi(\txt{1\over2}-\eta)\big)\over\sin
\big(\txt{1\over2}\pi(\txt{1\over2}+\eta)\big)}
-{2\over3}{d\over d\eta}\log{\sin\big(\txt{1\over3}
\pi(\txt{1\over2}-\eta)\big)\over\sin
\big(\txt{1\over3}\pi(\txt{1\over2}+\eta)\big)}.&(42.14)
}
$$
In this way we obtain $(42.10)$.
\medskip
\noindent
{\csc Notes:} The assertion in this section is a typical
instance of Selberg's fundamental discovery [32]. 
The similarities between $\zeta(s)$ and 
$\zeta_\varGamma(s)$ are striking. However, they are
in fact ostensible. Comparing the original 
and the pseudo-Euler products, one surmise that the
correct analogy should be found between $1/\zeta(s)$ and
$\zeta_\varGamma(s)$. Then we find grave differences.
The most salient is: Despite the fact that both 
have a simple zero at
$s=1$, the former has infinitely many poles
in the critical strip $0<\Re s<1$ and the latter has none at all.
Nevertheless, see the notes to the next section for another
view point. It should be added that strictly speaking
the analogue of the Riemann
hypothesis holds with $\xi_\varGamma(s)$ but not with
$\zeta_\varGamma(s)$ itself as the latter has complex zeros
coming from $\zeta(2s)$. It should be observed that the
factor $\zeta(2s)$ is due to parabolic elements. Hence, in
the case where the underlying discrete subgroup has a compact
fundamental domain unlike the full modular group, the corresponding
zeta-function of Selberg satisfies precisely 
the analogue of the Riemann hypothesis, as it is known that
such a discrete group lacks parabolic elements. 
\medskip
\noindent
{\bf 43. Weyl's law.} We shall prove the asymptotic formula:
$$
N_\varGamma(K)=\mathop{{\sum}^\r{u}}_{|\nu_V|\le K}1=
{1\over 12}K^2+O\big(K\log K\big).\eqno(43.1)
$$
This is the same as counting complex zeros of
$\xi_\varGamma(s)$ in the region $0<\Re s<1$,
$0<\Im s\le K$. We shall employ an argument from the
theory of the Riemann zeta-function. 
\par
First, we put
$\alpha={1\over2}+K$ and 
$\beta={1\over2}+2K$ in $(42.5)$ so that
$$
\eqalignno{
{1\over K^2}\mathop{{\sum}^\r{u}}_{K\le|\nu_V|\le 2K}1&\ll
\mathop{{\sum}^\r{u}}_{V\hfil} \bigg\{{1\over K^2+|\nu_V|^2}
-{1\over 4K^2+|\nu_V|^2}\bigg\}\cr
&={1\over2K}\Psi_\varGamma\left(\txt{1\over2}+K\right)
-{1\over 4K}\Psi_\varGamma
\left(\txt{1\over2}+2K\right).&(43.2)
}
$$
By the definition $(42.3)$--$(42.4)$ and an asymptotic
expansion of $(\Gamma'/\Gamma)(s)$, we have
$$
N_\varGamma(K)\ll K^2.\eqno(43.3)
$$
Also putting $\alpha={3\over2}+iK$ and 
$\beta={1\over2}+K$ in $(42.5)$, we have
$$
\mathop{{\sum}^\r{u}}_{V\hfil}\bigg\{
{1\over (1+iK)^2+|\nu_V|^2}-
{1\over K^2+|\nu_V|^2}\bigg\}=O(1).\eqno(43.4)
$$
Applying $(43.3)$ to this, we get
$$
\mathop{{\sum}^\r{u}}_{||\nu_V|-K|\le{1\over2}K}
{1\over (1+iK)^2+|\nu_V|^2}=O(1),\eqno(43.5)
$$
and thus
$$
\mathop{{\sum}^\r{u}}_{||\nu_V|-K|\le{1\over2}K}
{1\over 1+i(K-|\nu_V|)}=O(K).\eqno(43.6)
$$
Hence $(43.3)$ is improved to
$$
N_\varGamma(K+1)-N_\varGamma(K)\ll K.\eqno(43.7)
$$
This implies in turn that for any bounded $s$
$$
\Psi_\varGamma(s)=
\mathop{{\sum}^\r{u}}_{||\nu_V|-t|\le{1\over2}}
{1\over s-{1\over2}-\nu_V}
+O\big(t\log 2t\big),\quad \Im s=t\ge1.\eqno(43.8)
$$
\par
We may assume hereafter that the parameter
$K$ satisfies
$$
\inf_V||\nu_V|-K|\gg K^{-1}.\eqno(43.9)
$$
Then, by the functional equation $(42.13)$, we have
$$
N_\varGamma(K)={1\over\pi}\arg\xi_\varGamma
\!\left(\txt{1\over2}+iK\right).\eqno(43.10)
$$
The argument starts with the value $0$ at $+\infty+iK$ and
continuously varies along the line $\Im s=K$. By $(42.6)$ we have
$$
\eqalignno{
N_\varGamma(K)&={1\over\pi}\Im
\int_0^{iK}W\big(\eta+\txt{1\over2}\big)d\eta+
{1\over\pi}\arg\zeta_\varGamma\!\left(\txt{1\over2}+iK\right)
+O\big(K\log K\big)\cr&={1\over12}K^2+
{1\over\pi}\arg\zeta_\varGamma\!\left(\txt{1\over2}+iK\right)
+O\big(K\log K\big),&(43.11)
}
$$
in which the error term is due to the $\Gamma$-factor and
$\zeta(2s)$. The argument of
$\zeta_\varGamma(s)$ is defined in the same way
as that of $\xi_\varGamma(s)$, and we consider the
estimation of
$$
\eqalignno{
\log\zeta_\varGamma\left(\txt{1\over2}+iK\right)&
=\int_\infty^{1\over2}
{\zeta_\varGamma'\over\zeta_\varGamma}(u+iK)du\cr
&=\int_{3\over2}^{1\over2}
{\zeta_\varGamma'\over\zeta_\varGamma}(u+iK)du
+O(1).&(43.12)
}
$$
Here the left side is defined to be the result of
the analytic continuation of the integral of $(\zeta'_\varGamma/
\zeta_\varGamma)(s)$ along the above
horizontal line. By the definition $(42.3)$ and the
approximation $(43.8)$, we have
$$
{\zeta_\varGamma'\over\zeta_\varGamma}(u+iK)
=\mathop{{\sum}^\r{u}}_{||\nu_V|-K|\le{1\over2}}
{1\over u-{1\over2}+i(K-|\nu_V|)}+O\big(K\log K\big),
\eqno(43.13)
$$
where we have used a well-known bound of $(\zeta'/\zeta)(s)$
for $\Re s\ge1$. Hence, in view of $(43.7)$ we obtain
$$
\arg\zeta_\varGamma(\txt{1\over2}+iK)\ll K\log K.\eqno(43.14)
$$
We end the proof of $(43.1)$.
\medskip
\noindent
{\csc Notes:} The function $\Xi_\varGamma(t)=
\xi_\varGamma\big({1\over2}+it\big)$ is an analogue of
the $\Xi$-function in the theory of $\zeta(s)$. This is real
for real $t$. It oscillates wildly as $(43.1)$ implies; in fact the
the number of real zeros in the interval $|t-K|\le B\log K$ is
more than $K\log K$, provided $B$ is sufficiently large.
Still all zeros of $\Xi_\varGamma(t)$ are on the real axis
without any single exception. Amazing, indeed.
\medskip
\noindent
{\bf 44. Pseudo-prime number theorem.} We shall prove
the asymptotic formula
$$
\pi_\varGamma(x)=\sum_{\hbox{$\scr\varGamma$\ssrm-pseudo-prime$\,\scr{<x}$}}1
=\int_2^x {du\over\log u}
+O\left(x^{3/4}(\log x)^{-1/2}\right).\eqno(44.1)
$$
Following the traditional treatment of the distribution
of prime numbers, we consider
$$
\psi_\varGamma(x)=\sum_{n=1}^\infty
\sum_{N(\varpi)<x^{1/n}}\log N(\varpi).
\eqno(44.2)
$$
The relevant generating function is $(\zeta'_\varGamma/
\zeta_\varGamma)(s)-(\zeta'_\varGamma/
\zeta_\varGamma)(s+1)$. We are, however, unable to
apply Perron's inversion formula, because of the difficulty
to get any efficient bound for the number of
pseudo-primes in a given unit interval.
Thus, we shall take a detour, by employing the Riesz mean of
$\psi_\varGamma(x)$:
$$
\tilde\psi_\varGamma(x)=\int_1^x\psi_\varGamma(y){dy\over y}.
\eqno(44.3)
$$
It should be noted that the lower limit of
integration is due to the fact that $N(\varpi)>1$ by definition.
We have, for any $\tau>0$,
$$
\tau^{-1}\left(\tilde\psi_\varGamma(x)
-\tilde\psi_\varGamma(xe^{-\tau})\right)
\le \psi_\varGamma(x)\le\tau^{-1} \left(
\tilde\psi_\varGamma(xe^\tau)-
\tilde\psi_\varGamma(x)\right).\eqno(44.4)
$$
On the other hand, we have
$$
\eqalignno{
\tilde\psi_\varGamma(x)&={1\over 2\pi i}
\int_{(2)}\bigg({\zeta'_\varGamma\over
\zeta_\varGamma}(s)-{\zeta'_\varGamma\over
\zeta_\varGamma}(s+1)\bigg)x^s{ds\over s^2}\cr
&={1\over 2\pi i}\int_{(2)}{\zeta'_\varGamma\over
\zeta_\varGamma}(s)x^s{ds\over s^2}+O\big((\log x)^2\big),
&(44.5)
}
$$
where the error term is the result of moving the
contour to $\Re s=(\log x)^{-1}$ in the
relevant part in the first line; note that $(42.9)$.
The integral in the second line is
$$
\eqalignno{
&{1\over 2\pi i}
\int_{2-iK}^{2+iK}{\zeta'_\varGamma\over
\zeta_\varGamma}(s)x^s{ds\over s^2}+O(1)\cr
=&{1\over2\pi i}\int_{2-iK}^{2+iK}
\Psi_\varGamma(s)x^s{ds\over s^2}
+{1\over2\pi i}\int_{2-iK}^{2+iK}
\Psi^1_\varGamma(s)x^s{ds\over s^2}+O(1),&(44.6)
}
$$
where $K\approx x^3$ is chosen so that $(43.9)$ is
fulfilled, and $\Psi^1_\varGamma(s)=(\zeta'_\varGamma/
\zeta_\varGamma)(s)-\Psi_\varGamma(s)$. The
second integral on the right side is, by $(42.3)$, equal to
$x+O\big(x^{1/2}\big)$. In fact, the part composed of
logarithmic derivatives of the $\Gamma$-function
can be estimated by shifting the path to $\Re s={1\over4}$,
say, and the part containing $(\zeta'/\zeta)(2s)$ can be
estimated by shifting the path to $\Re s={1\over2}$.
Hence, we have
$$
\tilde\psi_\varGamma(x)=x+{1\over2\pi i}\int_{2-iK}^{2+iK}
\Psi_\varGamma(s)x^s{ds\over s^2}
+O\big(x^{1/2}\big).\eqno(44.7)
$$
We move the contour to $C_K$ which is the result of
connecting $2-iK$, $-{1\over2}-iK$,
$-{1\over2}+iK$, $2+iK$ by straight segments:
$$
\eqalignno{
\tilde\psi_\varGamma(x)&=x+
{1\over2\pi i}\int_{C_K}
\Psi_\varGamma(s)x^s{ds\over s^2}
+2\Re\mathop{{\sum}^\r{u}}_{|\nu_V|\le K}{x^{1/2+ \nu_V}
\over{\big({1\over2}+ \nu_V}\big)^2}+O\big(x^{1/2}\big)\cr
&=x+{1\over2\pi i}\int_{C_K}
\Psi_\varGamma(s)x^s{ds\over s^2}
+O\big(x^{1/2}\log x\big),&(44.8)
}
$$
where we have applied $(43.7)$. The last integral can be
ignored. In fact, to estimate the part over $\Re s=-{1\over2}$ one
may use the bound $\Psi_\varGamma(s)\ll |s|\log|s|$ which
follows from $(43.7)$--$(43.8)$; and the segment
over $\Im s=\pm K$ is divided into pieces of length $1/K$ and
we use $(43.7)$--$(43.9)$. With this, we return to $(44.4)$, and
set $\tau=x^{-1/4}(\log x)^{1/2}$. We obtain
$$
\psi_\varGamma(x)=x+O\left(x^{3/4}(\log x)^{1/2}\right).
\eqno(44.9)
$$
We end the proof of $(44.1)$.
\medskip
\noindent
{\csc Notes:} The Riemann hypothesis for $\zeta(s)$ implies 
an error term in the prime number theorem that
is equivalent to the validity of the hypothesis. Hence
it appears natural to expect that a similar correspondence
should hold for $\zeta_\varGamma(s)$ and the pseudo-prime
number theorem for the group $\varGamma$. This is, however,
still a challenging problem. The difficulty is in the
fact that there are too many complex zeros
to effectively deal with. See Iwaniec [12] for the first 
significant step.
\vskip 0.7cm
\centerline{\bf REFERENCES}
\bigskip
\item{[1]} M. Baruch and Z. Mao.
Bessel identities
in Waldspurger correspondence, the archime\-dean theory.
Israel J. Math., {\bf145} (2005), 1--82.
\item{[2]} R.W. Bruggeman. Fourier coefficients 
of cusp forms. Invent.\ math., {\bf 45} (1978), 1--18.
\item{[3]} ---. {\it Fourier Coefficients of
Automorphic Forms\/}.  Lect.\ Notes in Math., {\bf 865},
Springer-Verlag, Berlin 1981.
\item{[4]} ---. {\it Families of Automorphic Forms\/}. Birkh\"auser 
Verlag, Berlin 1994.
\item{[5]} R.W. Bruggeman and Y. Motohashi. 
A new approach to the
spectral  theory of the fourth moment of 
the Riemann zeta-function.
J. reine angew.\ Math., {\bf 579} (2005), 75--114.
\item{[6]} J.W. Cogdell and I.I. Piatetskii-Shapiro.  {\it The
Arithmetic and Spectral Analysis of 
Poin\-car\'e Series\/}. Academic
Press, San Diego 1990.
\item{[7]} J. Delsarte. Sur le gitter fuchsien. C.R. Acad.\ Sci.\
Paris, {\bf 214} (1942), 147--149.
\item{[8]} I.M. Gel'fand, M.I. Graev and  I.I.
Pyatetskii-Shapiro. {\it Representation Theory and Automorphic
Functions\/}. W.B. Saunders Company, Philadelphia 1969.
\item{[9]} I.S. Gradshteyn and I.M. Ryzhik. {\it Tables of
Integrals, Series and Products\/}.  
Academic Press, San Diego 1979.
\item{[10]} T. Hawkins. {\it Emergence of the Theory of
Lie Groups. An Essay in the History of
Mathematics 1869--1926\/}. Springer-Verlag, 
New York 2000.
\item{[11]} D.A. Hejhal. {\it The Selberg Trace Formula 
for ${\rm PSL}(2,{\Bbb{R}})$\/}. Vol.\ 1. 
Lect.\ Notes in Math., {\bf548}, Springer-Verlag, 
Berlin 1976; Vol.\ 2. ibid, {\bf 1001}, 1983.
\item{[12]} H. Iwaniec. Prime geodesic theorem. J. reine angew.,
{\bf 349} (1984), 136--159.
\item{[13]} H. Jacquet. Fonctions de 
Whittaker associ\'ees aux groupes
de Chevalley. Bull.\ Soc.\ Math.\ France, {\bf95} (1967),
243--309.
\item{[14]} H. Jacquet and R.P. Langlands. {\it Automorphic
Forms on $\r{GL}(2)$\/}. Lect.\ Notes in Math., {\bf 114},
Springer-Verlag, Berlin 1970.
\item{[15]} M. Jutila and Y. Motohashi. 
Uniform bound for Hecke
$L$-functions. Acta Math., {\bf 195} (2005), 61--115.
\item{[16]}  A.A. Kirillov. On $\infty$-dimensional unitary
representations of the group 
of second-order matrices with elements
from a locally compact field. 
Soviet Math.\ Dokl., {\bf 4} (1963), 748--752.
\item {[17]} N.V.  Kuznetsov. Petersson hypothesis 
for forms of weight
zero and Linnik  hypothesis. Preprint:
Khabarovsk Complex Res.\ Inst., East Siberian 
Branch Acad.\ Sci.\ USSR,
Khabarovsk, 1977. (Russian). This is partly translated into
English: Petersson hypothesis for parabolic forms 
of weight zero and Linnik hypothesis. 
Sums of Kloosterman sums.  Math.\ USSR--Sb., 
{\bf 39} (1981), 299--342.
\item{[18]} N.N. Lebedev. {\it Special Functions \& Their
Applications\/}. Dover Publ.\ Inc., Mineola, New York 1972.
\item{[19]} H. Maass. \"Uber eine neue Art von 
nichtanalytischen automorphen Funktionen
und die Bestimmung Dirichletscher Reihen 
durch Funktionalgleichungen. 
Math.\ Ann., {\bf 121} (1949), 141--183.
\item{[20]} ---. {\it Modular Functions 
of One Complex Variable\/}. Lect.\ Math.\ Phys.,
{\bf 29}, Tata IFR, Bombay 1964/83.
\item{[21]} Y. Motohashi. The binary additive divisor problem.
Ann.\ Sci.\ \'Ecole Norm.\ Sup.\ $4^e$ s\'{e}rie, {\bf 27} 
(1994), 529--572.
\item{[22]} ---. {\it Spectral Theory of the Riemann
Zeta-Function\/}. Cambridge Tracts in Math., {\bf 127},
Cambridge Univ.\ Press, Cambridge 1997.
\item{[23]} ---. A note on the mean value of the zeta 
and $L$-functions.\ Part XII. Proc.\ Japan Acad., 
{\bf 78}A (2002), 36--41; 
Part XIV. ibid, {\bf 80}A (2004), 28--33;
Part XV.  ibid, {\bf 83}A (2007), 73--78. 
\item{[24]} ---. Mean values of zeta-functions via
representation theory. Proc.\ Symp.\ Pure Math., AMS, 
{\bf 75} (2006), 257--279.
\item{[25]} ---. Sums of Kloosterman sums revisited. 
In: {\it The Conference on $L$-Functions\/}. World Scientific, 
Singapore 2007, pp.\ 141--163.
\item{[26]} ---. The Riemann zeta-function and the 
Hecke congruence subgroups. RIMS Kyoto Univ.\ Kokyuroku, 
{\bf 958} (1996), 166--177; Part II. J. Res.\ Inst.\ Sci.\ Tech.\ 
Nihon Univ., {\bf 119} (2009), 29--64
\item{[27]} ---. {\it Analytic Number Theory\/}. Vol.\ 1.
{\it Distribution of Prime Numbers\/}.\
Asakura Books, Tokyo 2009; Vol.\ 2.
{\it Zeta Analysis\/}. ibid 2011. (Japanese). An enlarged
English edition is under preparation.
\item{[28]} R. Penrose. 
{\it The Road to Reality. A Complete Guide to the Laws of
the Universe\/}. Vintage Books, New York 2007.
\item{[29]} H. Petersson. \"Uber eine Metrisierung der ganzen 
Modulformen.  Jber.\ Deutsch.\ Math.\ Verein., {\bf 49}
(1939), 49--75.
\item{[30]} F. Riesz and B. Sz-Nagy. {\it Funtional Analysis\/}.
Frederick Ungar Publ.\ Co., New York 1955.
\item{[31]} W. Roelcke. Das Eigenwertproblem der
automorphen Formen in der hyperbolischen Ebene. 
Teil I. Math.\ Ann., 
{\bf 167} (1966), 292--337; Teil II. ibid.,
{\bf 168} (1967), 261--324.
\item{[32]} A. Selberg. Harmonic analysis and discontinuous 
groups in weakly symmetric Riemannian 
spaces with applications to Dirichlet series. 
J. Indian Math.\ Soc., {\bf20} (1956), 47--87.
\item{[33]} ---. On the estimation of Fourier coefficients
of modular forms.  Proc.\ Symp.\ Pure Math., 
AMS, {\bf 8} (1965), 1--15.
\item{[34]} S.A.  Stepanov. The number of points of a hyperelliptic
curve over a prime field. Izv.\ Akad.\ Nauk SSSR ser.\ Mat.,
{\bf 33} (1969), 1171--1181. (Russian)
\item{[35]} E.C. Titchmarsh. {\it Fourier Integrals\/}.
Clarendon Press, Oxford 1967.
\item{[36]} N.Ja.\ Vilenkin and A.U. Klimyk.
{\it Representation of Lie Groups and Special 
Functions\/}. Vol.\ 1.
Kluwer Acad.\ Publ., Dordrecht 1991; This is an
enlarged edition of  N.Ja.\ Vilenkin: {\it Special Functions and
the Theory of Group Representations\/}.  Amer.\ Math.\ Soc.,
Providence 1968.
\item{[37]} G.N. Watson. {\it A Treatise on the 
Theory of Bessel Functions\/}.  
Cambridge Univ.\ Press, Cambridge 1996.
\item{[38]} E.T. Whittaker and G.N. Watson. 
{\it A Course of Modern Analysis\/}. Cambridge Univ.\ Press, 
Cambridge 1927.
\bigskip
\medskip
\noindent
\font\small=cmr8
{\small\noindent 
Department of Mathematics, Nihon University,
Surugadai, Tokyo 101-8308, Japan
}
\hfill\def\ymzeta
{\font\brm=cmr17 at 30pt\font\ssrm=cmr5 at 4pt
\font\sssrm=cmr5 at2.5pt
{{\brm O}\raise 9pt\hbox{\hskip -22pt
$\hfil\raise3pt\hbox{\sssrm KH}\atop\hbox{
{\ssrm Y}$\zeta$\hskip-1pt{\ssrm M}}$}}}
\ymzeta

\bye